\documentclass[aos,preprint]{imsart}

\usepackage[hang]{subfigure}
\usepackage{color}
\usepackage{latexsym}
\usepackage{mathtools}
\usepackage{amsmath}
\usepackage{amssymb}
\usepackage[utf8]{inputenc}
\usepackage{textcomp}
\usepackage{graphicx}
\usepackage{epsfig}
\usepackage{dsfont}
\usepackage{hyperref}
\usepackage{bbm}
\usepackage{enumitem}
\usepackage{here}

\usepackage{color}
\usepackage[dvipsnames]{xcolor}

\RequirePackage{amsthm,amsmath,amsfonts,amssymb}
\RequirePackage[authoryear,numbers]{natbib}
\RequirePackage[colorlinks,citecolor=blue,urlcolor=blue]{hyperref}
\RequirePackage{graphicx}

\usepackage{multirow}

\startlocaldefs

\newtheorem{lem}{Lemma}
\newtheorem{assumption}{Assumption}

\renewcommand{\epsilon}{\varepsilon}

\newcommand{\argmin}{\mathop{\rm arg~min}\limits}
\providecommand{\skakko}[1]{\left(#1\right)}
\providecommand{\mkakko}[1]{\left\{#1\right\}}
\providecommand{\lkakko}[1]{\left[#1\right]}
\newcommand{\abs}[1]{\left\lvert#1\right\rvert}
\usepackage{bm} 
\theoremstyle{plain}
\newtheorem{theorem}{Theorem}[section]

\theoremstyle{remark}
\newtheorem{remark}{Remark}
\theoremstyle{remark}

\newtheorem{example}{Example}

\endlocaldefs

\begin{document}

\begin{frontmatter}
	\title{
Parametric generalized Spectrum \\for Heavy-Tailed Time Series
}
	\runtitle{Parametric generalized Spectrum for Heavy-Tailed Time Series}
	
\begin{aug} 
\author[A]{\fnms{Yuichi} \snm{Goto}\ead[label=e1]{yuichi.goto@math.kyushu-u.ac.jp}}
\and
\author[B]{\fnms{Gaspard} \snm{Bernard}\ead[label=e2]{gaspardbernard@as.edu.tw}}

\address[A]{\qq{Faculty of Mathematics}, Kyushu University\printead[presep={,\ }]{e1}}
\address[B]{\qq{Institute of Statistical Science}, Academia Sinica\printead[presep={,\ }]{e2}}
\end{aug}

\begin{abstract}
Recently, several spectra have emerged, designed to encapsulate the distributional characteristics of non-Gaussian stationary processes. 
This article introduces parametric families of generalized spectra based on the characteristic function, alongside inference procedures
enabling $\sqrt{n}$-consistent estimation of the unknown parameters in a broad class of parametric models.
These spectra capture non-linear dependencies without requiring that the underlying stochastic processes satisfy any moment assumptions.
Crucially, this approach facilitates frequency domain analysis for heavy-tailed time series, including possibly non-causal Cauchy autoregressive models and discrete-stable integer-valued autoregressive models. To the best of our knowledge, the latter models have not been studied theoretically in the literature.
By estimating parameters across both causal and non-causal parameter spaces, our method automatically identifies the causal or non-causal structure of Cauchy autoregressive models.
Furthermore, our estimator does not depend on smoothing parameters since it is based on the integrated periodogram associated with the generalized spectrum.
As applications, we develop goodness-of-fit tests, moving average unit-root tests, and tests for non-invertibility. 
We study the finite-sample performance of the proposed estimators and tests via Monte Carlo simulations, and apply the methodology to estimation and forecasting of a measles count dataset. {We evaluate finite-sample performance using Monte Carlo simulations and illustrate the practical value of the procedure with an application to measles case-count estimation and forecasting.}
\end{abstract}

\end{frontmatter}

\section{Introduction}
Spectral methods are widely employed in the analysis of time series data and have found applications in various real-world scenarios. 
In the context of Gaussian processes, the classical spectrum, based on second-order moments, completely captures the distributional characteristics of the time series. Consequently, it stands as one of the most successful techniques for analyzing Gaussian time series data. 
However, spectral methods that rely on the classical spectrum, although valid for non-Gaussian processes, suffer from several limitations. 
First, they fail to capture the distributional characteristics of non-Gaussian processes, as they can only detect linear dependencies and not non-linear ones. 
Second, these methods require certain moment assumptions to hold, { which is in practice {often problematic} when dealing with dataset {known for exhibiting} heavy-tailed behaviors such as financial data.}

To address these issues, several new spectra have been introduced.
One notable advancement is the introduction of the {\it generalized spectral density}, based on the characteristic function, as proposed by \cite{h99}. \cite{h99} established the rates of convergence of the kernel estimator of the spectrum with respect to the integrated mean squared error while
\cite{fp17} and \cite{fp18} developed a test for independence based on this spectrum. Tests for invertibility of vector autoregressive {(AR)}  models with non-Gaussian innovations based on the generalized spectrum have been proposed by \cite{cce17}.
The spectra based on the joint distribution function and copula have been explored by \cite{h00}, later named {\it Laplace spectral density kernel} and {\it copula spectral density kernel}, respectively, by \cite{dhkv13}. 
\cite{h00} and \cite{lr11}  proposed goodness-of-fit tests for the {Laplace spectral density}, whereas
\cite{dhkv13} and \cite{kvdh16} studied the problem of estimating {these} spectra. 
\cite{l08}, \cite{l12}, \cite{h13}, and \cite{l21} { studied the Laplace spectrum and quantile spectrum in relation to the copula spectral density kernel.}
Additionally, \cite{vvd18} proposed spectra based on Spearman's $\rho$ and Kendall's $\tau$. 
For a comprehensive review on spectral methods, refer to \cite{v20}.

Nevertheless, the majority of literature has primarily focused on the development of nonparametric methods, except for \cite{cce17} and \cite{v22}, leaving the question of explicit closed-form expressions for parametric generalized spectra largely unexplored.
Notably, \cite{v22} studied the parameter estimation of a non-causal and non-invertible parametric infinite-order moving average (MA($\infty$)) model with i.i.d.\ or martingale difference errors based on the generalized spectrum of {residuals}.
An important difference with our approach is that \cite{v22} does not { leverage the closed} form of generalized spectra and { its study does} not encompass moving average (MA) models with unit roots since the reciprocal of the linear filter of {an} infinite-order MA model is used for the calculation of residual from observations. {Another crucial difference between our results and \cite{v22} is that, although the existence of the generalized spectrum itself does not require moment assumptions, the consistency results for the parametric estimators proposed in \cite{v22} still rely on rather stringent finite-moment conditions.
In contrast, our procedure accommodates heavy-tailed innovations.}

Our primary aim in this article is to derive explicit forms for the generalized spectra of processes for which the classical definition of the spectra cannot be used.
This includes non-causal and/or non-invertible processes as well as processes exhibiting infinite variance. 
Subsequently, we tackle the question of estimating the unknown parameters in some models { where the closed-form expressions for the} parametric generalized spectra allows it. Our parameter estimation strategy relies on minimizing a least squares criterion { quantifying the discrepancy} between the parametric model and the periodogram of the generalized spectrum.
It is worth mentioning that our estimation method does not require the use of any smoothing parameter for kernel estimates, a distinct advantage.

Traditionally, the periodogram, while widely used, lacks consistency as an estimator of the spectral density. Consequently, the smoothed periodogram, incorporating  some bandwidth parameter, is typically employed to achieve consistency. However, the choice for the bandwidth parameters significantly impacts the estimates. 
To circumvent these difficulties, it is common to consider {\it functionals} of the spectra instead.
This approach is developed in, for example, Section 7.6 of \cite{b81}, \cite{d85}, and Section 6 of \cite{tk00}. 
Notably, the integrated copula spectral density kernel has been recently studied by \cite{gkvvdh22}. Furthermore, \cite{c88} and \cite{pd13} discussed  inference based on some least squares difference for the classical spectral density and the time-varying spectral density, respectively.

Our methodology enables the application of spectral methods to heavy-tailed time series models, encompassing (continuous) MA and AR models with error processes following stable distributions and non-negative integer-valued MA {(INMA)} and non-negative integer-valued {AR} {(INAR)} models  with error processes following discrete stable distributions. 

Constructing heavy-tailed non-negative integer-valued time series presents a unique challenge. {In practice, count-valued datasets tend to exhibit a few distinctive features. In particular, they are often both over-dispersed (i.e., the variance substantially exceeds the mean) and zero-inflated (i.e., zeros are overrepresented in the dataset).}
Nevertheless, there exists an extensive literature on non-negative integer-valued time series, including popular models such as the integer-valued autoregressive (INAR) model \cite{m85, aa91, l97} and the integer-valued generalized autoregressive conditional heteroscedasticity (INGARCH) model \cite{flo06}.
Comprehensive reviews of recent developments can be found in \cite{dfhllp21}.

INGARCH models require at least that the conditional mean exists, since the model is defined through the autoregressive structure of the conditional expectation. 
\cite{g20} explored an INGARCH model based on the beta-negative binomial distribution, which accommodates heavy tails; however, finite variance is still assumed to guarantee the finiteness of the asymptotic variance of the estimator.
Conversely, the INAR model seems to be able to accommodate heavy tails by assuming heavy-tailed disturbance processes, but the lack of a zero-median property due to its non-negative integer-valued nature makes the {direct} application of the least absolute deviation method problematic. { On this note, a major difficulty in applying LAD to the INAR model stems from the stepwise nature of the conditional median in the integer-valued setting. In particular, the conditional median of $Z_t \mid Z_{t-1}=z$ is not linear in $z$, which makes centering to enforce a zero-median property more delicate. Note that this issue is specific to settings that combine integer-valued models with the absence of a finite first moment.}
Therefore, to our best knowledge, non-negative integer-valued time series models with heavy-tails have not been developed so far. { One of the aims of this contribution is to propose methods that can be readily applied in such scenarios.}

{As a compelling application of our method, we address the analysis of heavy-tailed integer-valued data. Figure \ref{fig:ts_hill} displays the time series and Hill plots for the measles dataset obtained from \cite{tscount}. The empirical evidence strongly suggests the presence of heavy tails, with even the mean being undefined, as hinted by the Hill estimator plot displayed in Figure \ref{fig:ts_hill}.
\begin{figure}[htbp]
\centering
\includegraphics[width=0.8\textwidth]{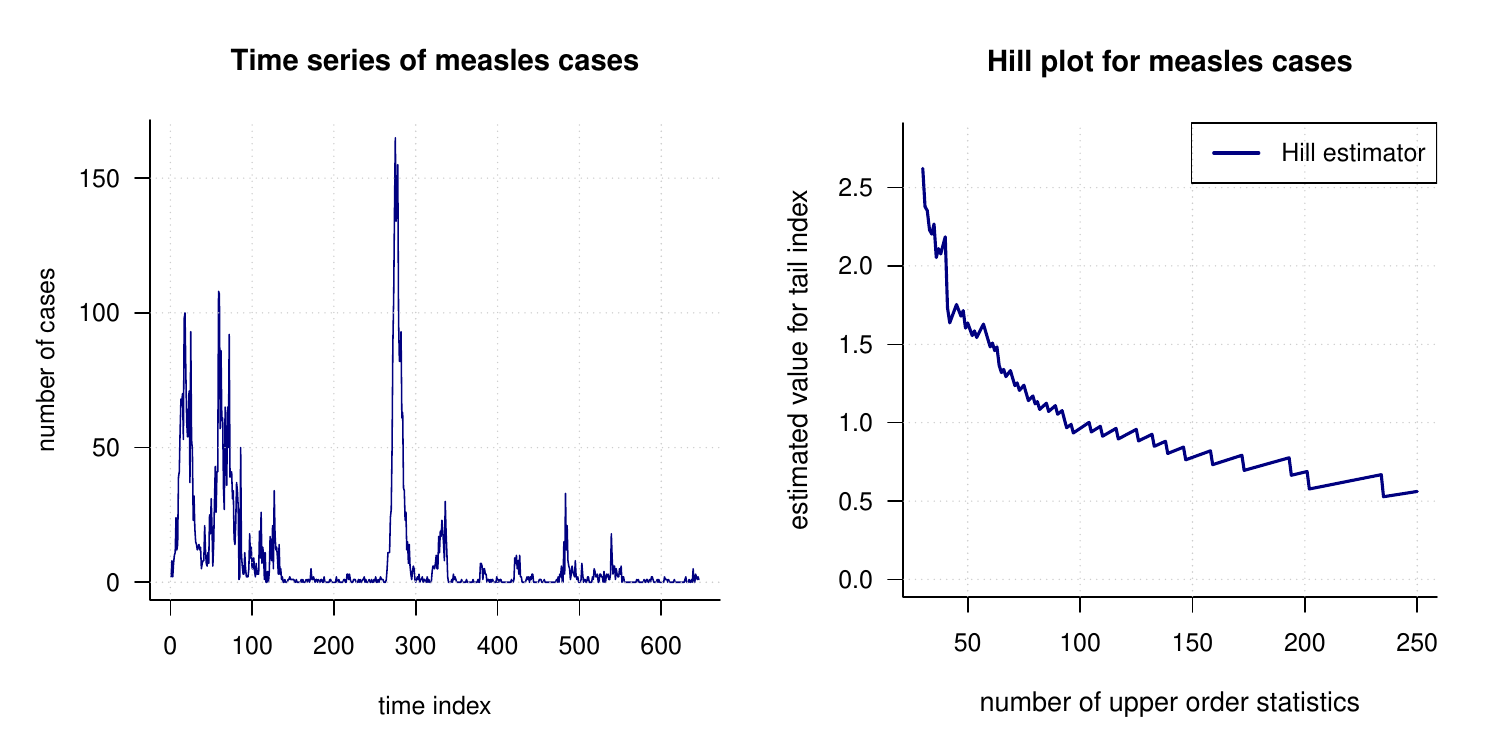}
\caption{Time series plot (left) and Hill plot (right) of the measles dataset.}
\label{fig:ts_hill}
\end{figure}

 It should also be noted that generalized spectrum is also linked to the notion of {\it distance correlations}, but we do not discuss this fact in detail here. Further details can be found in \cite{dmmw18, efp19, wd22}, and references therein.

It is noteworthy that non-causal and non-invertible autoregressive moving average (ARMA) processes possess a causal invertible representation with the same (classical) spectral density  as the original process \cite[p.125--127]{bd09}. 
{In particular}, for Gaussian ARMA processes, the causality and invertibility is required to ensure identifiability. 
For non-Gaussian ARMA processes, for $k\geq3$, $k$-th order spectrum { is identifiable} for non-causal and non-invertible ARMA processes \cite[Theorem 1]{vl18}.
However, the $k$-th order spectrum {requires} at least finite $k$-th order moments { to exist}, { which is an assumption that cannot be satisfied under heavy tails}.
{\cite{hv25} identified non-causal and non-invertible ARMA models by first estimating their causal and invertible representation, and then determining the optimal parameters that minimize a chosen criterion across all possible representations.
Unlike this approach, our proposed method directly estimates parameters across the causal and non-causal parameter  space. Consequently, the resulting estimates automatically identify the causal or non-causal structure of the model.}


{We first derive closed-form expressions for the generalized spectrum of  processes in several relevant parametric heavy-tailed time series models, including the integer-valued models mentioned earlier. We then study parameter identifiability and propose root-$n$ consistent estimators. Unlike conventional spectral approaches, our assumptions and proofs extend to generalized spectral densities that may fail to be H\"{o}lder continuous---a feature of discrete stable processes---thereby ensuring the validity of our methodology in heavy-tailed settings where standard regularity conditions break down.}
 
The organization of the rest of the paper is as follows. 
Section \ref{sec2} reviews the generalized spectrum and its fundamental properties.
Section \ref{sec3} derives explicit forms of the generalized spectrum for several models and discusses the identifiability of parametric families.
Section \ref{sec4} presents estimation methods for parametric families of the generalized spectra and investigates their statistical properties.
Section \ref{sec5} is devoted to goodness-of-fit tests for given parametric families.
Section \ref{sec6} develops unit-root and non-invertibility tests for MA models.
Section \ref{sec7} demonstrates the finite-sample performance of the proposed methods.
Section \ref{sec8} illustrates the utility of our proposed methods via real data analysis.
{Supplemental material} provides  proofs of all results and additional figures and simulation results.

\section{Generalized spectrum}\label{sec2}

{ We start by defining the generalized spectrum and stating a few of its general properties for a strictly stationary process $\{Z_t\}$.}
The generalized spectrum (of order two) is defined, for $\lambda\in\mathbb R$ and $(u,v)\in\mathbb R^2$, as
\begin{align}\label{gs}
f(\lambda;u,v)\coloneqq \frac{1}{2\pi}\sum_{\ell=-\infty}^\infty{\rm Cov}\skakko{e^{i uZ_{t+\ell}},e^{-ivZ_{t}}}\exp\skakko{-i\ell\lambda},
\end{align}
where $i\coloneqq \sqrt{-1}$, under the summability condition
$\sum_{\ell=-\infty}^\infty\left|{\rm Cov}\skakko{e^{i uZ_{t+\ell}},e^{-ivZ_{t}}}\right|<\infty$ for all $(u,v)\in\mathbb R^2$ holds. Note that we adopt the convention that $${\rm Cov}\skakko{e^{i uZ_{t+\ell}},e^{-ivZ_{t}}}\coloneqq {\rm E}\skakko{e^{iuZ_{t+\ell}+ivZ_{t}}}-
{\rm E}\skakko{e^{iuZ_{t+\ell}}}{\rm E}\skakko{e^{ivZ_{t}}}.$$
Since \eqref{gs} is based on characteristic functions, {the construction does not require that $\{Z_t\}$ satisfy any moment assumption}. { Moreover, unlike the classical spectrum based on the autocovariance function, the generalized spectrum can capture 
{nonlinear} serial dependence.}
{ It can be readily checked that, from its definition,} the generalized spectrum satisfies the following properties, for any $(u,v)\in\mathbb R^2$ and $\lambda\in \mathbb R$,
\begin{align*}
f(\lambda;u,v)=&f(\lambda+2\pi;u,v),\quad
f(\lambda;u,v)=f(-\lambda;v,u),\\
\overline{f(\lambda;u,v)}=&f(-\lambda;-u,-v),\quad
f(\lambda;u,-u)\in\mathbb R.
\end{align*}

It is worth mentioning that the pairwise time reversibility of $Z_t$---meaning that the joint distributions of $(Z_t, Z_{t+\ell})$ and $(Z_{t+\ell}, Z_t)$ are identical for any $\ell \in \mathbb{Z}$---is equivalent to the spectral symmetry condition $f(\lambda; u, v) = f(\lambda; v, u)$ for all $(u,v)\in\mathbb R^2$ and $\lambda\in \mathbb R$.

{Note that} if $\{Z_t\}$ is an i.i.d.\ process, \eqref{gs} reduces to
$$
f_{\rm i.i.d.}(\lambda;u,v)\coloneqq {\rm Cov}\skakko{e^{i uZ_{t}},e^{-ivZ_{t}}}/(2\pi)
$$
and, thus, the generalized spectrum for any i.i.d.\ process is a constant function with respect to $\lambda\in[-\pi,\pi]$ for fixed $(u,v)\in\mathbb R^2$.

{We now introduce the} generalized spectrum of order $k$, defined for $(\lambda_1,\ldots,\lambda_{k-1})\in[0,2\pi]^{k-1}$ and $(u_1,\ldots,u_k)\in\mathbb R^k$ as
\begin{align*}
&f(\lambda_1,\ldots,\lambda_{k-1};u_1,\dots,u_k)\\
\coloneqq &
\frac{1}{(2\pi)^{k-1}}\sum_{\ell_1,\dots,\ell_{k-1}=-\infty}^\infty
{\rm Cum}\skakko{e^{iu_1Z_{t+\ell_1}},\cdots,e^{iu_{k-1}Z_{t+\ell_{k-1}}},e^{iu_kZ_{t}}}
\exp\skakko{-i\sum_{j=1}^{k-1}\ell_j\lambda_j},
\end{align*}
where ${\rm Cum}(Z_1,\cdots,Z_k)$ is the cumulant of order $\ell$ defined as
$${\rm Cum}(Z_1,\cdots,Z_k)\coloneqq 
\sum_{(\nu_1,\ldots,\nu_p)}(-1)^{p-1}(p-1)!\;{\rm E}\skakko{\prod_{j\in\nu_1}Z_{\nu_1}}\ldots{\rm E}\skakko{\prod_{j\in\nu_p}Z_{\nu_p}}.
$$
The summation $\sum_{(\nu_1,\ldots,\nu_p)}$ extends over all partitions $(\nu_1,\ldots,\nu_p)$ of $\{1,2,\cdots,k\}$ (see \citealt[p.19]{b81}). { We assume that $\{Z_t\}$ satisfies the following necessary summability} condition
$$\sum_{\ell_1,\dots,\ell_{k-1}=-\infty}^\infty\left|
{\rm Cum}\skakko{e^{iu_1Z_{t+\ell_1}},\cdots,e^{iu_{k-1}Z_{t+\ell_{k-1}}},e^{iu_kZ_{t}}}\right|<\infty\quad\text{for all $(u_1,\ldots,u_k)\in\mathbb R^k$}.$$ 

It is worth noting that Lemma 2.1. of \cite{gzkc23} implies that the geometrically $\alpha$-mixing strictly stationary process satisfies, for any $q\in\mathbb N$ and any $k\in\mathbb N$,
\begin{align}\nonumber
\sum_{\ell_1,\dots,\ell_{k-1}=-\infty}^\infty&
\skakko{1+\sum_{j={1}}^{k-1}\left| { \ell_j}\right|^q}\\\label{as_sum}
&\times
\sup_{(u_1,\ldots,u_k)\in\mathbb R^k}
\left|
{\rm Cum}\skakko{e^{iu_1Z_{t+\ell_1}},\cdots,e^{iu_{k-1}Z_{t+\ell_{k-1}}},e^{iu_kZ_{t}}}\right|<\infty.
\end{align}

{ When considering various time series models in Section \ref{sec3} where stationarity and geometric $\alpha$-mixing properties have not been established in the literature, we will typically provide a similar property, ensuring the existence of generalized spectra.}
\section{Parametric generalized spectrum}\label{sec3}
For specific stochastic processes associated with classical time series models, the explicit form of the generalized spectrum is of particular interest. 
In this section, we present several important examples of such parametric spectra.

\subsection{Continuous-valued models}
We first {consider various natural} continuous-valued time series models.
Throughout this section and the subsequent analysis, we denote by $B$ the \emph{backshift operator}. We start by deriving a closed form for the generalized spectra in the Cauchy MA(1) and AR(1) models.

\begin{example} 
The possibly non-invertible MA(1) model with Cauchy innovations and coefficient $a \in \mathbb{R}_0\coloneqq \mathbb{R}\setminus \{0\}$ is defined as 
$$Z_t=(1-a^{-1}B)\epsilon_t,$$ where $\{\epsilon_t\}$ follows an i.i.d.\ Cauchy distribution with scale parameter $\delta>0$.
The generalized spectrum of $\{Z_t\}$ is given by
\begin{align*}
&f_{\bm \theta}(\lambda;u,v)\\
=&
-\frac{1}{2\pi}
\exp\skakko{-\delta(|u|+|v|)(|a^{-1}|+1)}
(2\cos\lambda+1)+\frac{1}{2\pi}
\exp\skakko{-\delta|u+v|(|a^{-1}|+1)}\\
&+\frac{1}{2\pi}
\exp\skakko{-\delta(|va^{-1}|+|v+a^{-1}u|+|u|)}
\exp\skakko{-i\lambda}\\
&+\frac{1}{2\pi}
\exp\skakko{-\delta(|ua^{-1}|+|u+a^{-1}v|+|v|)}
\exp\skakko{i\lambda}.
\end{align*}

\end{example}

{It is crucial to establish identifiability of the model parameters from this parametric family of generalized spectra. This property is what allows us to perform consistent parameter estimation.}

\begin{lem}\label{IdentifiableCMA}
{In the Cauchy MA(1) model}, the vector of parameters ${\bm \theta}\coloneqq ( a, \delta)^\top$ is identifiable { from the generalized spectra}.
\end{lem}


\begin{example}
The possibly non-causal AR(1) model with Cauchy innovations and coefficient $a \in \mathbb R\setminus\{-1,1\}$ is defined as
$$Z_t=a Z_{t-1} + \epsilon_t,$$
where $\{\epsilon_t\}$ follows an i.i.d.\ centered Cauchy distribution with scale parameter $\delta>0$. 
The generalized spectrum of this model is given by
$$f_{\bm \theta}(\lambda, u,v)=\frac{1}{2\pi} \sum_{\ell=-\infty}^{\infty} {C}_{\ell,a,\delta}(u,v)\;e^{-i\ell\lambda},$$
where, for $a$ such that $|a|<1$,
\begin{align}\label{eq:causalAR}
C_{\ell,a,\delta}(u,v)
\coloneqq 
\begin{cases}
\exp\skakko{-\frac{\delta}{1-|a|}
\skakko{
|ua^{\ell}+v|
+
|u|(1-|a|^\ell)
}}
-
\exp\skakko{-\frac{\delta(|u|+|v|)}{1-|a|}} & \ell\geq0\\
C_{|\ell|,a,\delta}(v,u)
 & \ell<0.\\
\end{cases}
\end{align}
and 
for $a$ such that $|a|>1$,
\begin{align}\nonumber
&C_{\ell,a,\delta}(u,v)\\\label{eq:noncausalAR}
\coloneqq &
\begin{cases}
\exp\skakko{-\frac{\delta|a^{-1}|}{1-|a^{-1}|}
\skakko{
|va^{-\ell}+u|
+
|v|(1-|a^{-\ell}|)
}}
-
\exp\skakko{-\frac{\delta(|u|+|v|) |a^{-1}|}{1-|a^{-1}|}}
& \ell\geq0\\
C_{|\ell|,a,\delta}(v,u) 
& \ell<0.\\
\end{cases}
\end{align}
\end{example}

\begin{lem}\label{IdentifiableCAR}
In the Cauchy AR(1) model, the vector of parameters ${\bm \theta}\coloneqq ( a, \delta)^\top$ is identifiable { from the generalized spectra}.
\end{lem}

Note that, in $C_{\ell,a,\delta}(u,v)$ for $\ell>0$,
$u$ and $v$ correspond to  the future $Z_{t+\ell}$
and the current $Z_t$, respectively.
The term $|ua^{\ell}+v|
$ in \eqref{eq:causalAR} illustrates that the coefficient $a^{\ell}$ affects the variable $u$ associated with the future, reflecting the causal direction of dependence.
In contrast, the term $|va^{-\ell}+u|
$ in \eqref{eq:noncausalAR} indicates  that 
the coefficient $a^{-\ell}$ affects the variable $v$ associated with the current observation, reflecting the non-causal direction of dependence.

\begin{lem}\label{CAR_mixing}
The Cauchy AR(1) model enjoys geometrically $\alpha$-mixing property.
\end{lem}

Next, we derive the generalized spectra for Gaussian MA(1) and AR(1) models. { We consider these models due to their extremely classical nature.} 


\begin{example} The MA(1) model with Gaussian innovations defined as
$$Z_t=(1-a^{-1}B)\epsilon_t,$$ where $\{\epsilon_t\}$ follows an i.i.d.\ centered normal distribution with variance $\sigma^2$. The generalized spectral density is given by
\begin{align*}
f_{\bm \theta}(\lambda;u,v)
=&
\frac{1}{2\pi}
\exp\skakko{-\frac{\sigma^2}{2}(u^2+v^2)(a^{-2}+1)}\\
&\times\lkakko{
\exp\mkakko{-\sigma^2 uv (a^{-2}+1)}-1
+
2\mkakko{\exp\skakko{\sigma^2 uv a^{-1}}-1}
\cos \lambda},
\end{align*}
where $\bm \theta\coloneqq (a,\sigma^2)^\top$ such that $|a|>1$ and $\sigma^2>0$.
\end{example}
\begin{lem}\label{IdentifiableGMA}
In the Gaussian MA(1) model, the vector of parameters ${\bm \theta}\coloneqq ( a, \sigma^2)^\top$ is identifiable { from the generalized spectra}.
\end{lem}

\begin{example}
The AR(1) model with Gaussian  innovations and coefficient $a$ such that $|a|<1$ is defined as
$$Z_t=a Z_{t-1} + \epsilon_t,$$
where $\{\epsilon_t\}$ follows an i.i.d.\ centered normal distribution with variance $\sigma^2$. 
The generalized spectrum of this model is given by
$$f_{{\bm \theta}}(\lambda, u,v)=\frac{1}{2\pi} \sum_{\ell=-\infty}^{\infty} {C}_{\ell,a,\sigma^2}(u,v)\;e^{-i\ell\lambda},$$
where
\begin{align*}
C_{\ell,a,\sigma^2}(u,v)
\coloneqq 
e^{-\frac{\sigma^2}{2(1-a^2)}(u^2+v^2)} \skakko{ e^{-\frac{\sigma^2}{1-a^2} uv a^{|\ell|}} - 1 }.
\end{align*}
\end{example}

\begin{lem}\label{IdentifiableGAR}
In the Gaussian AR(1) model, the vector of parameters ${\bm \theta}\coloneqq ( a, \sigma^2)^\top$ is identifiable { from the generalized spectra}.
\end{lem}

An important observation is that the spectra of Gaussian MA and AR models satisfy the symmetry condition $f(\lambda,u,v)=f(\lambda,v,u)$, reflecting the time-reversibility of the process.
This implies that the roles of the future variable (associated with $u$) and the current variable (associated with $v$) can be interchanged without altering the spectrum.

\subsection{Discrete-valued models}

Next, we consider heavy-tailed integer-valued time series models.
When constructing a discrete time series model with infinite variance, various discrete-value distributions with infinite variance are available, such as beta-negative binomial and zeta distributions. However, the characteristic functions of these distributions are often quite complex. One discrete-valued distribution with infinite variance that possesses a simple characteristic function is the discrete stable distributions introduced by \cite{sv79}.

The random variable $W$ is said to follow a discrete stable distribution with scale parameter $\delta>0$ and exponent $\alpha\in(0,1]$ if the probability generating function (pgf) is given by, for $x\in\mathcal C$ such that $|x|\leq 1$, $ {\rm E}\skakko{x^{W}}=\exp\{-\delta(1-x)^\alpha\}$. Thus, the characteristic function of the discrete stable distribution is given by $ {\rm E}\skakko{e^{isW}}=\exp\{-\delta(1-e^{is})^\alpha\}$. If $\alpha=1$, $W$ follows a Poisson distribution.
The tail behavior is characterized by \citet[Theorem 1]{cs98} as ${\rm P}(W \ge n) = O(n^{-\alpha})$, which implies that the fractional moment $\mathbb{E}\skakko{W^m}$ is finite if and only if $m < \alpha$.

{ The parameters of the discrete-stable innovations admit straightforward interpretations for count-valued data. The scale parameter $\delta$ directly affects the probability that an innovation equals zero and is therefore closely linked to the degree of zero inflation in the series. By contrast, the exponent $\alpha$ governs tail heaviness and is thus directly related to the over-dispersion of the process. Discrete-stable innovations therefore provide a flexible mechanism for modeling the distinctive features of count-valued time series. We consider INMA and INAR models defined through the use of the {\it binomial thinning operator}, denoted ``$\circ$'' and defined as 
\begin{align*}
p\circ Z_t=
\begin{cases}
\sum_{i=1}^{Z_t} Y_{i,t}& \text{if } Z_t >0\\
0& \text{if } Z_t=0,
\end{cases}
\end{align*}
 for some stochastic process $Z_t$ and with $\{Y_{i,t}\}$ i.i.d.\ (w.r.t. $i$ and $t$) Bernoulli random variables with parameter $p$. The binomial thinning operators aims at generalizing the notion of regression coefficient to integer-valued models, mimicking their effect on the expectation, since ${\rm E}(p\circ Z_t)=p\,{\rm E}(Z_t)$ (granted that ${\rm E}(Z_t)$ is finite).}

\begin{figure}[htbp]
\centering
\includegraphics[width=0.8\textwidth]{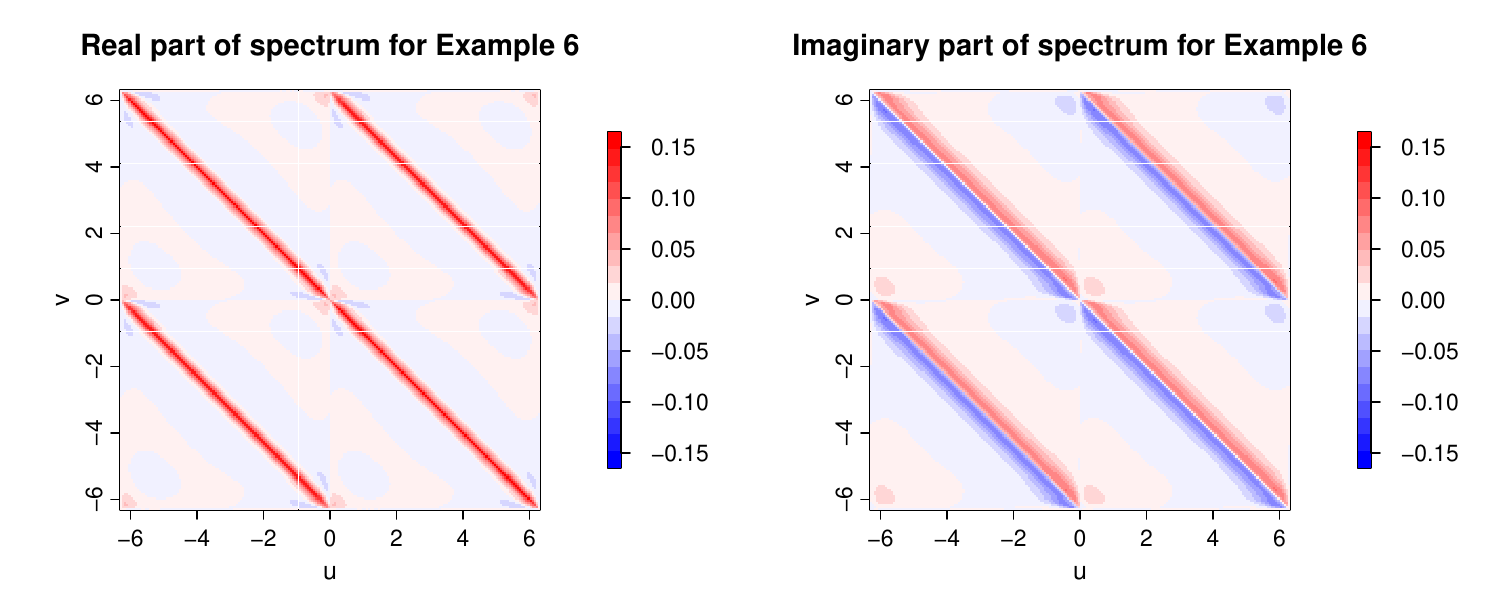}
  \caption{Heatmaps of the real part (left panel) and the imaginary part (right panel) of the spectrum for the causal integer-valued AR(1) model in Example~\ref{IntegerAR(1)}, with parameters $\delta = 2$, $p = 0.3$, $\alpha = 0.7$, and $\lambda = 0.785$.}
  \label{Fig_spec_AR}
\end{figure}

\begin{example}\label{IntegerMA(1)}
The integer-valued MA(1) model with discrete stable innovations is given by
$$Z_t=p\circ \epsilon_{t-1} + \epsilon_t,$$
where ``$\circ$'' { is the binomial thinning operator defined earlier}. $\{\epsilon_t\}$ follows i.i.d.\ discrete stable distribution with parameters $\alpha\in(0,1]$  and {\it scale parameter} $\delta>0$, and
$\{Y_{i,k}\}$ is independent of $\{\epsilon_t\}$.
Note that the characteristic function of Bernoulli distribution is given by ${\rm E}\big(e^{isY_{i,k}}\big)=1-p+pe^{is}$. The generalized spectrum of $\{ Z_t\}$ can be expressed as 
\begin{align*}
f_{\bm \theta}&(\lambda;u,v)\\
\coloneqq &
-\frac{1}{2\pi}
\exp\skakko{-\delta
\lkakko{(p^\alpha+1)
\mkakko{(1-e^{iu})^\alpha
+(1-e^{iv})^\alpha}}}(2\cos\lambda+1)\\
&+\frac{1}{2\pi}
\exp\skakko{-\delta(p^\alpha+1)(1-e^{i(u+v)})^\alpha}\\
&+\frac{1}{2\pi}
\exp\skakko{-\delta \lkakko{p^\alpha(1-e^{iv})^\alpha
+\mkakko{1- e^{iv}(1-p+pe^{iu})}^\alpha
+ (1-e^{iu})^\alpha}}e^{-i\lambda}\\
&+\frac{1}{2\pi}
\exp\skakko{-\delta \lkakko{p^\alpha(1-e^{iu})^\alpha
+\mkakko{1- e^{iu}(1-p+pe^{iv})}^\alpha
+ (1-e^{iv})^\alpha}}e^{i\lambda}.
\end{align*}

\end{example}

{ The following lemma establishes the crucial identifiability property, allowing to conduct parameter estimation based on the generalized spectrum.}

\begin{lem}\label{IdentifiableINMA}
In the causal integer-valued MA(1) model, the vector of parameters ${\pmb \theta}\coloneqq (\delta, \alpha, p)$ is identifiable { from the generalized spectra}.
\end{lem}

\begin{example}\label{IntegerAR(1)}
    The integer-valued AR(1) with discrete stable innovations is given by $$Z_t= p\circ Z_{t-1} + \epsilon_t.$$
    The binomial thinning operator $\circ$ and the discrete stable distribution are defined as in the previous example. The generalized spectrum is 
    $$f_{\bm \theta}(\lambda, u,v)=\frac{1}{2\pi} \sum_{\ell=-\infty}^{\infty} {C}_{\ell,p,\delta,\alpha}(u,v)\;e^{-i\ell\lambda},$$
    with
\begin{align*}
{C}_{\ell,p,\delta,\alpha}(u,v)
\coloneqq &\exp\!\lkakko{
-\frac{\delta}{1- p^{\alpha}}
\mkakko{(1-e^{iu})^{\alpha}\, (1-p^{\alpha \ell})
+
{\skakko{1-e^{{i}v}\left(1-p^{\ell}+p^{\ell}e^{iu}\right)}^{\alpha}}
}}\\
&-
\exp\lkakko{- \frac{\delta}{1-p^\alpha}\mkakko{(1-e^{iu})^\alpha + (1-e^{{i}v})^\alpha}}
\quad 
\text{for $\ell\geq0$}\\
{C}_{\ell,p,\delta,\alpha}(u,v)
\coloneqq &
{C}_{|\ell|,p,\delta,\alpha}(v,u) 
\quad 
\text{for $\ell<0$}
.
\end{align*}
{Here, we used the explicit form of the joint pgf given by:
\begin{align}\label{eq:pgf_INAR1}
{\rm E}\skakko{u^{X_t}}
=&\exp\mkakko{- \delta\frac{(1-u)^\alpha}{1-p^\alpha}}
\end{align}
\text{and }
\begin{align}\label{eq:pgf_INAR2}
{\rm E}\skakko{u^{X_{t+\ell}}v^{X_{t}}}
=&
\exp\mkakko{-\delta \frac{(1-v(1-p^\ell(1-u)))^\alpha}{1-p^\alpha}
-\delta (1-u)^\alpha\frac{1-p^{\alpha\ell}}{1-p^\alpha}}\quad \text{for $\ell\geq0$}.
\end{align}
}
\end{example}

{
Figure \ref{Fig_spec_AR} shows the plot of this spectrum.
The proofs of \eqref{eq:pgf_INAR1} and \eqref{eq:pgf_INAR2} are deferred to Subsection \ref{Proof_pgf_INAR}.}
 {The next lemma guarantees stationarity and the existence of the generalized spectrum.}
 
\begin{lem}\label{lem:stationarity}
If $p\in(0,1)$, the integer-valued AR(1) model admits a unique strictly stationary solution, which is ergodic.
{Moreover, the integer-valued AR(1) model enjoys geometrically $\beta$-mixing property.}
\end{lem}

\begin{remark}\label{rem:INAR}
From the proof of Lemma \ref{lem:stationarity},
it is interesting to note that the integer-valued AR(1) process $Z_t$
also follows a discrete stable distribution with exponent
$\alpha \in (0,1]$ and scale parameter $\delta/(1-p^\alpha)$.
\end{remark}

{ Even if such result can be found in the literature, we chose to include here a short and explicit proof for the sake of clarity. Note also at this point that it has been shown--see for instance Corollary 5 in \cite{szu24} or remark 3 in \cite{bkp13}--that the process is in fact {\it geometrically ergodic}, which by Theorem 3.7 of \cite{brad05} directly implies geometrically $\beta$-mixing property.
}

 The following lemma provides the usual crucial identification result.
\begin{lem}\label{IdentifiableINAR}
In the integer-valued AR(1) model, the vector of parameters ${\pmb \theta}\coloneqq (\delta, \alpha, p)$ is identifiable { from the generalized spectra}.
\end{lem}

\section{Estimation of unknown parameter}\label{sec4}
In the previous section, we introduced the parametric families of generalized spectra. This section is dedicated to the estimation of the unknown parameters of parametric models, based on the generalized spectrum.
The parametric family $\mathcal F$ is defined, for some positive integer $d$ and parameter space $\bm{\Theta}\subset \mathbb R^d$,  as
$$\mathcal F\coloneqq \mkakko{f_{\bm{\theta}}(\lambda;u,v);\bm{\theta}\in\bm{\Theta}}.$$
We assume the following about the parametric family of generalized spectra and about the parameter space.
\begin{assumption}\label{as}\;
\begin{enumerate}
\item[$(a)$] The parameter space $\bm{\Theta}$ is a compact subset of $\mathbb R^d$.

\item[$(b)$] The parametric generalized spectrum $f_{\bm{\theta}}(\lambda;u,v)\in\mathcal F$ is identifiable, that is, for given constant $L>0$ , $f_{\bm{\theta}}(\lambda;u,v)=f_{\bm{\theta^\prime}}(\lambda;u,v)$ for all $\lambda\in[0,2\pi]$ and $(u,v)\in\mathbb [-L,L]^2$ implies $\bm{\theta}=\bm{\theta^\prime}$.



\item[$(c)$]
Let $\mathcal S\subset[-L,L]^2$ be a  closed set such that its $\varepsilon$-neighborhood $\mathcal S_\varepsilon\coloneqq \{(u,v)\in[-L,L]^2;{\rm dist}((u,v),\mathcal S)\leq \varepsilon\}$ satisfies 
${\rm Leb}(\mathcal S_\varepsilon)=O(\varepsilon)$ as $\varepsilon\downarrow0$, where ${\rm Leb}(\mathcal S_\varepsilon)$ is the Lebesgue measure of $\mathcal S_\varepsilon$
and ${\rm dist}((u,v),\mathcal S)\coloneqq \inf_{(x,y)\in \mathcal S}\|(u,v)-(x,y)\|_2$.

Moreover, there exist constants
$\alpha,\beta,\gamma\in(0,1]$, $\tau\geq 0$ and $C^\prime>0$ such that,
for every $\varepsilon>0$,
\begin{align}\label{as:ub}
|f_{\bm\theta}(\lambda;u,v)|\leq C^\prime \quad
\text{for all $(\lambda, u, v, \bm{\theta})\in[0,2\pi]\times[-L,L]^2\times\bm{\Theta}$}
\end{align}
and, for all
$(\lambda,u,v),(\lambda',u',v')\in[0,2\pi]\times([-L,L]^2\setminus\mathcal S_\varepsilon)$
and all $\bm\theta,\bm\theta'\in\bm\Theta$,
\begin{align}\label{as:ph}
|f_{\bm\theta}(\lambda;u,v)-f_{\bm\theta'}(\lambda';u',v')|
\le
C^\prime\,\varepsilon^{-\tau}
\Bigl(
|\lambda-\lambda'|^{\alpha}
+\|(u,v)-(u',v')\|^{\beta}
+\|\bm\theta-\bm\theta'\|^{\gamma}
\Bigr).
\end{align}
We say that $f_{\bm\theta}(\lambda;u,v)$ is
\emph{Uniformly Bounded and Partially H\"{o}lder continuous except for $\mathcal S$},
abbreviated as \textrm{UBPH($\mathcal S,\tau,\alpha,\beta,\gamma$)},
if the above conditions hold.
\end{enumerate}
\end{assumption}

{The assumptions (a) and (b) are extremely common.
{ Assumption (c) is necessary to approximate an integral by its corresponding Riemann sum at the proper rate. Usually, to obtain the proper rate of convergence,
H\"{o}lder continuity is typically assumed.}
However, for the generalized spectra in integer-valued MA and AR models,
the derivatives with respect to $u$ or $v$  behave like
$u^{\alpha-1}$, $v^{\alpha-1}$, or $(u+v)^{\alpha-1}$ for some exponent
$\alpha\in(0,1]$.
This behavior is not compatible with the assumption of uniform H\"{o}lder continuity, { as the derivatives of the generalized spectra are unbounded around} the points
$u=0$, $v=0$, and $u+v=0$ (mod $2\pi$).
To address this issue, we decompose the domain into two parts:
(i) a \emph{good region}, where H\"{o}lder continuity holds uniformly, and
(ii) a \emph{bad region}, defined as an $\varepsilon$-neighborhood of the
singular set, whose Lebesgue measure is sufficiently small { ($\mathcal{S}_{\varepsilon}$ in Assumption 1. (c))}.
This decomposition allows us to control the contribution of the singularities
while retaining quantitative bounds on the approximation error. { In Assumption \ref{as}, the exponents $\alpha,\, \beta,\, \gamma$ quantify the complexity of the inference problem associated with the global regularity away from the singularity, while the exponent $\tau$ quantifies the complexity associated with the blow-up rate locally around the singularity.}
The generalized spectra for the integer-valued MA(1) and AR(1) models satisfy the assumption (c).
}

{ The proposed estimation procedure is based on minimizing (with respect to ${\bm\theta}$) a suitable criterion. As our criterion, we consider the following $L_2$-type divergence function:}
\begin{align*}
D(f,f_{\bm{\theta}})
\coloneqq \int_{0}^{2\pi}\int_{-L}^L\int_{-L}^L|f(\lambda;u,v)-f_{\bm{\theta}}(\lambda;u,v))|^2{\rm d}u{\rm d}v{\rm d}{\lambda}.
\end{align*} 

\begin{remark}
The choice of integrating over $[-L,L]^2$ with respect to $(u,v)$ can seem a little odd at first sight but is quite natural given that the integer-valued examples presented in Section \ref{sec3} are periodic. As the focus of this article is integer-valued models, we then adopt this definition for the sake of convenience. { In all the examples considered in this article, we will chose $L=\pi$.}
\end{remark}
\begin{remark}
There are instances where spectra can be zero at some $\lambda,u,v$, rendering impossible the use of the contrast-type divergences based on $f/f_{\bm{\theta}}$, including the Whittle-type divergence function. For instance, 
{
$f(\lambda,u,0)=f(\lambda,0,v)=0$ for all $(\lambda,u,v)$ for the integer-valued MA(1) and AR(1) models and the MA(1) and AR(1) models with centered Cauchy errors.}
 In this sense, our method is significantly more robust than the Whittle-likelihood methods. These methods rely on the assumption that the spectral density is bounded away from zero, which is often imposed in classical spectral analysis.
\end{remark}

To approximate $D(f,f_{\bm{\theta}})$, we consider 
\begin{align}\label{D_n}
D_n( I_n,f_{\bm{\theta}})
\coloneqq \frac{8\pi L^2}{nM_n^2}\sum_{j=1}^{n-1}\sum_{i_1,i_2=1}^{M_n}| I_n(\lambda_j;u_{i_1},v_{i_2})-f_{\bm{\theta}}(\lambda_j;u_{i_1},v_{i_2}))|^2,
\end{align}
where $\lambda_j=2\pi j/n$, $u_{i_1}\coloneqq -L+{2Li_1}/{M_n}, v_{i_2}\coloneqq -L+{2Li_2}/{M_n}$, { $M_n \in \mathbb{N}_0$ some parameter to be specified in practice} and, for any $\lambda\in[0,2\pi]$ and $(u,v)\in\mathbb R^2$,
$$I_n(\lambda;u,v)\coloneqq \frac{1}{2\pi n}d_n(\lambda;u)d_n(-\lambda;v), \quad
d_n(\lambda;u)\coloneqq\sum_{t=1}^ne^{iuZ_t}e^{-it\lambda}.
$$
The next lemma shows that $D_n( I_n,f_{\bm{\theta}})$ has a non-negligible bias but first, we need to introduce some extra notation. Consider sets $\mathcal{S}$ and $\mathcal{S}_{\epsilon}$ as in Assumption \ref{as} ($c^\prime$). We say that the generalized spectrum of order 2, $f(\lambda, u,v)$ is \emph{Partially H\"{o}lder continuous except for $\mathcal{S}$} (${\rm PH}(\mathcal{S},\tilde{\tau},\tilde{\alpha},\tilde{\beta})$) with exponent $\tilde{\alpha}$ associated to $\lambda \in [0, 2\pi)$ and with exponent $\tilde{\beta}$ associated to $(u,v)\in[-L,L]^2$ if there exists $C'>0$ such that
\begin{align*}
|f(\lambda;u,v)-f(\lambda';u',v')|
\le
C^\prime\,\varepsilon^{-\tilde{\tau}}
\Bigl(
|\lambda-\lambda'|^{\tilde{\alpha}}
+\|(u,v)-(u',v')\|^{\tilde{\beta}}
\Bigr).
\end{align*}

\begin{lem}\label{lem1}
 Consider a { family of parametric generalized spectra} satisfying Assumption~\ref{as} and a stationary process ${Z_t}$ satisfying condition \eqref{as_sum} for $q=1$. Assume that the generalized spectrum of order two of $\{ Z_t\}, f(\lambda, u, v)$ is ${\rm PH}(\mathcal{S},\tilde{\tau},\tilde{\alpha},\tilde{\beta})$. Assume moreover that, as $n\to\infty$, $M_n \rightarrow \infty$.
Then, for any $\bm{\theta}\in\bm{\Theta}$, 
$D_n( I_n,f_{\bm{\theta}})$ converges in probability to $\tilde D(f,f_{\bm{\theta}})$ as $n\to\infty$
where
\begin{align*}
\tilde D(f,f_{\bm{\theta}})
\coloneqq 
D(f,f_{\bm{\theta}})+\int_{0}^{2\pi}\int_{-L}^L\int_{-L}^L
f(\lambda;u,-u)f(-\lambda;v,-v){\rm d}u{\rm d}v{\rm d}{\lambda}.
\end{align*}
\end{lem}

As the bias $\int_{0}^{2\pi}\int_{-L}^L\int_{-L}^L
f(\lambda;u,-u)f(-\lambda;v,-v){\rm d}u{\rm d}v{\rm d}{\lambda}$ is independent of $\bm{\theta}$, $D_n( I_n,f_{\bm{\theta}})$ can be used for the estimation of the unknown parameters $\bm{\theta}$. We define the true value ${\bm{\theta}}_0$ and its estimator $\hat{\bm{\theta}}_n$ by
\begin{align}\label{Estimator}
{\bm{\theta}}_0\coloneqq  
\argmin_{\bm{\theta}\in\bm{\Theta}}
 D(f,f_{\bm{\theta}})
\skakko{=\argmin_{\bm{\theta}\in\bm{\Theta}}\tilde D(f,f_{\bm{\theta}})} \text{ and }
\hat{\bm{\theta}}_n\coloneqq 
\argmin_{\bm{\theta}\in\bm{\Theta}}
D_n\skakko{ I_n,f_{\bm{\theta}}}.
\end{align}

\begin{remark}
Our estimator is closely related to \cite{ky02} who explored parameter estimation based on the joint characteristic function ${\rm E}\big(e^{i\sum_{j=0}^pu_jZ_{t+j}}\big)$ for a strictly stationary and ergodic process $\{Z_t\}$ and a fixed constant $p$.
Their estimator can capture the joint distribution of $(Z_t,\ldots,Z_{t+p})$. In contrast, our estimator is based on the bivariate characteristic functions of $(Z_{t+\ell},Z_t)$ for all $\ell\in\mathbb Z$.
\end{remark}

\begin{remark}
\cite{km94} investigated the self-normalized classical periodogram for MA($\infty$) process, when the disturbance belongs to the domain of attraction of a stable law. Subsequently,  \cite{mgkt95} estimated parameters of ARMA model whose error process belongs to the domain of attraction of a stable law through the use of the self-normalized periodogram.
\end{remark}

We assume from now on that $D(f,f_{\bm{\theta}})$ has a unique minimum at ${\bm{\theta}}_0\in\mathbb \bm{\Theta}$. This condition is fulfilled {if the model is correctly specified or, in other words,} $f\in\mathcal F$, {but it can also be satisfied under certain forms of misspecification}. Then, the following theorem shows that $\hat{\bm{\theta}}_n$ is a consistent estimator of ${\bm{\theta}}_0$.

\begin{theorem}\label{thm1}
 Consider a { family of generalized spectra} satisfying Assumption \ref{as} and a stationary process ${Z_t}$ satisfying condition \eqref{as_sum} for $q=1$. Assume that the generalized spectrum of order two of $\{ Z_t\}, f(\lambda, u, v)$ is ${\rm PH}(\mathcal{S},\tilde{\tau},\tilde{\alpha},\tilde{\beta})$. Assume moreover that, as $n\to\infty$, $M_n \rightarrow \infty$ and that $D(f,f_{\bm{\theta}})$ possesses a unique minimum. The estimator $\hat{\bm{\theta}}_n$ converges in probability to ${\bm{\theta}}_0$ as $n\to\infty$.
\end{theorem}
In the proof of Theorem \ref{thm1}, {we show that the stochastic equicontinuity of $D_n( I_n,f_{\bm{\theta}})-\tilde D(f,f_{\bm{\theta}})$ and the pointwise consistency of Lemma \ref{lem1} is sufficient to yield some uniform laws of large numbers of the form
$\sup_{\bm{\theta}\in\bm{\Theta}}
\left|D_n( I_n,f_{\bm{\theta}})-\tilde D(f,f_{\bm{\theta}})\right|
\stackrel{\rm P}{\rightarrow} 0$ as $n\to\infty$.

{ 
As usual, obtaining the asymptotic normality requires some more stringent regularity assumptions. We show that under the assumption that $f_{\bm{\theta}}$ is twice continuously differentiable with respect to ${\bm{\theta}}$ with derivatives of {UBPH}, $\hat{\bm{\theta}}_n$ is asymptotically normal with the usual parametric convergence rates}.

\begin{assumption}\label{as2}\;
\begin{enumerate}


\item[($d$)]
The parametric spectrum $f_{\bm{\theta}}(\lambda;u,v)$ is twice differentiable with respect to $\bm{\theta}$.
Its derivatives
{
$\frac{\partial^2}{\partial \bm{\theta}\partial \bm{\theta}^\top}
f_{\bm\theta}(\lambda;u,v)$ and
$
\frac{\partial}{\partial \bm{\theta}}
f_{\bm\theta}(\lambda;u,v)
$ are UBPH($\mathcal S,\tau,\alpha,\beta,\gamma$), where the absolute value on the left-hand sides of
\eqref{as:ub} and \eqref{as:ph}
is replaced by the norm $\|\cdot\|$.
}

\item[($e$)]
The true value ${\bm{\theta}}_0$ belongs to the interior of ${\bm{\Theta}}$.

\item[($f$)]
The matrix ${\bm J}$ defined in \eqref{J} is invertible.
\end{enumerate}
\end{assumption}

{It can be easily checked that} the generalized spectra for the INMA(1) and INAR(1) models satisfy Assumption 2 ($d$).

\begin{remark}
 The differentiability condition of $f_{\bm{\theta}}$ with respect to $\bm{\theta}$ 
excludes MA and AR processes with Cauchy errors studied in Section \ref{sec3}. 
A possible way to handle non-differentiable parametric models is to replace the 
absolute value function $|x|$ with a smooth approximation. For example, one may 
consider the sigmoid-based function
$$
\phi(x)\coloneqq 2\log(1+e^{kx})/k - x - 2\log(2)/k,$$
where $k>0$. 
Its derivative
$
\phi'(x)={2e^{kx}}/{(1+e^{kx})}-1
$
provides a smooth approximation of the sign function.
{Taking $k$ arbitrarily large, we expect the value of $D_n$ defined in \eqref{D_n} to be nearly identical with and without smoothing.} In the simulation study, we do not implement such a smoothing
modification and simply use the original absolute value function. Nevertheless, the empirical results appear to support the asymptotic
normality of the estimator.
\end{remark}

We are now ready to state the asymptotic normality result for the proposed estimator $\hat{\bm{\theta}}_n$.

\begin{theorem}\label{thm2}
 Consider a { family of parametric generalized spectra} satisfying Assumptions~\ref{as} and~\ref{as2} and a stationary process ${Z_t}$ satisfying condition \eqref{as_sum} for $q=1$. Assume that the generalized spectrum of order four of $\{ Z_t\}$ is {\it partially H\"{o}lder} (${\rm PH}$ as in Theorem \ref{thm1}) with arbitrary exponents (note that it implies that second order generalized spectrum $f(\lambda, u, v)$ is also partially H\"{o}lder). Assume that $f(\lambda, u, v)$ is ${\rm PH}(\tilde{\tau}, \tilde{\alpha}, \tilde{\beta})$ with exponents satisfying
 $$\frac{\min\{\alpha, \tilde{\alpha}, \beta, \tilde{\beta}\}}{\max\{\tau, \tilde{\tau}\}+1}> 1/2.$$
 Assume finally that $n\to\infty$, $n^{-1/2}M_n \rightarrow \infty$ and that $D(f,f_{\bm{\theta}})$ admits a unique minimizer. 
 
 Then, under these assumptions, 
 $$\sqrt n(\hat{\bm{\theta}}_n-{\bm{\theta}}_0)\stackrel{\cal L}{\rightarrow} \mathcal{N}(0,\bm J^{-1}\bm I \bm J^{-1})$$ 
 as $n\to\infty$. 
Here, ${\bm J}$ and $\bm I$ are defined as follows:
\begin{align}\nonumber
{\bm J}\coloneqq &
-\int_{0}^{2\pi}\int_{-L}^L\int_{-L}^L
f(\lambda;u,v)\frac{\partial^2}{\partial \bm{\theta}\partial \bm{\theta}^\top}\overline{f_{\bm{\theta}}(\lambda;u,v)}\Big|_{\bm{\theta}={{\bm{\theta}}_0}}
{\rm d}u{\rm d}v{\rm d}{\lambda}\\\nonumber
&-\int_{0}^{2\pi}\int_{-L}^L\int_{-L}^L
\frac{\partial^2}{\partial \bm{\theta}\partial \bm{\theta}^\top}f_{\bm{\theta}}(\lambda;u,v)\Big|_{\bm{\theta}={{\bm{\theta}}_0}}\overline{f(\lambda;u,v)}
{\rm d}u{\rm d}v{\rm d}{\lambda}\\\label{J}
&+\int_{0}^{2\pi}\int_{-L}^L\int_{-L}^L
\frac{\partial^2}{\partial \bm{\theta}\partial \bm{\theta}^\top}
\skakko{f_{\bm{\theta}}(\lambda;u,v)\overline{f_{\bm{\theta}}(\lambda;u,v)}}\Big|_{\bm{\theta}={{\bm{\theta}}_0}}
{\rm d}u{\rm d}v{\rm d}{\lambda}.
\end{align}
and
\begin{align*}
\bm I\coloneqq &(I_{st})_{s,t=1,\ldots,d}, \quad I_{st}\coloneqq L_{1_{st}}+L_{2_{st}}+L_{3_{st}}+L_{4_{st}}.
\end{align*}
Denoting by $\delta(x)$ the Dirac delta function centered in $x$ and with the convention $\textbf{d}{\bm w}={\rm d}u{\rm d}u^\prime{\rm d}v{\rm d}v^\prime{\rm d}\lambda{\rm d}\lambda^\prime$, $L_{1_{st}},\, L_{2_{st}},\, L_{3_{st}},\,L_{4_{st}}$ are defined as follows:
\begin{align}\nonumber
L_{1_{st}}\coloneqq 
&{2\pi}
\int_{(0,2\pi)^2}\int_{[-L,L]^4}
\Big[
\delta(\lambda-\lambda^\prime)
f(\lambda;u,v^\prime)
f(-\lambda;v,u^\prime)
\\\nonumber
&+
\delta(2\pi-\lambda-\lambda^\prime)
f(\lambda;u,u^\prime)
f(-\lambda;v,v^\prime)
+
f(\lambda,-\lambda,\lambda^\prime;u,v,u^\prime,v^\prime)\Big]\\\nonumber
&\times
\frac{\partial}{\partial {\theta_s}}\overline{f_{\bm{\theta}}(\lambda;u,v)}\Big|_{\bm{\theta}={{\bm{\theta}}_0}}\frac{\partial}{\partial {\theta_t}}\overline{f_{\bm{\theta}}(\lambda^\prime;u^\prime,v^\prime)}\Big|_{\bm{\theta}={{\bm{\theta}}_0}}
\textbf{d}{\bm w},
\\\nonumber
L_{2_{st}}\coloneqq 
&2\pi
\int_{(0,2\pi)^2}\int_{[-L,L]^4}
\Big[
\delta(\lambda-\lambda^\prime)
f(\lambda;u,-u^\prime)
f(-\lambda;v,-v^\prime)\\\nonumber
&+
\delta(2\pi-\lambda-\lambda^\prime)
f(\lambda;u,-v^\prime)
f(-\lambda;v,-u^\prime)
+
f(\lambda,-\lambda,-\lambda^\prime;u,v,-u^\prime,-v^\prime)\Big]\\\nonumber
&\times
\frac{\partial}{\partial {\theta_s}}\overline{f_{\bm{\theta}}(\lambda;u,v)}\Big|_{\bm{\theta}={{\bm{\theta}}_0}}\frac{\partial}{\partial {\theta_t}}{f_{\bm{\theta}}(\lambda^\prime;u^\prime,v^\prime)}\Big|_{\bm{\theta}={{\bm{\theta}}_0}}
\textbf{d}{\bm w},
\\\nonumber
L_{3_{st}}\coloneqq 
&2\pi
\int_{(0,2\pi)^2}\int_{[-L,L]^4}
\Big[
\delta(\lambda-\lambda^\prime)
f(\lambda;-v,v^\prime)
f(-\lambda;-u,u^\prime)\\\nonumber
&+
\delta(2\pi-\lambda-\lambda^\prime)
f(\lambda;-v,u^\prime)
f(-\lambda;-u,v^\prime)
+
f(-\lambda,\lambda,\lambda^\prime;-u,-v,u^\prime,v^\prime)\Big]\\\nonumber
&\times
\frac{\partial}{\partial {\theta_s}}{f_{\bm{\theta}}(\lambda;u,v)}\Big|_{\bm{\theta}={{\bm{\theta}}_0}}\frac{\partial}{\partial {\theta_t}}\overline{f_{\bm{\theta}}(\lambda^\prime;u^\prime,v^\prime)}\Big|_{\bm{\theta}={{\bm{\theta}}_0}}
\textbf{d}{\bm w},
\\\nonumber
L_{4_{st}}\coloneqq 
&2\pi
\int_{(0,2\pi)^2}\int_{[-L,L]^4}
\Big[
\delta(\lambda-\lambda^\prime)
f(\lambda;-v,-u^\prime)
f(-\lambda;-u,-v^\prime)
\\\nonumber
&
+
\delta(2\pi-\lambda-\lambda^\prime)
f(\lambda;-v,-v^\prime)
f(-\lambda;-u,-u^\prime)
\\\nonumber
&+
f(-\lambda,\lambda,-\lambda^\prime;-u,-v,-u^\prime,-v^\prime)\Big]\\\nonumber
&
\times
\frac{\partial}{\partial {\theta_s}}{f_{\bm{\theta}}(\lambda;u,v)}\Big|_{\bm{\theta}={{\bm{\theta}}_0}}\frac{\partial}{\partial {\theta_t}}{f_{\bm{\theta}}(\lambda^\prime;u^\prime,v^\prime)}\Big|_{\bm{\theta}={{\bm{\theta}}_0}}
\textbf{d}{\bm w}.
\end{align}
\end{theorem}
{
\begin{remark}
        In Theorem \ref{thm2}, note that the integration over frequencies in $L_{1_{st}},\, L_{2_{st}},\, L_{3_{st}},\,L_{4_{st}}$ is taken on the open domain $(0, 2\pi)^2$, which implies that the Dirac delta functions $\delta(\cdot)$ do not take values on the boundary of the domain.
\end{remark}

\begin{remark}
  In Theorem \ref{thm2}, the regularity condition $$\frac{\min\{\alpha, \tilde{\alpha}, \beta, \tilde{\beta}\}}{\max\{\tau, \tilde{\tau}\}+1}> 1/2$$
  has a natural interpretation. If the generalized spectrum is smoother away from the singularity (larger $\alpha,\, \tilde{\alpha},\, \beta,\, \tilde{\beta}$), then one can tolerate a faster blow-up of the derivative near the singularity (larger $\tau,\, \tilde{\tau}$), and conversely. In other words, the condition quantifies a trade-off between global smoothness and local behavior around the singularity in our inference problem.
If Hölder continuity holds globally (i.e., $\tau=\tilde \tau=0$), the condition reduces to $\min\{\alpha, \tilde{\alpha}, \beta, \tilde{\beta}\} > 1/2$. Furthermore, if the related generalized spectra are Lipschitz continuous (corresponding to $\alpha = \tilde{\alpha} = \beta = \tilde{\beta} = 1$), this condition is trivially satisfied.

  \end{remark}
 }

\section{Goodness-of-fit test}\label{sec5}
In statistical analysis, it is often of interest to assess whether a fitted model is appropriate. 
In this section, we consider a goodness-of-fit test. 
More precisely, for a given parametric family 
$\mathcal{F}=\{ f_{\bm{\theta}}(\lambda;u,v) : \bm{\theta}\in\bm{\Theta} \}$, 
we test the composite hypothesis
\[
H_1: f \in \mathcal{F} 
\quad \text{versus} \quad 
K_1: f \notin \mathcal{F}.
\]

We observed in the previous section that the estimator $D_n(I_n,f_{\bm{\theta}})$ of the distance $D(f,f_{\bm{\theta}})$ 
suffers from the bias term
\[
\int_{0}^{2\pi}\int_{-L}^L\int_{-L}^L
f(\lambda;u,-u)\, f(\lambda;-v,v)\,{\rm d}u\,{\rm d}v\,{\rm d}\lambda,
\]
which is difficult to estimate by nonparametric methods. 
One might consider using the integrated version of 
$I_n(\lambda_j;u,-u) I_n(\lambda_j;-v,v)$ as an estimator of the bias term. 
However, Lemma~\ref{lem_cum} shows that
\begin{align*}
&{\rm E}\left[ I_n(\lambda_j;u,-u) I_n(\lambda_j;-v,v) \right]\\
=&
f(\lambda_j;u,-u) f(\lambda_j;-v,v)
+
f(\lambda_j;u,v) f(\lambda_j;-v,-u)
+
O\left({1}/{n}\right),
\end{align*}
which indicates that 
$I_n(\lambda_j;u,-u) I_n(\lambda_j;-v,v)$ contains a non-negligible bias. 
To address this problem, we add an adjustment term $A_n$, where $A_n$ is non-negative and the bias arising $D_n + A_n$ is estimable. The non-negative restriction prevents potential power loss under some alternatives.
Thus, we propose the following test statistic.
\begin{align*}
T_{n}= D_n(I_n, f_{\hat{\bm \theta}_n})+ A_n({I}_n,f_{\hat{\bm \theta}_n})-B_n(I_n),
\end{align*}
where
\begin{align*}
A_n( I_n,f_{\bm{\theta}})
\coloneqq &
\frac{4\pi L^2}{nM_n^2}\sum_{j=1}^{n-1}\sum_{i_1,i_2=1}^{M_n}
\left| I_n(\lambda_j;u_{i_1},-u_{i_1})-f_{\bm{\theta}}(\lambda_j;u_{i_1},-u_{i_1})
\right.\\
&\quad\qquad\qquad\qquad\left.
+I_n(\lambda_j;{-}v_{i_2},v_{i_2})-f_{\bm{\theta}}(\lambda_j;{-}v_{i_2},v_{i_2})\right|^2,
\end{align*}
and
\begin{align*}
B_n(I_n)
\coloneqq &
\frac{8\pi L^2}{nM_n^2} \sum_{j=1}^{n-1}\sum_{i_1,i_2=1}^{M_n} \bigg(
\frac{{I}^2_n(\lambda_j;u_{i_1},-u_{i_1})+I^2_n(\lambda_j;{-}v_{i_2},v_{i_2})}{4} 
\\
&\quad\qquad+ 
\frac{I_n(\lambda_j;u_{i_1},v_{i_2})\bar{I}_n(\lambda_j;u_{i_1},v_{i_2})+ I_n(\lambda_j;u_{i_1},-u_{i_1}){I}_n(\lambda_j;{-}v_{i_2},v_{i_2})}{2}
\bigg).
\end{align*}
The term $B_n$ corresponds to the bias correction term.
We design $T_{n}$ as the estimator of the quantity 
$D(f,f_{\bm{\theta}_0}) + A(f,f_{\bm{\theta}_0})$,
where
\begin{align*}
A(f,f_{\bm{\theta}})\coloneqq 
\frac{1}{2}
\int_{0}^{2\pi}\int_{-L}^L\int_{-L}^L
&| f(\lambda_j;u,-u)-f_{\bm{\theta}}(\lambda_j;u,-u)\\
\qquad\qquad
&+f(\lambda_j;{-}v,v)-f_{\bm{\theta}}(\lambda_j;{-}v,v)|^2
{\rm d}u\,{\rm d}v\,{\rm d}\lambda.
\end{align*}
The adjustment term $A_n$ induces the following biases:
$$
\int_{0}^{2\pi}\int_{-L}^L\int_{-L}^L
f(\lambda;u,v) f(\lambda;-v,-u) 
+ 
\frac{1}{2}
f^2(\lambda;u,-u)
+
\frac{1}{2}
f^2(\lambda;-v,v)
{\rm d}u\,{\rm d}v\,{\rm d}\lambda.
$$
Among these, the second and third bias terms
are directly estimable.
On the other hand, the biases $$
\int_{0}^{2\pi}\int_{-L}^L\int_{-L}^L
f(\lambda;u,-u) f(\lambda;-v,v)
{\rm d}u\,{\rm d}v\,{\rm d}\lambda
$$
and
$$
\int_{0}^{2\pi}\int_{-L}^L\int_{-L}^L
f(\lambda;u,v) f(\lambda;-v,-u)
{\rm d}u\,{\rm d}v\,{\rm d}\lambda$$
arising from $D_n(I_n, f_{\hat{\bm \theta}_n})$ and  $A_n( I_n,f_{\hat{\bm \theta}_n})$, respectively, cannot be estimated individually but can be estimated jointly.

The hypothesis can now be rephrased as follows:
\[
\tilde H_1: D_0(f, f_{{\bm \theta}_0}) + A(f,f_{{\bm \theta}_0})=0
\quad \text{versus} \quad 
\tilde K_1: D_0(f, f_{{\bm \theta}_0}) + A(f,f_{{\bm \theta}_0})>0.
\]
From the proof of Theorem \ref{thm2}, we see that 
\begin{align*}
D_n( I_n,f_{\hat{\bm{\theta}}_n})
=&
D_n( I_n,f_{{\bm{\theta}}_0})
+
\frac{\partial}{\partial \bm{\theta}}
D_n(I_n,f_{{\bm{\theta}}})\Big|_{\bm{\theta}={{\bm{\theta}}_0}}
(\hat{\bm{\theta}}_n-\bm{\theta}_0)\\
&+
\frac{1}{2}
\skakko{\hat{\bm{\theta}}_n-{\bm{\theta}}_0}^\top
\frac{\partial^2}{\partial \bm{\theta}\partial \bm{\theta}^\top}
D_n(I_n,f_{{\bm{\theta}}})\Big|_{\bm{\theta}={{\bm{\theta}}_n^*}}
\skakko{\hat{\bm{\theta}}_n-{\bm{\theta}}_0}\\
=&
D_n( I_n,f_{{\bm{\theta}}_0})
+
O_p\skakko{\frac{1}{n}}.
\end{align*}
Thus, the plug-in effect is negligible for $\sqrt n D_n( I_n,f_{\hat{\bm{\theta}}_n})$.
Assume that $A(f,f_{\bm{\theta}})$ has a unique minimum at ${\bm{\theta}}_0\in \mathring{{\bm{\Theta}}}$ to ensure ${\partial}/(\partial \bm{\theta}) A(f,f_{\bm{\theta}})\mid_{\bm{\theta}=\bm{\theta}_0}=0$.
Then, by the same argument, the plug-in effect is negligible for $\sqrt n A_n(I_n,f_{\hat{\bm{\theta}}_n})$.
Therefore, we obtain
\begin{align*}
T_{n}= D_n(I_n, f_{{\bm \theta}_0})+ A_n({I}_n,f_{{\bm \theta}_0})-B_n(I_n)+
O_p\skakko{\frac{1}{n}}.
\end{align*}
Then, the cumulant argument gives that 
$\sqrt n \skakko{T_{n}- D_0(f, f_{{\bm \theta}_0}) - A(f,f_{{\bm \theta}_0})}$ converges in distribution to a normal distribution as $n\to\infty$.
However, the asymptotic variance is very complex and difficult to obtain a precise estimate. Therefore, we employ a sub-sampling method; that is, we reject $H_1$ whenever $p_{n,1}\leq\varphi$,
where 
\begin{align*}
p_{n,1}\coloneqq \frac{1}{n-b+1}\sum_{t=1}^{n-b+1}\mathbb I
\left\{
\sqrt b T_{b,t}
>
\sqrt n T_{n}
\right\}
\end{align*}
and $T_{b,t}$ is defined as $T_{n}$ but based on the subsample $\{Z_t,\ldots,Z_{t+b-1}\}$.

\section{Hypothesis tests for unit-root and invertibility}\label{sec6}
In this section, we propose one-sided and two-sided hypotheses tests based on the subsampling method. 
Applications of the theory developed in this section include testing for the existence of MA unit roots and for the non-invertibility for non-causal  MA($d_1+d_2$) models with Cauchy residuals. {Here, we make the assumption that the parameters of the model can be properly identified by the parametric generalized spectrum.  Closed forms of generalized spectrum for Cauchy MA(1) models have been provided in Section \ref{sec3}, as well as an identifiability result.}

Theorem~\ref{thm2} gives the asymptotic distribution of the proposed estimator.  However, for constructing the asymptotically valid tests (and confidence intervals), the estimation of the asymptotic variance is required, which turns out to be a non-trivial problem.
Such a problem has been studied in \cite{tpk96}, where variance estimation involved some kernel-based consistent estimator of 
$$\int_{-\pi}^\pi\int_{-\pi}^\pi \phi(\lambda,\lambda^\prime)f_{\rm L_2}(-\lambda,\lambda^\prime,-\lambda^\prime){\rm d}\lambda{\rm d}\lambda^\prime,$$ where $\phi(\cdot,\cdot)$ is a continuous function and $f_{\rm L_2}(\cdot,\cdot,\cdot)$ is the classical fourth order spectrum (see Proposition 1 of \citealt{tpk96}). As mentioned in the Introduction, selecting an appropriate bandwidth parameter for the smoothed periodogram poses a challenge. 
Moreover, tests and confidence intervals based on the estimated asymptotic variance often exhibit unsatisfactory finite-sample performance. 
To address this issue, resampling methods have been investigated by several authors.
Recently, \cite{mpk20} introduced the hybrid periodogram bootstrap for spectral means and ratio statistics based on the classical spectrum under mild conditions, highlighting the failure of the multiplicative periodogram bootstrap to capture terms involving fourth-order spectra in variance. \cite{kp23} applied this method to Whittle likelihood estimation. Although the assumptions are mild, they { still require for} the classical spectrum to be non-zero, which is not necessarily satisfied by generalized spectra. 
In particular, we confirmed earlier that the generalized spectral can be zero, making this approach unsuitable for our specific problem.
In contrast, subsampling methods do not require such a condition and are applicable to a broad range of  $\alpha $-mixing processes. See \cite{prw99} for details. Note also that  \cite{gkvvdh22} constructed confidence bands and hypothesis tests based on the subsampling method for the integrated copula spectral density kernel.
In the following section, we introduce subsampling-based hypothesis tests.

\subsection{Two-sided hypothesis}
We start by considering the classical two-sided hypothesis testing problem {formulated as}
\begin{align*}
H_2:{{\theta}}_{0i}=\kappa_i
\quad \text{versus} \quad 
K_2:{{\theta}}_{0i}\neq\kappa_i,
\end{align*}
where $\bm{\theta}_0=({\theta}_{01},\ldots,{\theta}_{0d})^\top$. Then, the subsampling-based $p$-value can be defined as
\begin{align*}
p_{n,2}\coloneqq \frac{1}{n-b+1}\sum_{t=1}^{n-b+1}\mathbb I
\left\{
\left|
\sqrt b \skakko{{\hat{\theta}}_{b,t}^{(i)}-\kappa}\right|
>
\left|\sqrt n \skakko{{\hat{\theta}}_{n}^{(i)}-\kappa}\right|
\right\},
\end{align*}
where ${\hat{\theta}}_{b,t}^{(i)}$ is the subsample estimator based on $\{Z_t,\ldots,Z_{t+b-1}\}$. 
{The following result establishes the asymptotic validity and consistency of this test.}

\begin{theorem}\label{thm3}
{Assume that the assumptions of Theorem \ref{thm2} hold. Consider the strictly stationary process $\{Z_t\}$ is $\alpha$-mixing with vanishing $\alpha$-mixing coefficient, and that the sub-sampling block length $b$ satisfies $b/n\to0$ as $n\to\infty$.} 
The test which rejects $H_2$ in favor of $K_2$ whenever $p_{n,2}\leq\varphi$ has  asymptotic {size} $\varphi$ and is {consistent}.
\end{theorem}

We demonstrate how the above can be applied to construct the unit-root test.
\begin{example}\label{ex:unitroot}
Consider MA($d_1+d_2$) models with Cauchy residuals for any nonnegative integers $d_1$ and $d_2$, that $\{ Z_t\}$ satisfies
$$Z_t=
(1-\xi_1^{-1}B)\cdots(1-\xi_{d_1}^{-1}B)
(1-\xi_{d_1+1}B^{-1})\cdots(1-\xi_{d_1+d_2}B^{-1})\epsilon_t,$$
where $\epsilon_{t}$ follows i.i.d.\ centered Cauchy distribution with scale parameter $\delta>0$.
The parameter space $\bm\Theta$ for $(\xi_1,\ldots,\xi_{d_1+d_2},\delta)$ is set to be $[L^{-1},L]^{d_1+d_2+1}$. 
{We moreover assume that Assumptions \ref{as} (b) and \ref{as2} (f) hold}. {Testing for the presence of an MA unit root in the model comes to considering the following null hypotheses and alternatives.} 
\begin{align*}
H_{2,{\rm ex}}:\xi_i=1
\quad \text{versus} \quad 
K_{2,{\rm ex}}:\xi_i\neq1.
\end{align*}
Denote the subsample estimator based on $\{Z_t,\ldots,Z_{t+b-1}\}$ and the full sample estimator by  $\hat \xi_{b,t}^{(i)}$ and $\hat\xi_{n}^{(i)}$, respectively. 
Our proposed test rejects $H_{2,{\rm ex}}$ whenever $p_{n,2,{\rm ex}}\leq\varphi$, where
\begin{align*}
p_{n,2,{\rm ex}}\coloneqq \frac{1}{n-b+1}\sum_{t=1}^{n-b+1}\mathbb I
\left\{
\left|\sqrt b \skakko{\hat\xi_{b,t}^{(i)}-1}\right|
>
\left|\sqrt n \skakko{\hat\xi_{n}^{(i)}-1}\right|
\right\}.
\end{align*}

{It is essential to note that in cases where data undergoes differencing before analysis, the risk of over-differencing---which leads to the presence of an MA unit root---must be carefully considered.} 
In this regard, \cite{vl18} and \cite{v22} assumed the absence of MA unit roots in the context of the estimation of ARMA parameters based on the higher-order classical spectrum and the generalized spectrum, respectively.

\end{example}

\subsection{One-sided hypothesis}
{Next, we consider for $i \in \{1, \ldots, d \}$ the one-sided hypothesis testing problem formulated by}
\begin{align*}
H_3:{\theta}_{0i}\leq\kappa_i
\quad \text{versus} \quad 
K_3:{\theta}_{0i}>\kappa_i.
\end{align*}
The subsampling based $p$-value can be defined as
\begin{align*}
p_{n,3}\coloneqq &\frac{1}{n-b+1}\sum_{t=1}^{n-b+1}\mathbb I
\left\{
\sqrt b \skakko{\hat{\theta}_{b,t}^{(i)}-\kappa_i}
>
\sqrt n \skakko{\hat{\theta}_{n}^{(i)}-\kappa_i}
\right\},
\end{align*}
where ${\hat{\bm\theta}}_{b,t}\coloneqq (\hat{\theta}_{b,t}^{(1)},\ldots,\hat{\theta}_{b,t}^{(d)})^\top$ is the subsample estimator based on $\{Z_t,\ldots,Z_{t+b-1}\}$ and $\bm{\hat{\theta}}_{n}\coloneqq (\hat{\theta}_{n}^{(1)},\ldots,\hat{\theta}_{n}^{(d)})$ is the full sample estimator.

\begin{theorem}\label{thm4}
{Assume that the assumptions of Theorem \ref{thm2} hold. Consider the strictly stationary process $\{Z_t\}$ is $\alpha$-mixing with vanishing $\alpha$-mixing coefficient, and the sub-sampling block length $b$ satisfies $b/n\to0$ as $n\to\infty$.} The test which rejects $H_3$ in favor of $K_3$ whenever $p_{n,3}\leq\varphi$ has { asymptotic} level $\varphi$ and is consistent.
\end{theorem}

Similarly, we can test the hypothesis 
\begin{align*}
H_4:{\theta}_{0i}\geq\kappa_i
\quad \text{versus} \quad 
K_4:{\theta}_{0i}<\kappa_i.
\end{align*}
by the test which rejects  $H_4$ whenever $p_{n,4}\leq\varphi$, where
\begin{align*}
p_{n,4}\coloneqq &\frac{1}{n-b+1}\sum_{t=1}^{n-b+1}\mathbb I
\left\{ 
\sqrt b \skakko{\kappa_i-\hat{\theta}_{b,t}^{(i)}}
>
\sqrt n \skakko{\kappa_i-\hat{\theta}_{n}^{(i)}}
\right\}.
\end{align*}

We conclude this section by demonstrating how our procedures can be used to test for non-invertibility.
In practice, practitioners typically aim to reject the null hypothesis of non-invertibility for the observed process.

\begin{example}
Consider MA($d_1+d_2$) models with Cauchy residuals in Example \ref{ex:unitroot}.
If $|\xi_1|,\ldots,|\xi_{d_1+d_2}|>1$,  we have that 
\begin{align*}
\epsilon_t=&\frac{1}{(1-\xi_1^{-1}B)\cdots(1-\xi_{d_1}^{-1}B)
(1-\xi_{d_1+1}B^{-1})\cdots(1-\xi_{d_1+d_2}B^{-1})}Z_t\\
=&
\skakko{\sum_{j=0}^\infty\xi_1^{-j}B^j}
\cdots
\skakko{\sum_{j=0}^\infty\xi_{d_1}^{-j}B^j}
\skakko{-\sum_{j=0}^\infty\xi_{d_{1+1}}^{-j-1}B^{j+1}}
\cdots
\skakko{-\sum_{j=0}^\infty\xi_{d_{1}+d_2}^{-j-1}B^{j+1}}Z_t
\end{align*}
and thus the process is invertible. Therefore, the test for non-invertibility can be formulated as
\begin{align*}
H_{3,{\rm ex}}:|\xi_i|\leq1 
\quad \text{versus} \quad 
K_{3,{\rm ex}}:|\xi_i|>1.
\end{align*}
Then we propose the test rejects $H_{3,{\rm ex}}$ whenever $p_{n,3,{\rm ex}}\leq\varphi$,
\begin{align*}
p_{n,3,{\rm ex}}\coloneqq \frac{1}{n-b+1}\sum_{t=1}^{n-b+1}\mathbb I
\left\{
\sqrt b \skakko{\left|\hat\xi_{b,t}^{(i)}\right|-1}
>
\sqrt n \skakko{\left|\hat\xi_{n}^{(i)
}\right|-1}
\right\}.
\end{align*}
Note that \cite{cce17} considered a test for invertibility for vector ARMA model based on the generalized spectrum for residuals but they do not consider the unit-root case.
\end{example}

\section{Simulation studies}\label{sec7}

In this section, we illustrate the validity of our mathematical results through some simulation studies. We cover both estimation and hypothesis testing. One of the main goals is to show that our asymptotic results are of practical use, highlighting that a reasonably large sample size is enough for our proposed procedures to perform well. Our simulation setups are designed to illustrate the robustness to the presence of heavy-tailed innovations.
\subsection{Estimation} 

\subsubsection{Real-valued models} We first consider the continuous case with heavy-tailed innovations. More precisely, we consider the MA$(1)$ and AR$(1)$ processes with Cauchy innovations, defined in Section \ref{sec3}. In the case of the AR$(1)$ model, we consider both the causal and non-causal cases, as the generalized spectra have different closed form expressions in these two cases. As we showed in Section \ref{sec3} that the parameters $(a, \delta)$ of the families of generalized spectrum were identifiable in these two models, we estimate simultaneously all the components of the vectors of parameters. In the case of the AR$(1)$ model and for the sake of reducing computation time, we only consider the frequencies associated to $\ell\in \{-2,-1,0,1,2\}$ in the generalized spectrum.

For the (invertible) MA$(1)$ model, we generate $2000$ {time series} of {length} $n=500$ with parameter values $(a, \delta)=({0.8^{-1}}, 1)$. For the AR$(1)$ model, we generate $2000$ {time series} of {length} $n=500$ with parameter values $(a, \delta)=(0.7, 2)$ in the causal case and $(a, \delta)=(1.3, 2)$ in the non-causal case. These choices of parameters yield the same marginal distributions for $Z_t$ in the AR$(1)$ model, our aim being to illustrate that the proposed estimators distinguish between causal and non-causal processes, even in this challenging scenario. For each of the $2000$ {replications} we compute the estimator $\hat{\pmb \theta}_n=(\hat{a}_n, \hat{\delta}_n)^{\top}$ described in \eqref{Estimator} with parameter choices $L=3.14$ and $M_n=30$. In the case of the AR$(1)$ models, we {do} not assume that $\abs{a}$ is larger (or smaller) than {one} and run the minimization process for $f_{{\pmb \theta}}(\lambda, u, v)$, defined for all possible ${a} \in \mathbb{R} \setminus ([-1.1,-0.9]\cup [0.9,1.1])$.
For each of these three scenarios, we plot the histogram of $\hat{\pmb \theta}_n$, compared to the target Gaussian distribution. Due to the complicated nature of the asymptotic variance, we use estimators based on $\hat{\pmb \theta}_n$ to compute the mean and empirical variance of the Gaussian distribution. 
{Since the results for the MA(1), causal AR(1), and non-causal AR(1) models are  similar, we only present the figure for the non-causal AR(1) model here and defer the results for the MA(1) and causal AR(1) models to the Supplemental material.}
Inspection of Figure \ref{Fig3bis} {tends} to indicate that the limiting distribution of the proposed estimator is indeed asymptotically Gaussian. We choose to simply provide here the histograms of $\hat{\pmb \theta}_n$ instead of $\sqrt{n}(\hat{\pmb \theta}_n-{\pmb \theta}_0)$, as the actual dispersion of the estimator in the finite-sample scenario considered is easier to appreciate this way. Note that we choose to plot ${\hat{a}_n^{-1}}$ in the MA$(1)$ case, to maintain consistency between the AR and MA plots. 
Our method achieved 100\% accuracy in distinguishing between causal and non-causal models based on the estimated value of $\hat a_n$ within our experimental setup.

\begin{figure}[!htbp]
\vspace{-0mm}
\centering
\includegraphics[width=0.8\textwidth]{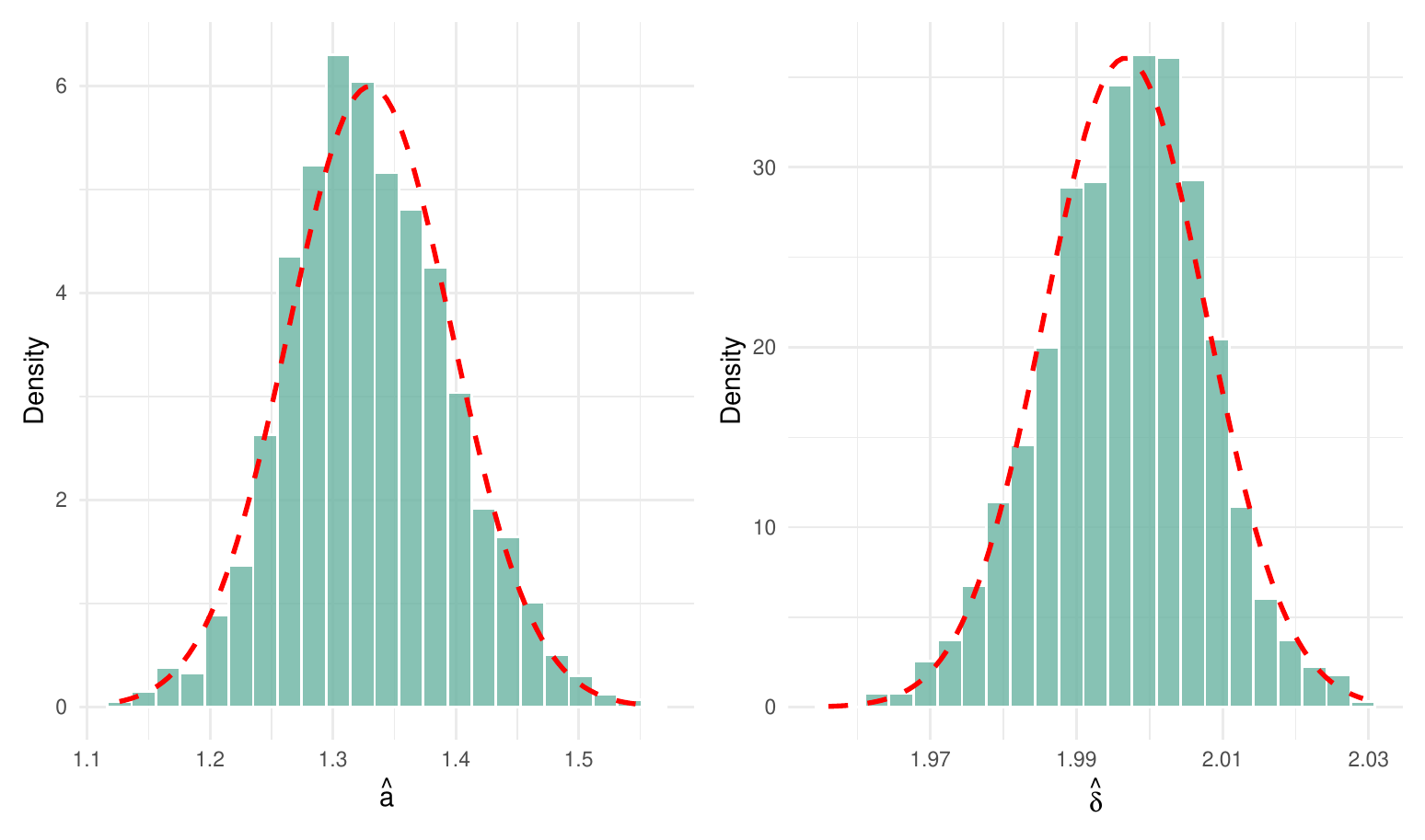}
\vspace*{-0mm}
  \caption{Histograms of ${\hat{a}_n}$ and $\hat{\delta}_n$ (in green) with respect to limiting Gaussian distribution (in red) in the non-causal AR$(1)$ model with Cauchy innovations. $2000$ replication, $n=500$ and $(a, \delta)=(1.3, 2)$.}
  \label{Fig3bis}
\end{figure}

\subsubsection{Integer-valued models} We now consider heavy-tailed integer-valued models, that is, the INMA$(1)$ and INAR$(1)$ processes with discrete stable innovations, defined in Section \ref{sec3}. Again, it has been shown in Section \ref{sec3} that the parameters $(\alpha, p, \delta)$ of the families of generalized spectra were identifiable in these two models. In the case of the INAR$(1)$ model we again only consider the frequencies associated to $\ell\in \{-2,-1,0,1,2\}$.

It has been observed that in augmented GARCH models, estimation of an exponent parameter (partially) determining the heaviness of the tails of the conditional distribution is surprisingly hard. It has been shown in \cite{hz2011} that even if asymptotic normality of the estimators could be derived, finite-$n$ performance could be rather problematic due to the flatness of the objective function. As, in the integer-valued models we consider, $\alpha$ is also an exponent parameter tied to the heaviness of the tails, we suspect a similar phenomenon could be observed when trying to estimate it. In augmented GARCH models, it has been proposed by \cite{fz2025} to restrict the search of $\alpha$ to a few values representing different practically relevant scenarios. Following this rationale, we consider $\alpha \in \{0.3,0.7, 0.9\}$, corresponding respectively to the presence of {\it ultra heavy-tailed}, {\it strongly heavy-tailed} and {\it moderately heavy-tailed} innovations. We do not impose any restriction of this type on the parameters $p$ and $\delta$. For the INMA$(1)$ model and INAR$(1)$ models, we generate $2000$ {time series} of {length} $n=500$ with both parameter values taken as $(\alpha, p, \delta)=(0.7, 0.3, 2)$. We assume that $\alpha \in \{0.3,0.7, 0.9\}$, $p\in (0,1)$, $\delta \in \mathbb{R}^+_0$ and, for each of the $2000$ {replications} we compute the estimator $\hat{\pmb \theta}_n=(\hat{\alpha}_n, \hat{p}_n, \hat{\delta}_n)^{\top}$ described in \eqref{Estimator} with parameter choices $L=3.14$ and $M_n=30$. In the INAR$(1)$ and INMA$(1)$ models, we achieve a 100\% selection rate with $\hat{\alpha}_n=0.7$ for each sample. We then plot the histogram of $(\hat{p}_n, \hat{\delta}_n)$, compared to the target Gaussian distribution (with Monte-Carlo mean and variance). 
{Since the results for the INMA(1) and INAR(1) models are similar, we only present the figure for the INAR(1) model here and defer the result for the INMA(1) model to the Supplemental material.}
Inspection of Figure \ref{Fig5} {tends} to indicate that the limiting distribution of the proposed estimator is asymptotically Gaussian. We should highlight here that $n=500$ is a rather reasonable sample size and that due to the complex nature of the optimization task and the heaviness of the tails, it is not very surprising that the histograms do not perfectly match the target Gaussian density function. With respect to such observations, the fit appears to be here very satisfactory.

\begin{figure}[!htbp]
\vspace{-0mm}
\centering
\includegraphics[width=0.8\textwidth]{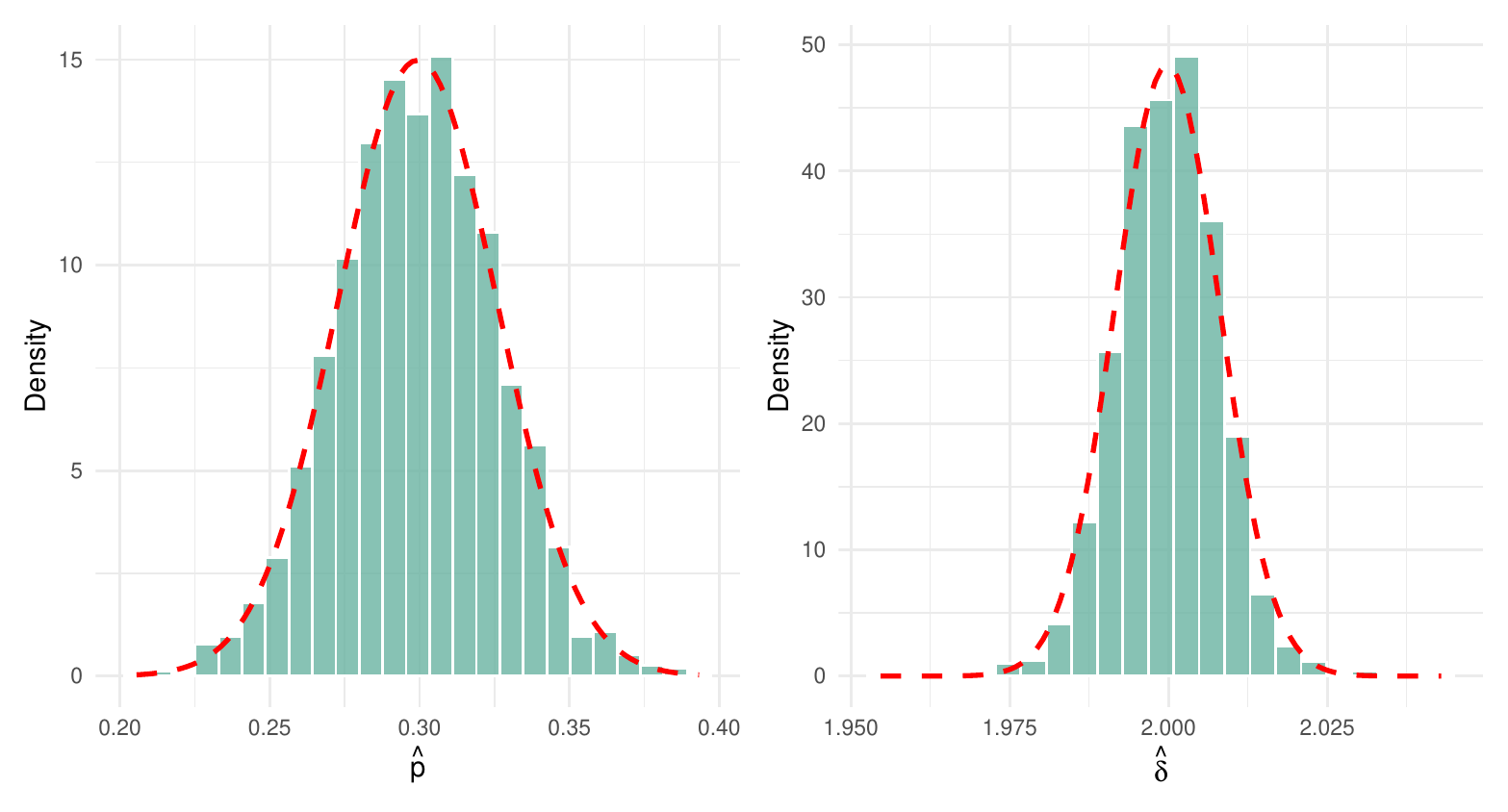}
\vspace*{-0mm}
  \caption{Histograms of ${\hat{p}_n}$ and $\hat{\delta}_n$ (in green) with respect to limiting Gaussian distribution (in red) in the INAR$(1)$ model with discrete stable innovations. $2000$ replication, $n=500$ and $(p, \delta)=(0.3, 2)$. The rate of selection of the proper $\alpha=0.7$ is a hundred percent.}
  \label{Fig5}
\end{figure}

\subsection{Hypothesis testing}

It is well known that subsampling based methods can have rather problematic behavior in small samples scenarios. The main potential limitation of the asymptotic results derived in Section \ref{sec6} is that they are likely not to translate very well in practical scenarios, where the sample is usually not extremely large. The aim of this subsection is then to investigate if, even with a relatively small sample, the hypothesis tests we proposed still behave reasonably closely to what is expected given the asymptotic results.

First, we study the goodness-of-fit test {for the null hypothesis that the time series belongs to the parametric family of INAR(1) models with discrete stable errors (see Section \ref{sec3} for a closed form), against the alternative that it does not.}
We first simulate $400$ {time series} of {length} $n \in \{100,300\}$ following the INAR$(1)$ model with parameters $(\alpha, p, \delta)=(0.7, 0.3, 2)$, {corresponding to the null hypothesis}, and then {generate} $400$ {time series} of {length} $n \in \{100,300\}$ following the time-varying INAR$(1)$  model with parameters
$$(\alpha, p, \delta)=
\begin{cases}
(0.7, 0.3, 2), & t\le n/2,\\
(0.7, 0.7, 2),  & t> n/2,
\end{cases}
$$
{corresponding to the alternative.}
We implement the goodness-of-fit test described in Section \ref{sec5}. 
We consider subsampling blocks of length $25$ when $n=100$ and of size $30$ when $n=300$. We also decide to compute the estimator of Section \ref{sec4} with $M_n=30$. {The choice of rather small block sizes aims at making it easier to illustrate the consistency of the test under the alternative considered. Indeed, larger block sizes seemed to lead to a slightly slower convergence speed of the empirical rejection probability to the asymptotic power. Our aim here is primarily to demonstrate that in practical scenarios with reasonably small $n$, the proposed test are of practical use.} The parameter $L$ is set to $3.14$. We performed the test at the asymptotic level $\alpha=0.05$.

Then, we consider the unit-root test {for the null hypothesis that the time series has an MA unit root against the alternative that it does not.} We first simulate $400$ {time series} of {length} $n \in \{100,300\}$ following the MA$(1)$ model with Cauchy innovation (see Section \ref{sec3}) with parameters $(a, \delta)=(1,2)$, {corresponding to the null hypothesis}, and then $400$ {time series} of {length} $n \in \{100,300\}$ with parameters $(a, \delta)=(0.3,2)$, {corresponding to the alternative hypothesis}. We implement the unit-root test described in Section \ref{sec6}.
The various parameters in the tests are set in the same manner as for the goodness-of-fit tests described earlier {but this time we consider blocks of length $50$, as illustrating the consistency of the test with larger block sizes does not appear to be complicated here}.

For each of the scenarios, the empirical rejection {probabilities} are indicated in Table \ref{Table1}. The inspection of the table tends to indicate that, even with a very small sample size $n=100$, the proposed tests are reasonably close to what is expected when considering their limiting distribution. Indeed, the empirical rejection probabilities under the null are satisfyingly close to the level constraint, while both tests seem to enjoy some power against the considered alternatives.

\begin{table}[ht]
\centering
\caption{Empirical rejection {probabilities} ($\alpha=0.05$)}
\label{Table1}
\begin{tabular}{lcccc}
\hline
 & null ($n=100$) & alternative ($n=100$)  & null ($n=300$) & alternative ($n=300$)  \\
\hline
goodness-of-fit test & 0.105 & 0.260 & 0.0650 & 0.660 \\
MA unit-root test & 0.0900 & 1.000 & 0.0550 & 1.000 \\
\end{tabular}
\end{table}

\section{Empirical Study}\label{sec8}
To illustrate the practical utility of the proposed model, we analyze a measles cases count dataset ($n=646$), which exhibits empirical evidence of an infinite mean; see Figure~\ref{fig:ts_hill}. This dataset represents the weekly count of cases of measles observed in the Lander of North Rhine-Westphalia (Germany) between January 2001 and May 2013 and is available in the {\it tscount} R package (see \citealt{tscount}). It should be noted that quick observation of Figure~\ref{fig:ts_hill} tends to indicate that the dataset exhibits typical behavior of integer-valued time series, with zero-inflation, over-dispersion and bubble-like phenomenon. Due to the way infectious diseases propagate, it is extremely natural to assume an underlying autoregressive structure. Motivated by the behavior of the Hills estimator that hints towards the presence of heavy-tailed innovaation, we fit the discrete-stable INAR(1) model to the data.

\subsection{Parameter estimation and goodness-of-fit}{
We divided the sample into a training set that contains the first $n_{\rm train}=400$ data points and a testing set that contains the last $n_{\rm test}=246$ data points. The model parameters are estimated using the proposed estimation procedure \eqref{Estimator} with parameters $L=3.14$ and $M_n=30$. The estimated parameters are
$$
\hat{\delta} = 0.283, \quad \hat{\alpha} = 0.364, \quad \text{and} \quad \hat{p} = 0.560.
$$

According to Remark~\ref{rem:INAR}, the estimated process $\{Z_t\}$ follows a discrete stable  distribution.
In particular, the resulting marginal distribution of $Z_t$ is discrete stable with exponent { $\hat{\alpha}_{\rm marginal} = 0.364$ and the scale parameter $\hat{\delta}_{\rm marginal} = 1.487$.}
Notably, the estimated exponent is substantially smaller than unity, indicating that the process admits only fractional moments. This observation justifies our departure from standard Poisson or negative binomial INAR models, which implicitly assume the existence of all moments. { The fact that the scale parameter $\hat{\delta}_{\rm marginal}=1.487$ of the marginal distribution of $Z_t$ is rather small is in line with the presence of zero-inflation in the dataset.}

We then conduct the goodness-of-fit test introduced in Section~\ref{sec5}. We perform this test using the training set only. We use a block length of size $40$ and the exact same estimation procedure as before with $L=3.14$ and $M_n=30$. The resulting p-value is $p_{n,1} = 0.271
$. Hence, at a conventional significance level such as $0.05$, we do not reject the null hypothesis $H_1$, which provides no statistical evidence to reject the discrete stable INAR(1) model. Note that the sample size, block sizes and various parameter values are reasonably close to the simulation setups that have been studied in Section \ref{sec7}. We then have good reason to believe in the validity of our conclusion.

\subsection{Prediction}

We next assess the one-step-ahead predictive performance of our model, using the training subsample of size $n_{\rm {test}} = 246$ and the estimated parameters described in the previous subsection. Since the mean of the process is not well defined, we adopt the median-based predictor
$$
\hat{Z}_{t|t-1}
=
\operatorname{median}\!\left(\mathrm{Bin}(Z_{t-1},\hat{p})
+\mathrm{DS}(\hat{\delta}, \hat{\alpha})\mid Z_{t-1}\right),
$$
where $\operatorname{median}(\mathcal D\mid Z)$ denotes the median of a distribution $\mathcal D$, conditionally to the random variable $Z$,
$\mathrm{Bin}(Z_{t-1},\hat{p})$ denotes the binomial distribution with number of trials $Z_{t-1}$ and success probability $\hat{p}$, and
$\mathrm{DS}(\hat{\delta}, \hat{\alpha})$ denotes the discrete stable distribution with scale parameter $\hat{\delta}$ and exponent $\hat{\alpha}$. Note here that the lack of linearity of the conditional median with respect to past observations makes it impossible to express the predictor as a simple shift of the median of the innovation, as it would be the case with an expectation-based predictor.
The predictive accuracy is evaluated using the mean squared prediction error (MSPE),
$$
\mathrm{MSPE}
=
\frac{1}{n_{\text{test}}}
\sum_{t=n_{\text{train}}+1}^{n}
(\hat{Z}_{t\mid t-1} - Z_t)^2=9.959
.
$$

The value of the MSPE seems reasonably small, accounting for the over-dispersion of the process which tends to generate large errors. For for the sake of interpretability we plot the predicted one-step ahead time series with respect to the true values. Inspection of Figure \ref{Fig6} indicates that our predictor behaves in a satisfactory way.
\begin{figure}[!htbp]
\vspace{-0mm}
\centering
\includegraphics[width=0.8\textwidth]{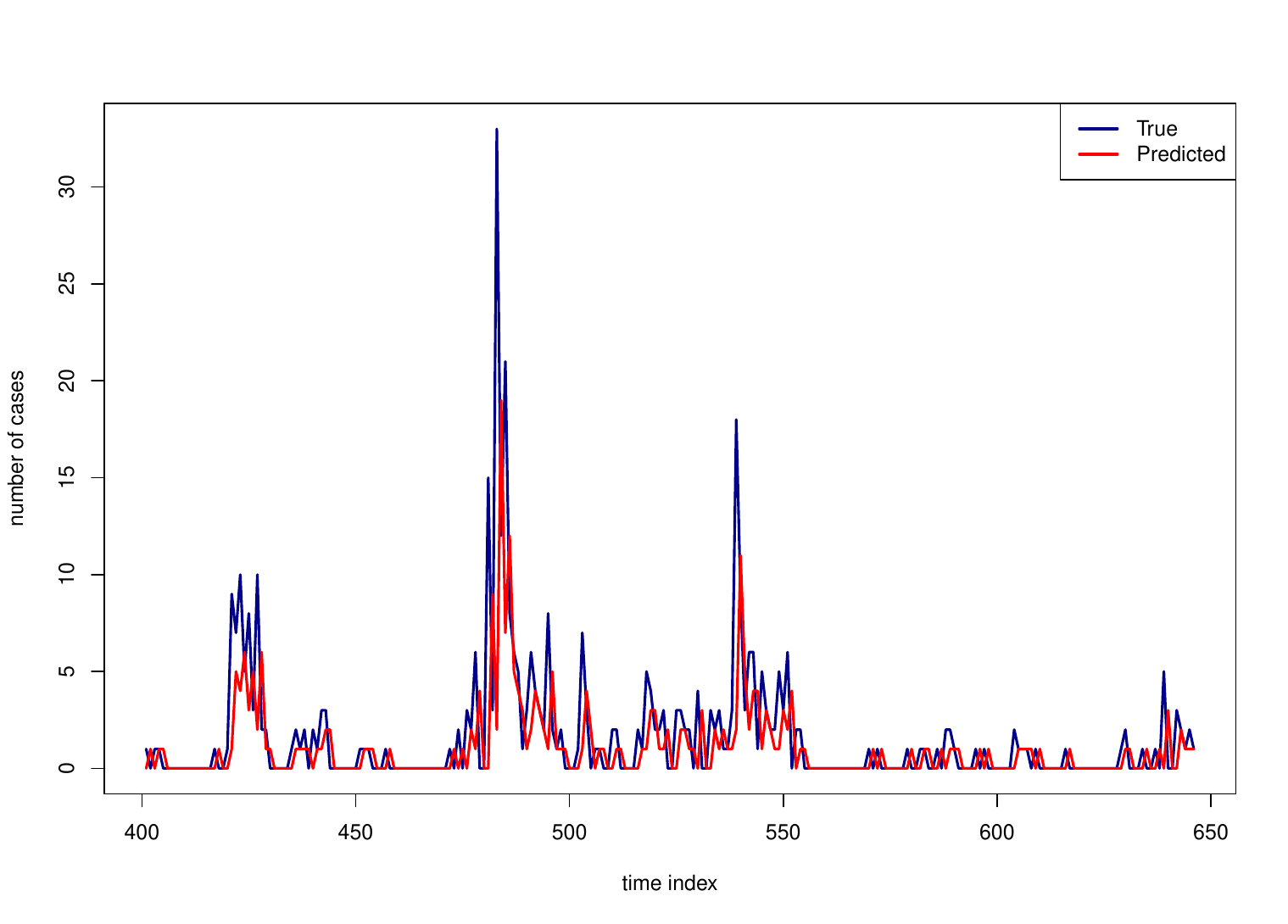}
\vspace*{-0mm}
  \caption{One-step ahead predicted values of the measles data count time series (in red) and true value of the time series (in blue). The testing set size is $n_{\rm test}=246$.}
  \label{Fig6}
\end{figure}

Now, we compare our predictor to two natural parametric candidates based on the INAR$(1)$ model. The first one is a one-step ahead predictor that assumes Poisson innovations, while the second one assumes a negative binomial model. As both of these models implicitly assume finite moments of order $1$, it makes sense to consider not only a median-based predictor (as in our previous construction) but also a classical expectation-based one. We estimate the parameters using the spINAR R package (see \cite{spinar} for more details). The MSPE for each method are displayed in Table \ref{Table2}. We observe that our predictor performs significantly better than the competitors. The slightly worst performances of expectation-based predictors are in line with the fact that heavy-tailed models {provide a better fit to the data set}.

\begin{table}[ht]
\centering
\caption{Performances of the various predictors}
\label{Table2}
\begin{tabular}{lcc}
\hline
model & predictor & MSPE \\
\hline
discrete stable INAR$(1)$ & median-based & 9.959
\\
Poisson INAR$(1)$ &  median-based & 22.500  \\
negative binomial INAR$(1)$  &  median-based & 30.500  \\
Poisson INAR$(1)$  &  expectation-based & 22.711
\\
negative binomial INAR$(1)$ &  expectation-based & 49.102
\\
\hline
\end{tabular}
\end{table}
}

\begin{funding}
This research was supported by JSPS Grant-in-Aid for Early-Career Scientists, Grant Number JP23K16851 (Y.G.), and the Research Fellowship Promoting International Collaboration of the Mathematical Society of Japan (Y.G.).
Additionally, the authors used ChatGPT and Gemini for proofreading and to enhance the clarity of the manuscript during the writing process. 
The authors maintain full responsibility for the content and accuracy of the final manuscript.
\end{funding}

\begin{supplement}
\stitle{Supplement to ``Parametric generalized spectrum for heavy-Tailed time series''}
\sdescription{
Appendix A contains the proofs of the equations, lemmas, and theorems.
Appendix B presents additional plots of the generalized spectra.
Additional simulation results for Section \ref{sec7} are reported in Appendix C.
}
\end{supplement}

\bibliographystyle{apalike}
\bibliography{ref}

@article{hv25,
  title={Non-causal and non-invertible ARMA models: Identification, estimation and application in equity portfolios},
  author={Hecq, Alain and Velasquez-Gaviria, Daniel},
  journal={J. Time Ser. Anal.},
  volume={46},
  number={2},
  pages={325--352},
  year={2025},
  publisher={Wiley Online Library}
}

@article{spinar,
  title={FspINAR: An R Package for Semiparametric and Parametric Estimation and Bootstrapping of Integer-Valued Autoregressive (INAR) Models},
  author={Faymonville, Maxime and Riffo, Javiera and Rieger, Jonas and Jentsch, Carsten},
  journal={arXiv:2401.14239},
  year={2024}
}

@article{fz2025,
  title={Finite moments testing in a general class of
nonlinear time series models},
  author={Franck, Christian and Zakoïan, Jean-Michel},
  journal={Bernoulli},
  volume={31},
  pages={2649--1674},
  year={2025}
}

@article{hz2011,
  title={Asymptotic properties of LS and QML estimators for a class of
nonlinear GARCH processes},
  author={Hamadeh, Tawfik and Zakoïan, Jean-Michel},
  journal={J. Statist. Plann. Inference},
  volume={141},
  pages={488--507},
  year={2011}
}

@article{tscount,
  title={tscount: An {R} package for analysis of count time series following generalized linear models},
  author={Liboschik, Tobias and Fokianos, Konstantinos and Fried, Roland},
  journal={J. Stat. Softw.},
  volume={82},
  pages={1--51},
  year={2017}
}

@book{r00,
  title={Gaussian and non-Gaussian linear time series and random fields},
  author={Rosenblatt, Murray},
  year={2000},
  publisher={Springer Science \& Business Media}
}

@article{m88,
  title={Mixing properties of ARMA processes},
  author={Mokkadem, Abdelkader},
  journal={Stochastic Process. Appl.},
  volume={29},
  number={2},
  pages={309--315},
  year={1988},
  publisher={Elsevier}
}

@book{lp08,
  title={Markov chains and mixing times},
  author={Levin, David A and Peres, Yuval},
  year={2017},
  publisher={Amer. Math. Soc.}
}

@article{bkp13,
  title={Heavy-Tailed Branching Process with Immigration},
  author={Basrak, Bojan and Kulik, Rafal and Palmowski, Zbigniew},
  journal={Stoch. Models},
  volume={28},
  number={},
  pages={413--434},
  year={2013},
  publisher={}
}

@article{szu24,
  title={Ergodic properties of subcritical multitype Galton-Watson processes with immigration},
  author={Sz\"ucs, G\'abor},
  journal={ESAIM: Probab. and Statist.},
  volume={28},
  number={},
  pages={350--365},
  year={2024},
  publisher={}
}

@article{brad05,
  title={Basic Properties of Strong Mixing Conditions. A Survey and Some Open Questions},
  author={Bradley, Richard C.},
  journal={Probab. Survey},
  volume={2},
  number={},
  pages={107--144},
  year={2005},
  publisher={JSTOR}
}

@article{cce17,
  title={Testing for fundamental vector moving average representations},
  author={Chen, Bin and Choi, Jinho and Escanciano, Juan Carlos},
  journal={Quant. Econ.},
  volume={8},
  number={1},
  pages={149--180},
  year={2017},
  publisher={Wiley Online Library}
}

@article{km94,
  title={Some limit theory for the self-normalised periodogram of stable processes},
  author={Kl{\"u}ppelberg, Claudia and Mikosch, Thomas},
  journal={Scand. J. Stat.},
  volume={21},
  number={4},
  pages={485--491},
  year={1994},
  publisher={JSTOR}
}

@article{mgkt95,
  title={Parameter estimation for ARMA models with infinite variance innovations},
  author={Mikosch, Thomas and Gadrich, Tamar and Kluppelberg, Claudia and Adler, Robert J},
  journal={Ann. Statist.},
  volume={23},
  number={1},
  pages={305--326},
  year={1995},
  publisher={JSTOR}
}

@article{ky02,
  title={Empirical characteristic function in time series estimation},
  author={Knight, John L and Yu, Jun},
  journal={Econom. Theory},
  volume={18},
  number={3},
  pages={691--721},
  year={2002},
  publisher={Cambridge University Press}
}

@article{cs98,
  title={Discrete stable random variables},
  author={Christoph, Gerd and Schreiber, Karina},
  journal={Stat. Probab. Lett.},
  volume={37},
  number={3},
  pages={243--247},
  year={1998},
  publisher={Elsevier}
}

@article{sv79,
  title={Discrete analogues of self-decomposability and stability},
  author={Steutel, Fred W and van Harn, Klaas},
  journal={Ann. Probab.},
  pages={893--899},
  year={1979},
  publisher={JSTOR}
}

@article{mpk20,
author = {Marco Meyer and Efstathios Paparoditis and Jens-Peter Kreiss},
title = {Extending the validity of frequency domain bootstrap methods to general stationary processes},
volume = {48},
journal = {Ann. Statist.},
number = {4},
publisher = {Institute of Mathematical Statistics},
pages = {2404 -- 2427},
year = {2020}
}

@article{kp23,
  title={Bootstrapping Whittle estimators},
  author={Kreiss, J-P and Paparoditis, Efstathios},
  journal={Biometrika},
  volume={110},
  number={2},
  pages={499--518},
  year={2023},
  publisher={Oxford University Press}
}

@article{tpk96,
  title={Nonparametric approach for non-Gaussian vector stationary processes},
  author={Taniguchi, Masanobu and Puri, Madan L and Kondo, Masao},
  journal={J. Multivariate Anal.},
  volume={56},
  number={2},
  pages={259--283},
  year={1996},
  publisher={Elsevier}
}

@Manual{R,
    title = {R: A Language and Environment for Statistical Computing},
    author = {{R Core Team}},
    organization = {R Foundation for Statistical Computing},
    address = {Vienna, Austria},
    year = {2021},
    url = {https://www.R-project.org/},
  }

@article{dw16,
  title={Detecting relevant changes in time series models},
  author={Dette, Holger and Wied, Dominik},
  journal={J. R. Stat. Soc. Ser. B Methodol.},
  volume={78},
  number={2},
  pages={371--394},
  year={2016},
  publisher={Oxford University Press}
}

@Book{prw99,
  author    = {Dimitris N. Politis and Joseph P. Romano and Michael Wolf},
  title     = {Subsampling},
  year      = {1999},
  publisher = {Springer},
  location  = {New York},
}

@book{bd09,
  title={Time series: theory and methods},
  author={Brockwell, Peter J and Davis, Richard A},
  year={2009},
  publisher={Springer science \& business media}
}

@article{vl18,
  title={Frequency domain minimum distance inference for possibly noninvertible and noncausal ARMA models},
  author={Velasco, Carlos and Lobato, Ignacio N},
  journal={Ann. Statist.},
  volume={46},
  number={2},
  pages={555--579},
  year={2018},
  publisher={JSTOR}
}

@article{efp19,
  title={An updated literature review of distance correlation and its applications to time series},
  author={Edelmann, Dominic and Fokianos, Konstantinos and Pitsillou, Maria},
  journal={Int. Stat. Rev.},
  volume={87},
  number={2},
  pages={237--262},
  year={2019},
  publisher={Wiley Online Library}
}

@article{dmmw18,
  title={Applications of distance correlation to time series},
  author={Davis, Richard A and Matsui, Muneya and Mikosch, Thomas and Wan, Phyllis},
  year={2018},
  journal={Bernoulli},
    volume={24},
  number={4A},
  pages={3087--3116}
}

@article{wd22,
  title={Goodness-of-fit testing for time series models via distance covariance},
  author={Wan, Phyllis and Davis, Richard A},
  journal={J. Econometrics},
  volume={227},
  number={1},
  pages={4--24},
  year={2022},
  publisher={Elsevier}
}

@article{fp18,
  title={Testing independence for multivariate time series via the auto-distance correlation matrix},
  author={Fokianos, Konstantinos and Pitsillou, M},
  journal={Biometrika},
  volume={105},
  number={2},
  pages={337--352},
  year={2018},
  publisher={Oxford University Press}
}

@article{fp17,
  title={Consistent testing for pairwise dependence in time series},
  author={Fokianos, Konstantinos and Pitsillou, Maria},
  journal={Technometrics},
  volume={59},
  number={2},
  pages={262--270},
  year={2017},
  publisher={Taylor \& Francis}
}

@article{l97,
	author = {Latour, Alain},
	journal = {Adv. Appl. Probab.},
	number = {1},
	pages = {228--248},
	publisher = {Cambridge University Press},
	title = {The multivariate {GINAR} (p) process},
	volume = {29},
	year = {1997}}

@article{aa91,
	author = {M. A. Al-Osh and A. A. Alzaid},
	journal = {Comm. Statist. Theory Methods},
	number = {2},
	pages = {261-282},
	publisher = {Taylor & Francis},
	title = {Binomial autoregressive moving average models},
	volume = {7},
	year = {1991}}

@article{m85,
	author = {McKenzie, Ed},
	journal = {J. Am. Water Resour. Assoc.},
	number = {4},
	pages = {645--650},
	publisher = {Wiley Online Library},
	title = {Some simple models for discrete variate time series},
	volume = {21},
	year = {1985}}

@article{flo06,
	author = {Ferland, Ren{\'e} and Latour, Alain and Oraichi, Driss},
	journal = {J. Time Ser. Anal.},
	number = {6},
	pages = {923--942},
	publisher = {Wiley Online Library},
	title = {Integer-valued {GARCH} process},
	volume = {27},
	year = {2006}}

@Article{v20,
  author   = {von Sachs, Rainer},
  title    = {Nonparametric Spectral Analysis of Multivariate Time Series},
  number   = {1},
  pages    = {361-386},
  volume   = {7},
  journal  = {Annu. Rev. Stat. Appl.},
  year     = {2020},
}

@article{dfhllp21,
	author = {Davis, Richard A and Fokianos, Konstantinos and Holan, Scott H and Joe, Harry and Livsey, James and Lund, Robert and Pipiras, Vladas and Ravishanker, Nalini},
	journal = {J. Amer. Statist. Assoc.},
	number = {535},
	pages = {1533--1547},
	publisher = {Taylor \& Francis},
	title = {Count time series: A methodological review},
	volume = {116},
	year = {2021}}

@article{g20,
  title={Beta--negative binomial auto-regressions for modelling integer-valued time series with extreme observations},
  author={Gorgi, Paolo},
  journal={J. R. Stat. Soc. Ser. B Methodol.},
  volume={82},
  number={5},
  pages={1325--1347},
  year={2020},
  publisher={Oxford University Press}
}

@article{h13,
  title={Robust spectral analysis (ar{X}iv:1111.1965v2)},
  author={Hagemann, Andreas},
  journal = {Ar{X}iv {E}-prints},
  year={2013}
}

@ARTICLE{kvdh16,
  author = {Kley, Tobias and Volgushev, Stanislav and Dette, Holger and Hallin, Marc},
  title = {Quantile spectral processes: Asymptotic analysis and inference},
  journal = {Bernoulli},
  year = {2016},
  volume = {22},
  pages = {1770--1807},
  number = {3},
  publisher = {Bernoulli Society for Mathematical Statistics and Probability}
}

@ARTICLE{dhkv13,
  author = {Dette, Holger and Hallin, Marc and Kley, Tobias and Volgushev, Stanislav},
  title = {Of copulas, quantiles, ranks and spectra: An $L_1$-approach to spectral
	analysis},
  journal = {Bernoulli},
  year = {2015},
  volume = {21},
  pages = {781--831},
  number = {2},
  publisher = {Bernoulli Society for Mathematical Statistics and Probability}
}

@article{vvd18,
  title={Fourier analysis of serial dependence measures},
  author={Van Hecke, Ria and Volgushev, Stanislav and Dette, Holger},
  journal={J. Time Ser. Anal.},
  volume={39},
  number={1},
  pages={75--89},
  year={2018},
  publisher={Wiley Online Library}
}

@ARTICLE{lr11,
  author = {Lee, Junbum and Rao, Suhasini Subba},
  title = {The quantile spectral density and comparison based tests for nonlinear
	time series (ar{X}iv:1112.2759)},
	journal = {Ar{X}iv {E}-prints},
  year = {2012}
}

@article{l08,
  title={Laplace periodogram for time series analysis},
  author={Li, Ta-Hsin},
  journal={J. Amer. Statist. Assoc.},
  volume={103},
  number={482},
  pages={757--768},
  year={2008},
  publisher={Taylor \& Francis}
}

@ARTICLE{l12,
  author = {Li, Ta-Hsin},
  title = {Quantile periodograms},
  journal = {J. Amer. Statist. Assoc.},
  year = {2012},
  volume = {107},
  pages = {765--776},
  number = {498},
  publisher = {Taylor \& Francis}
}

@article{l21,
  title={Quantile-frequency analysis and spectral measures for diagnostic checks of time series with nonlinear dynamics},
  author={Li, Ta-Hsin},
  journal={J. R. Stat. Soc. Ser. C Appl. Stat.},
  volume={70},
  number={2},
  pages={270--290},
  year={2021},
  publisher={Wiley Online Library}
}

@article{gkvvdh22,
	author = {Goto, Yuichi and Kley, Tobias and Van Hecke, Ria and Volgushev, Stanislav and Dette, Holger and Hallin, Marc},
	journal = {Ann. Statist.},
	number = {6},
	pages = {3563--3591},
	publisher = {Institute of Mathematical Statistics},
	title = {The integrated copula spectrum},
	volume = {50},
	year = {2022}}

@article{c88,
  title={Weighted least squares estimators on the frequency domain for the parameters of a time series},
  author={Chiu, Shean-Tsong},
  journal={Ann. Statist.},
  pages={1315--1326},
  year={1988},
  publisher={JSTOR}
}

@article{pd13,
  title={Testing semiparametric hypotheses in locally stationary processes},
  author={Preuss, Philip and Vetter, Mathias and Dette, Holger},
  journal={Scand. J. Stat.},
  volume={40},
  number={3},
  pages={417--437},
  year={2013},
  publisher={Wiley Online Library}
}

@article{d85,
  title={Asymptotic normality of spectral estimates},
  author={Dahlhaus, Rainer},
  journal={J. Multivariate Anal.},
  volume={16},
  number={3},
  pages={412--431},
  year={1985},
  publisher={Elsevier}
}

@ARTICLE{h00,
  author = {Hong, Yongmiao},
  title = {Generalized spectral tests for serial dependence},
  journal = {J. R. Stat. Soc. Ser. B Methodol.},
  year = {2000},
  volume = {62},
  pages = {557--574},
  number = {3},
  publisher = {Wiley Online Library}
}

@ARTICLE{h99,
  author = {Hong, Yongmiao},
  title = {Hypothesis testing in time series via the empirical characteristic
	function: a generalized spectral density approach},
  journal = {J. Amer. Statist. Assoc.},
  year = {1999},
  volume = {94},
  pages = {1201--1220},
  number = {448},
  publisher = {Taylor \& Francis Group}
}

@article{d00,
  title={A likelihood approximation for locally stationary processes},
  author={Dahlhaus, Rainer},
  journal={Ann. Statist.},
  volume={28},
  number={6},
  pages={1762--1794},
  year={2000},
  publisher={Institute of Mathematical Statistics}
}

@article{v22,
  title={Estimation of time series models using residuals dependence measures},
  author={Velasco, Carlos},
  journal={Ann. Statist.},
  volume={50},
  number={5},
  pages={3039--3063},
  year={2022},
  publisher={Institute of Mathematical Statistics}
}

@book{v00,
  title={Asymptotic statistics},
  author={Van der Vaart, Aad W},
  volume={3},
  year={2000},
  publisher={Cambridge university press}
}

@article{gzkc23,
	author = {Goto, Yuichi and Zhang, Xuze and Kedem, Benjamin and Chen, Shuo},
  journal = {Ar{X}iv {E}-prints},
	title = {Residual Spectrum Applied in Brain Functional Connectivity (ar{X}iv:2305.1946)},
	year = {2023}}

@book{tk00,
  title={Asymptotic theory of statistical inference for time series},
  author={Taniguchi, Masanobu and Kakizawa, Yoshihide},
  year={2000},
  publisher={Springer Science \& Business Media}
}

@book{b81,
  title={Time series: Data Analysis and Theory},
  author={Brillinger, David R},
  year={1981},
  edition={},
  publisher={San Francisco: Holden-Day}
}

\newpage

\appendix
\begin{center}

{\Large  SUPPLEMENT TO}\\
{\Large ``PARAMETRIC GENERALIZED SPECTRUM\\ FOR HEAVY-TAILED TIME SERIES''}

\vspace{1.5em}

{\large
Yuichi Goto$^{1}$ \quad and \quad Gaspard Bernard$^{2}$
}

\vspace{1em}

{\small
$^{1}$Faculty of Mathematics, Kyushu University\\
\texttt{yuichi.goto@math.kyushu-u.ac.jp}

\vspace{0.5em}

$^{2}$Institute of Statistical Science, Academia Sinica\\
\texttt{gaspardbernard@as.edu.tw}
}

\end{center}
\vspace{1em}

This is the supplemental material for the article
``Parametric generalized spectrum for heavy-tailed time series''.
Appendix~\ref{sec9} contains the proofs of the equations, lemmas, and theorems.
Appendix~\ref{sec11} presents additional plots of the generalized spectra.
Appendix~\ref{sec12} reports additional simulation results for Section~7.

\section{Proofs}\label{sec9}

\subsection{Proof of Lemma \ref{IdentifiableCMA}}
If $f_{\bm \theta}(\lambda;u,v)=f_{{\bm \theta}'}(\lambda;u,v)$, all the Fourier coefficients must be equal. This implies that, considering the Fourier coefficient associated to $\exp\skakko{-i\lambda}$,
$$\frac{1}{2\pi}
\exp\skakko{-\delta(|va^{-1}|+|v+a^{-1}u|+|u|)}=\frac{1}{2\pi}
\exp\skakko{-\delta'(|v(a')^{-1}|+|v+(a')^{-1}u|+|u|)}.$$
Now, note that the left hand side of this identity is not differentiable on the set $$S_a=\{(u,v) \mid v=-a^{-1}u\} \cup \{(u,v) \mid v=0\} \cup \{(u,v) \mid u=0\}$$ and that the right hand side is not differentiable on the set $$S_{a'}=\{(u,v) \mid v=-(a')^{-1}u\} \cup \{(u,v) \mid v=0\} \cup \{(u,v) \mid u=0\}.$$ Obviously, the sets $S_a$ and $S_{a'}$ must coincide. Note that for $\{(u,v) \mid v=-a^{-1}u\}=\{(u,v) \mid v=-(a')^{-1}u\}$ to hold, we must have $a'=a$.

Considering now the Fourier coefficient associated to $(\cos(\lambda)+1)$ yields
$$-\frac{1}{2\pi}
\exp\skakko{-\delta (|u|+|v|)(|a^{-1}|+1)}=-\frac{1}{2\pi}
\exp\skakko{-\delta'(|u|+|v|)(|(a')^{-1}|+1)}.$$
For $u=v \ge 0$, it gives
$$-\frac{1}{2\pi}
\exp\skakko{-\delta2u(|a^{-1}|+1)}=-\frac{1}{2\pi}
\exp\skakko{-\delta'2u(|(a')^{-1}|+1)}.$$
Taking the right derivative with respect to $u \ge 0$ on both sides yields
$$\frac{2\delta(|a^{-1}|+1)}{2\pi}
\exp\skakko{-\delta2u(|a^{-1}|+1)}=\frac{2\delta'(|(a')^{-1}|+1)}{2\pi}
\exp\skakko{-\delta'2u(|(a')^{-1}|+1)}.$$
Evaluating this last expression in $u=0$ gives that $\delta(|a^{-1}|+1)=\delta'(|(a')^{-1}|+1)$. Using $a=a'$ gives $\delta=\delta'$.
\qed

\subsection{Proof of Lemma \ref{IdentifiableCAR}}
Suppose $f_{\bm \theta}(\lambda;u,v)
=f_{\bm \theta'}(\lambda;u,v)$ for $\bm \theta=(a,\delta)$ and $\bm \theta'=(a',\delta')$.
This implies
the equality of the corresponding Fourier coefficients
$C_{\ell,a,\delta}(u,v)=C_{ \ell,a',\delta'}(u,v)$ for any $\ell,u,v$.

Suppose that $|a|<1$ and $|a'|>1$.
$C_{0,a,\delta}(1,-1)=
C_{0,a',\delta'}(1,-1)$
gives
$\gamma=\gamma'$,
where
$\gamma={\delta}/{(1-|a|)}$ and $\gamma'= {\delta |a^{-1}|}/{(1-|a^{-1}|)}$.
Evaluating the identity $C_{1,a,\delta}(1,v)=C_{1,a',\delta'}(1,v)$ yields
\begin{align*}
|a+v|
+
1-|a|
=
|v{a'}^{-1}+1|
+
|v|(1-|{a'}^{-1}|).
\end{align*}
The left hand side is non-differentiable at $v=-a$, while the right hand side is non-differentiable at $v=0,-a'$.
This is contradiction.
Thus, $a$ and $a'$ must be in the same region, i.e., both $|a|,|a'|<1$ or both $|a|,|a'|>1$.

Suppose that $|a|,|a'|<1$. The equation $C_{0,a,\delta}(1,-1)=
C_{0,a',\delta'}(1,-1)$
yields
$\gamma=\gamma''$, where $\gamma''={\delta'}/{(1-|a'|)}$.
From the identity
$C_{1,a,\delta}(1,v)=C_{1,a',\delta'}(1,v)$, it follows that
$|a+v|
+
1-|a|
=
|a'+v|
+
1-|a'|
$.
The function on the left-hand side is non-differentiable at $v=-a$, whereas the function on the right-hand side is non-differentiable  at $v=-a'$. Their points of non-differentiability must coincide. Thus, $a = a'$. 
Finally, combined with $\gamma=\gamma'$, we obtain $\delta = \delta'$.

A similar argument establishes the identifiability for the case $|a|, |a'| > 1$
\qed

\subsection{Proof of Lemma \ref{CAR_mixing}}
The geometrically $\beta$-mixing property of the causal AR model follows by \citet[Theorem 1']{m88}.
For the non-causal AR model, $Z_t$ can be represented as $Z_t = - \sum_{j=1}^\infty b^j \epsilon_{t+j}$, where $b=a^{-1}$ and $|b|<1$, it suffices to check the conditions of \citet[Theorem 4.4.1]{r00}.

First, we verify the moment condition. For any $\gamma \in (0,1)$, using the inequality $|\sum x_i|^\gamma \leq \sum |x_i|^\gamma$, we have
${\rm E}[|Z_t|^\gamma] 
\leq \sum_{j=1}^\infty |b^\gamma|^{j} {\rm E}[|\epsilon_{t+j}|^\gamma] < \infty$.

Next, we check the condition on the probability density function. Let $p_{\delta}(x)$ denote the pdf of the Cauchy distribution with scale parameter $\delta > 0$. For any $x$, using Fubini's theorem, we have
\begin{align*}
\int_{-\infty}^\infty
\left|
p_\delta(\xi+x)-p_\delta(\xi)
\right|
{\rm d}\xi
\leq&
\int_{-\infty}^\infty
\int_{\xi}^{\xi+x}
\left|
{p_\delta^{\prime}}(y)
\right|
{\rm d}y
{\rm d}\xi\\
\leq&
\int_{-\infty}^\infty
\int_{-\infty}^\infty
\mathbb I_{\{\xi \leq y\leq \xi+x\}}
\left|
{p_\delta^{\prime}}(y)
\right|
{\rm d}\xi
{\rm d}y.
\end{align*}
Since $\mathbb I_{\{\xi \leq y\leq \xi+x\}}=\mathbb I_{\{y-x\leq \xi \leq y \}}$, the integration with respect to $\xi$ yields the length $|x|$.
Thus, we obtain
\begin{align*}
\int_{-\infty}^\infty
\int_{-\infty}^\infty
\mathbb I_{\{\xi \leq y\leq \xi+x\}}
\left|
{p_\delta^{\prime}}(y)
\right|
{\rm d}\xi
{\rm d}y
\leq
|x|
\int_{-\infty}^\infty
\left|
{p_\delta^{\prime}}(y)
\right|
{\rm d}y
=
\frac{2}{\pi \delta}|x|.
\end{align*}
This satisfies the variation condition (4.4.5) in \cite{r00}.

Now we evaluate the decay rate $W(k, \gamma) = \max\{A_k, B_k\}$, where
$A_k
\coloneqq 
\sum_{m=k}^\infty d_{m,\gamma}^{\frac{1}{1+\gamma}}$ 
and
$B_k = \sum_{m=k}^\infty 
d_{m,2}^{1/2}\max\{|\log d_{m,2}|,1\}^{1/2}$ with $d_{m,\gamma}
=\sum_{j=m+1}^\infty|b^\gamma|^{j}=\frac{|b^\gamma|^{m+1}}{1-|b^\gamma|}$.
The term $A_k$ is evaluated as
\begin{align*}
A_k
=
\frac{1}{(1-|b^\gamma|)^\frac{\gamma}{1+\gamma}}
\sum_{m=k}^\infty 
{{|b^{\frac{\gamma}{1+\gamma}}|^{m+1}}}
=
\frac{1}{(1-|b^\gamma|)^\frac{\gamma}{1+\gamma}}
\frac{|b^{\frac{\gamma}{1+\gamma}}|^k}{1-|b^{\frac{\gamma}{1+\gamma}}|}
.
\end{align*}

For the second term $B_k$, note that for sufficiently large $m$, $|\log d_{m,2}|=\left|\log \frac{|b^2|^{m+1}}{1-|b^2|}\right|>1$ and there is some constant $C>0$ such that $|\log d_{m,2}| \leq C m$.
Thus,
\begin{align*}
B_k = \sum_{m=k}^\infty 
(d_{m,2}|\log d_{m,2}|)^{1/2}
\leq &
(\frac{b^2 C}{1-|b^2|})^{1/2}
\sum_{m=k}^\infty 
m^{1/2}|b|^{m}\\
\leq &
(\frac{b^2 C}{1-|b^2|})^{1/2}
\sum_{m=k}^\infty 
m|b|^{m}\\
\leq &
(\frac{b^2 C}{1-|b^2|})^{1/2}
\frac{|b|^k (k - (k-1)|b|)}{(1-|b|)^2}.
\end{align*}
Thus, for sufficiently large $k$, there exists a constant $C'$ such that $\max\{A_k,B_k\}\leq C'|b^{\frac{\gamma}{1+\gamma}}|^k$.

Therefore, by \citet[Theorem 4.4.1]{r00}, we conclude that for sufficiently large $k$,
$\alpha(k) \leq C'' |b^{\frac{\gamma}{1+\gamma}}|^{k/2}$.
\qed

\subsection{Proof of Lemma \ref{IdentifiableGMA}}
Suppose $f_{\bm \theta}(\lambda;u,v)
=f_{\bm \theta'}(\lambda;u,v)$ for $\bm \theta=(a,\sigma^2)$ and $\bm \theta'=(a',{\sigma'}^2)$.
This implies
the equality of the corresponding Fourier coefficients: \begin{align}\nonumber
&\exp\skakko{-\frac{\sigma^2}{2}(u^2+v^2)(a^{-2}+1)}
\lkakko{
\exp\mkakko{-\sigma^2 uv (a^{-2}+1)}-1
}\\\label{eq:GMA_iden_1}
&=
\exp\skakko{-\frac{{\sigma'}^2}{2}(u^2+v^2)({a'}^{-2}+1)}
\lkakko{
\exp\mkakko{-{\sigma'}^2 uv ({a'}^{-2}+1)}-1
}
\end{align}
and
\begin{align}  \nonumber
&\exp\skakko{-\frac{\sigma^2}{2}(u^2+v^2)(a^{-2}+1)}
2\mkakko{\exp\skakko{\sigma^2 uv a^{-1}}-1}
\cos \lambda\\\label{eq:GMA_iden_2}  
&=
\exp\skakko{-\frac{{\sigma'}^2}{2}(u^2+v^2)({a'}^{-2}+1)}
2\mkakko{\exp\skakko{{\sigma'}^2 uv {a'}^{-1}}-1}
\cos \lambda.
\end{align}

Setting $v=u$, equations \eqref{eq:GMA_iden_1} and \eqref{eq:GMA_iden_2} reduce to
\begin{align*}
&A(u)\lkakko{
A(u)-1
}
=
A'(u)\lkakko{
A'(u)-1
}\\
&\text{and }
\quad
A(u)\mkakko{\exp\skakko{\sigma^2 u^2 a^{-1}}-1}
\cos \lambda
=
A'(u)
\mkakko{\exp\skakko{{\sigma'}^2 u^2 {a'}^{-1}}-1}
\cos \lambda,
\end{align*}
respectively, where
$A(u)=\exp\skakko{-{\sigma^2}u^2(a^{-2}+1)}$ and $A'(u)=\exp\skakko{-{\sigma'}^2u^2({a'}^{-2}+1)}
$.
The first equation implies $A(u)\equiv A'(u)$ or $A(u)+A'(u)\equiv1$.
Suppose $A(u)+A'(u)=1$ and $A(u)\neq A'(u)$.
Then, from $A(\sqrt 2) = A^2(1)$, $A^2(1)+{A'}^2(1)=1$ also holds true. Substituting $A'(1)=1 - A(1)$ into this equation yields
 $A(1)=0,1$.
This is impossible. Thus, $A(u)=A'(u)$, which implies ${\sigma}^2({a}^{-2}+1)={\sigma'}^2({a'}^{-2}+1)$.
From \eqref{eq:GMA_iden_2}, it holds that
$\sigma^2 a^{-1}
=
{\sigma'}^2 {a'}^{-1}$.
Since
$$\frac{{\sigma}^2({a}^{-2}+1)}{\sigma^2 a^{-1}}
=
\frac{{\sigma'}^2({a'}^{-2}+1)}{
{\sigma'}^2 {a'}^{-1}}
\quad
\Leftrightarrow
\quad
a+\frac{1}{a}=a'+\frac{1}{a'}.
$$
Since $x+1/x$ is strictly increasing for $x>1$, we deduce $a=a'$.
Consequently, $\sigma^2 = \sigma'^2$ holds.
\qed
\subsection{Proof of Lemma \ref{IdentifiableGAR}}
If we assume that $f_{{\bm \theta}}(\lambda, u,v)=f_{{\bm \theta}'}(\lambda, u,v)$, it implies the equality of the Fourier coefficients. I.e., for all $\ell \in \mathbb{Z}$, $C_{\ell,a,\sigma^2}(u,v)=C_{\ell,a',(\sigma')^2}(u,v)$. For $\ell=0, u=v$, this identity yields $$\exp\skakko{-\frac{2u^2\sigma^2}{1-a^2}}
-\exp\skakko{-\frac{u^2\sigma^2}{1-a^2}}
=
\exp\skakko{-\frac{2u^2(\sigma')^2}{1-(a')^2}}
-\exp\skakko{-u^2\frac{(\sigma')^2}{1-(a')^2}}.$$ Taking the derivative with respect to $u^2$ on both sides yields
\begin{align*}
&\frac{2\sigma^2}{1-a^2} \exp\skakko{-\frac{2u^2\sigma^2}{1-a^2}}
-\frac{\sigma^2}{1-a^2} \exp\skakko{-\frac{u^2\sigma^2}{1-a^2}}\\
=&
\frac{2(\sigma')^2}{1-(a')^2}\exp\skakko{-\frac{2u^2(\sigma')^2}{1-(a')^2}}
-\frac{(\sigma')^2}{1-(a')^2}\exp\skakko{-u^2\frac{(\sigma')^2}{1-(a')^2}}.
\end{align*}
Evaluating now this derivative in $u^2=0$ directly implies that 
\begin{align}\label{UsefulThing}
\frac{\sigma^2}{1-a^2}=\frac{(\sigma')^2}{1-(a')^2}.
\end{align}  

Now, for $\ell=1, u=v=1$, equality of the Fourier coefficient yields $$\exp\skakko{-\frac{2a\sigma^2}{1-a^2}}
-\exp\skakko{-\frac{\sigma^2}{1-a^2}}
=
\exp\skakko{-\frac{2a'(\sigma')^2}{1-(a')^2}}
-
\exp\skakko{-\frac{(\sigma')^2}{1-(a')^2}}.$$ Using \eqref{UsefulThing} we get that for $\ell=1$ and $u=v=1$, $$-\exp\skakko{-\frac{2\sigma^2a}{1-a^2}}
=
-\exp\skakko{-\frac{2(\sigma')^2a'}{1-(a')^2}},$$ which directly yields $a=a'$, using again \eqref{UsefulThing}. Using now $a=a'$ in \eqref{UsefulThing} yields $\sigma^2=(\sigma')^2$.
\qed

\subsection{Proof of Lemma \ref{IdentifiableINMA}}

First, note that the following equalities hold.

\begin{align*}
f(\lambda&;u,-u)\\
\coloneqq &
-\frac{1}{2\pi}
\exp\skakko{-\delta
\lkakko{(p^\alpha+1)
\mkakko{2(1-e^{iu})^\alpha}}}(2\cos\lambda+1)\\
&+
\frac{1}{2\pi}\exp\skakko{-\delta(p^\alpha+1)(1-e^{i2u})^\alpha}\\
&+\frac{1}{2\pi}\exp\skakko{-\delta \lkakko{(p^\alpha+1)(1-e^{iu})^\alpha
+\mkakko{1- (e^{iu}-pe^{iu}+pe^{i2u})}^\alpha}}2\cos\lambda.
\end{align*}
This readily implies that the following equalities hold:


%

%
%

\begin{align*}
f(\pi/2;\pi,-\pi)
\coloneqq 
&-\frac{1}{2\pi}
\exp\skakko{-\delta2^{\alpha+1}(p^\alpha+1)}
+\frac{1}{2\pi},\\
\quad
f(\pi/2;\pi/2,-\pi/2)
\coloneqq &
-\frac{1}{2\pi}
\exp\skakko{-2\delta
(p^\alpha+1)
(1-i)^\alpha}
+\frac{1}{2\pi}e^{-\delta2^\alpha(p^\alpha+1)},\\
f(2\pi/3;\pi,-\pi)
\coloneqq &
\frac{1}{2\pi}+\frac{1}{2\pi}\exp\skakko{-\delta \lkakko{(p^\alpha+1)2^\alpha
+(2- 2p)^\alpha}}.
\end{align*}
This directly yields identifiability. Indeed, consider ${\pmb \theta}\coloneqq (\delta, \alpha, p) \ne (\delta', \alpha', p')\coloneqq {\pmb \theta}'$ but assume that their respective spectral densities are such that $f_{{\pmb \theta}'}(\lambda;u,v)=f_{\pmb \theta}(\lambda;u,v)$ for all $(\lambda, u, v)$. Now, combining $f(\pi/2;\pi/2,-\pi/2)$ and $f(\pi/2;\pi,-\pi)$, we get that $\alpha=\alpha'$. Then, we consider $f(\pi/2;\pi,-\pi)$ and $f(2\pi/3;\pi,-\pi)$, we get that $(1-p')^{\alpha}p^{\alpha} -(1-p)^{\alpha}(p')^{\alpha}=(1-p)^{\alpha}-(1-p')^{\alpha}$. Assume $p < p'$ and deduce a contradiction. Identifiability of $\delta$ follows immediately from any one of the above displays.

\subsection{Derivation of (3)}
Since 
$Z_t= \sum_{j=0}^\infty a^j \epsilon_{t-j}$
and 
$Z_{t+\ell}
= a^{\ell}Z_{t}
+ \sum_{j=0}^{\ell-1} a^{j}\epsilon_{t+\ell-j}$  for $\ell\geq0$,
we find
$$
{\rm E}\skakko{e^{iuZ_t}}
=
\prod_{j=0}^\infty{\rm E}\skakko{e^{iua^j\epsilon_{t-j}}}
=
\prod_{j=0}^\infty
e^{-\delta|ua^j|}
=
e^{-\frac{\delta|u|}{1-|a|}}$$
and 
\begin{align*}
{\rm E}\skakko{e^{iuZ_{t+\ell}+ivZ_{t}}}
=&
{\rm E}\skakko{e^{i(ua^{\ell}+v)Z_{t}}}
{\rm E}\skakko{e^{iu\sum_{j=0}^{\ell-1}a^{j}\epsilon_{t+\ell-j}}}\\
=&
e^{-\frac{\delta|ua^{\ell}+v|}{1-|a|}}
e^{-\frac{\delta|u|(1-|a|^\ell)}{1-|a|}}\\
=&
e^{-\frac{\delta}{1-|a|}
\skakko{
|ua^{\ell}+v|
+
|u|(1-|a|^\ell)
}}.
\end{align*}
\qed

\subsection{Derivation of (4)}
Since 
$Z_{t-\ell}
= a^{-\ell}Z_{t}
- \sum_{j=1}^{\ell} a^{-j}\epsilon_{t-\ell+j}$ for $\ell\geq0$
and
$Z_t= -\sum_{j=1}^\infty a^{-j} \epsilon_{t+j}$
,
we find
\begin{align*}
{\rm E}\skakko{e^{iuZ_t}}
=&
\prod_{j=1}^\infty{\rm E}\skakko{e^{-iua^{-j}\epsilon_{t+j}}}\\
=&
\prod_{j=1}^\infty
e^{-\delta|ua^{-j}|}\\
=&
e^{-\frac{\delta|u| |a^{-1}|}{1-|a^{-1}|}}
\end{align*}
and 
\begin{align*}
{\rm E}\skakko{e^{iuZ_{t-\ell}+ivZ_{t}}}
=&
{\rm E}\skakko{e^{i(ua^{-\ell}+v)Z_{t}}}
{\rm E}\skakko{e^{-iu\sum_{j=1}^{\ell}a^{j}\epsilon_{t-\ell+j}}}\\
=&
e^{-\frac{\delta|ua^{-\ell}+v||a^{-1}|}{1-|a^{-1}|}}
e^{-\frac{\delta|u|(1-|a|^{-\ell})|a^{-1}|}{1-|a^{-1}|}}\\
=&
e^{-\frac{\delta|a^{-1}|}{1-|a^{-1}|}
\skakko{
|ua^{-\ell}+v|
+
|u|(1-|a^{-\ell}|)
}}.
\end{align*}\qed

\subsection{Proofs of Equations (5) and (6)}
\label{Proof_pgf_INAR}
For $u\in[0,1]$,  it holds that
\begin{align*}
&{\rm E}\skakko{u^{X_t}}\\
=&
{\rm E}\skakko{u^{\epsilon_t}}
{\rm E}\skakko{u^{p\circ X_{t-1}}}\\
=&
{\rm E}\skakko{u^{\epsilon_t}}
{\rm E}\skakko{{\rm E}\skakko{u^{p\circ X_{t-1}}\mid \mathcal F_{t-1}}}\\
=&
\exp\mkakko{-\delta(1-u)^\alpha}
{\rm E}\skakko{(1-p+pu)^{X_{t-1}}}\\
=&
\exp\mkakko{-\delta(1-u)^\alpha}
\exp\mkakko{-\delta (1-u)^\alpha p^{\alpha}}
{\rm E}\skakko{(1-p^2+p^2u)^{X_{t-2}}}
\\
=&
\exp\mkakko{-\delta(1-u)^\alpha}
\exp\mkakko{-\delta (1-u)^\alpha p^{\alpha}}
\exp\mkakko{-\delta (1-u)^\alpha p^{2\alpha}}
{\rm E}\skakko{(1-p^3+p^3u)^{X_{t-3}}}
\\
=&
\exp\mkakko{-\delta(1-u)^\alpha \sum_{h=0}^{k-1} p^{h\alpha}}
{\rm E}\skakko{(1-p^k+p^ku)^{X_{t-k}}}.
\end{align*}
By the stationarity of $X_t$ and the dominated convergence theorem, we obtain
$${\rm E}\skakko{u^{X_t}}
=\exp\mkakko{-\delta(1-u)^\alpha \sum_{h=0}^{\infty} p^{h\alpha}}
=\exp\mkakko{- \delta\frac{(1-u)^\alpha}{1-p^\alpha}}
.$$
For $\ell\in\{1,2,\ldots\}$, we observe
\begin{align*}
&{\rm E}\skakko{u^{X_{t+\ell}}v^{X_{t}}}\\
=&
{\rm E}\skakko{{\rm E}\skakko{u^{X_{t+\ell}}\mid \mathcal F_{t+\ell-1}}v^{X_{t}}}\\
=&
{\rm E}\skakko{{\rm E}\skakko{u^{p\circ X_{t+\ell-1}}\mid \mathcal F_{t+\ell-1}}v^{X_{t}}}
{\rm E}\skakko{u^{\epsilon_{t+\ell}}}\\
=&
{\rm E}\skakko{(1-p(1-u))^{X_{t+\ell-1}}v^{X_{t}}}
{\rm E}\skakko{u^{\epsilon_{t+\ell}}}\\
=&
{\rm E}\skakko{(1-p^2(1-u))^{X_{t+\ell-2}}v^{X_{t}}}
{\rm E}\skakko{u^{\epsilon_{t+\ell}}}
{\rm E}\skakko{(1-p(1-u))^{\epsilon_{t+\ell-1}}}\\
=&
{\rm E}\skakko{\skakko{v(1-p^\ell(1-u))}^{X_{t}}}
\prod_{h=0}^{\ell-1}
{\rm E}\skakko{(1-p^h(1-u))^{\epsilon_{t+\ell-h}}}\\
=&
\exp\mkakko{-\delta \frac{(1-v(1-p^\ell(1-u)))^\alpha}{1-p^\alpha}}
\prod_{h=0}^{\ell-1}
\exp\mkakko{-\delta p^{h\alpha}(1-u)^\alpha}\\
=&
\exp\mkakko{-\delta \frac{(1-v(1-p^\ell(1-u)))^\alpha}{1-p^\alpha}}
\exp\mkakko{-\delta (1-u)^\alpha \sum_{h=0}^{\ell-1} p^{h\alpha}}\\
=&
\exp\mkakko{-\delta \frac{(1-v(1-p^\ell(1-u)))^\alpha}{1-p^\alpha}
-\delta (1-u)^\alpha\frac{1-p^{\alpha\ell}}{1-p^\alpha}}.
\end{align*} \qed

\subsection{Proof of Lemma \ref{lem:stationarity}}
Let $G_t(\cdot)$ denote the pgf of $Z_t$.
For $x\in \mathbb C$ such that $|x|\leq 1$, we have
$
{\rm E}\skakko{x^{Z_t}\mid Z_{t-1}}
=
\skakko{1-p(1-x)}^{Z_{t-1}}
{\rm E}\skakko{x^{\epsilon_t}}$. Taking expectations yields
\begin{align*}
G_t\skakko{x}
=&
G_{t-1}\skakko{1-p(1-x)}
\exp\{-\delta(1-x)^\alpha\}\\
=&
G_0\skakko{1-p^t(1-x)}
\prod_{k=0}^{t-1}
\exp\mkakko{-\delta p^{k\alpha}(1-x)^\alpha}.
\end{align*}
Since $G_0$ is continuous on $|x|\le 1$ and $G_0(1)=1$, we obtain 
$G_t(x)\to G_{\infty}(x)\coloneqq \exp\{-\delta(1-x)^\alpha/(1-p^\alpha)\}$
as $t\to\infty$.

Now suppose $Z_{t-1}$ has the pgf $G_{\infty}(\cdot)$. Then
\begin{align*}
G_t\skakko{x}
=&
G_{\infty}\skakko{1-p(1-x)}
\exp\mkakko{-\delta(1-x)^\alpha}\\
=&
\exp\mkakko{-\frac{\delta p^\alpha(1-x)^\alpha}{1-p^\alpha}
-\delta(1-x)^\alpha}
=
G_{\infty}\skakko{x}.
\end{align*}
Therefore $G_\infty(\cdot)$ is invariant, and choosing $Z_0\sim G_\infty$ yields a strictly stationary solution.

To show the uniqueness, we suppose that $Z_t$ is a stationary solution with a pgf $G$.
Then, we have
\begin{align*}
G\skakko{x}
=&
G\skakko{1-p(1-x)}
\exp\mkakko{-\delta(1-x)^\alpha}\\
=&
G\skakko{1-p^t(1-x)}
\prod_{k=0}^{t-1}
\exp\mkakko{-\delta p^{k\alpha}(1-x)^\alpha}\\
=&
\exp\mkakko{-\frac{\delta(1-x)^\alpha}{1-p^\alpha}}\\
=&
G_{\infty}(x).
\end{align*}
Hence the strictly stationary distribution is unique.

Irreducibility and aperiodicity of the process follow directly from the property that the innovation $\epsilon_t$ has a full support on $\mathbb{N}_0$, ensuring $P(Z_t = j \mid Z_{t-1} = i) > 0$ for all $i, j \in \mathbb{N}_0$. 
The existence of a unique invariant pgf $G_{\infty}(x)$ satisfies $G_{\infty}(1)=1$, which indicates that the invariant measure is a proper probability distribution.   
This establishes the positive recurrence of the Markov chain (see \citealt[Theorem 21.13]{lp08}). 
Consequently, the unique strictly stationary solution is ergodic (\citealt[Theorem 21.16]{lp08}).

Next, we show the process enjoys the geometrically $\beta$-mixing property.
The $\beta$-mixing coefficient is defined as
\begin{align*}
\beta(h)\coloneqq &
{\rm}E\skakko{
\sup_{\mathbf A\in\mathcal B^\infty}
\left|
{\rm P}\skakko{(Z_{t+h},Z_{t+h+1},\ldots)\in \mathbf A \mid Z_t}
-
{\rm P}\skakko{(Z_{t+h},Z_{t+h+1},\ldots)\in \mathbf A}
\right|
}\\
=&
{\rm E}\skakko{
\left\|
Q^{h}_{Z_t}
-
Q
\right\|_{\rm TV}
},
\end{align*}
where $Q^h_{Z_t}(\cdot)$ denotes the conditional distribution of
$(Z_{t+h},Z_{t+h+1},\ldots)$ given $Z_t$, $Q(\cdot)$ denotes its unconditional distribution.
$\|\cdot\|_{\rm TV}$ is the total variation distance, defined as, for probability measures $Q$ and $Q^\prime$ on a measurable space $((\mathbb N\cup\{0\})^\infty, \mathcal B^\infty)$, 
$\|Q-Q^\prime\|_{\rm TV}
\coloneqq 
\sup_{\mathbf A\in \mathcal B^\infty}
|Q(\mathbf A)-Q^\prime(\mathbf A)|
$.

For the Markov kernel $K_{Z_0}(\bm A)\coloneqq {\rm Q}\skakko{(Z_0,Z_1,\ldots)\in\bm A\mid Z_0}$ for $\bm A\in\mathcal B^\infty$,
the h-step kernel $P^h_{Z_t}(z,A)\coloneqq {\rm P}\skakko{Z_{t+h}\in A |Z_t}$ for $A\subset\mathbb N\cup \{0\}$,
the invariant measure $\pi$ of $Z_t$ under stationarity, 
 the tower and Markov properties yield 
$$Q^h_{Z_t}(\mathbf A)
=\int K_{Z_0}(\bm A)P^h_{Z_t}(dx)
\quad
\text{and}
\quad
{\rm Q}(\bm A)
=\int K_{Z_0}(\bm A) \pi(dx).
$$
Therefore,
it holds that
$
{\rm E}\skakko{
\left\|
Q^{h}_{Z_t}
-
Q
\right\|_{\rm TV}
}
\leq
{\rm E}\skakko{
\left\|
P^h_{Z_t}
-
\pi
\right\|_{\rm TV}
}.
$

Let $Z_t^*$ be a random variable, independent of $Z_t$, with distribution $\pi$ and define
$Z_{t+h}^* = p \circ Z_{t+h-1}^* + \epsilon_{t+h}$ for any $h\in \mathbb N$.
Note that we assume  the counting sequences in the thinning operator and the error terms are identical for both $Z_{t+h}$ and $Z_{t+h}^*$.
Then, we observe, for any $A\subset\mathbb N\cup \{0\}$, 
\begin{align*}
\left|
P^h_{Z_t}(A)-\pi(A)
\right|
=
\left|
{\rm P}\skakko{Z_{t+h}\in A|Z_t}
-
{\rm P}\skakko{Z_{t+h}^*\in A|Z_t}
\right|
\leq
{\rm P}\skakko{Z_{t+h}\neq Z_{t+h}^*|Z_t}
\end{align*}
and 
\begin{align*}
{\rm E}\skakko{{\rm P}\skakko{Z_{t+h}\neq Z_{t+h}^*|Z_t}}
=
{\rm E}\skakko{{\rm P}\skakko{Z_{t+h}\neq Z_{t+h}^*|Z_t,Z_t^*}}
=
{\rm E}\skakko{{\rm P}\skakko{D_{t+h}>0|Z_t,Z_t^*}},
\end{align*}
where $D_{t}=|Z_{t}-Z_{t}^*|$. 
Since
$D_{t+h}$
follows conditionally on $D_{t}$ binomial distribution with number of trials $D_t$ and success probability $p^h$.
It holds that, for $r\in(0,\min\{\alpha,1\})$,
\begin{align*}
{\rm E}\skakko{{\rm P}\skakko{D_{t+h}>0|Z_t,Z_t^*}}
=&
{\rm E}\skakko{
1- (1-p^h)^{D_t}
}\\
\leq&
{\rm E}\skakko{
\min\{D_tp^h,1\}
}\\
\leq&
p^{rh}{\rm E}(D_{t}^r)\\
\leq&
2p^{rh}{\rm E}(Z_{t}^r).
\end{align*}
Here, we used Bernoulli's inequality.
Thus, we conclude that there exist constants $C>0$ and $\varrho\in(0,1)$ such that $\beta(h)\leq C\varrho^h$.
\qed

\subsection{Proof of Lemma \ref{IdentifiableINAR}} Let $f_{\bm{\theta}}(\lambda;u,v)\coloneqq f(\lambda;u,v)$ with
$\bm{\theta}=(p,\delta,\alpha)$.
Suppose that
$f_{\bm{\theta}}(\lambda;u,v)=f_{\bm{\theta}'}(\lambda;u,v)
$
for all $\lambda\in[0,2\pi]$ and $(u,v)\in[-L,L]^2$, where
$\bm{\theta}'=(p',\delta',\alpha')$.
This implies
$
C_{\ell,\bm{\theta}}(u,v)
=
C_{\ell,\bm{\theta}'}(u,v)
$ for all $\ell\in\mathbb Z$ and all $(u,v)\in[-L,L]^2$.
Define 
$S_{\bm{\theta}}^{(1)}(u)
\coloneqq 
\mathrm{E}\!\left[u^{X_t}\right]$ and $S_{\ell,\bm{\theta}}^{(2)}(u,v)
\coloneqq 
\mathrm{E}\!\left[u^{X_{t+\ell}}v^{X_t}\right]$.
Recall that these quantities admit the explicit representations given in
\eqref{eq:pgf_INAR1} and \eqref{eq:pgf_INAR2}.

\textbf{Step 1: Identification of $\alpha$}.
For $\ell=0$, we observe that
\begin{align*}
C_{0,\bm{\theta}}(u,-u)
&=
S_{0,\bm{\theta}}^{(2)}(e^{iu},e^{-iu})
-
\left|S_{\bm{\theta}}^{(1)}(e^{iu})\right|^2
=
1-
\left|S_{\bm{\theta}}^{(1)}(e^{iu})\right|^2.
\end{align*}
Hence,
\begin{equation}\label{eq:S1}
\left|S_{\bm{\theta}}^{(1)}(e^{iu})\right|
=
\sqrt{1-C_{0,\bm{\theta}}(u,-u)}.
\end{equation}
By \eqref{eq:S1}, we obtain
\[
g_{1,\bm{\theta}}
\coloneqq 
\frac{\log
\sqrt{1-C_{0,\bm{\theta}}(\pi/2,-\pi/2)}
}{\log
\sqrt{1-C_{0,\bm{\theta}}(\pi,-\pi)}
}
=
\frac{
\log\left|S_{\bm{\theta}}^{(1)}(e^{i\pi/2})\right|
}{
\log\left|S_{\bm{\theta}}^{(1)}(e^{i\pi})\right|
}
=
2^{-\alpha/2}
\cos\!\left(\frac{\alpha\pi}{4}\right).
\]
Therefore,
$g_{1,\bm{\theta}}=g_{1,\bm{\theta}'}$
implies
$\alpha=\alpha'$.

\textbf{Step 2: Identification of $p$.}
From \eqref{eq:S1}, we have
\[
C_{1,\bm{\theta}}(u,-u)
=
S_{1,\bm{\theta}}^{(2)}(e^{iu},e^{-iu})
-
1
+
C_{0,\bm{\theta}}(u,-u),
\]
which yields
\[
S_{1,\bm{\theta}}^{(2)}(e^{iu},e^{-iu})
=
1
-
C_{0,\bm{\theta}}(u,-u)
+
C_{1,\bm{\theta}}(u,-u).
\]
Then,
\[
g_{2,\bm{\theta}}
\coloneqq 
\frac{
\log\!\left(
1
-
C_{0,\bm{\theta}}(u,-u)
+
C_{1,\bm{\theta}}(u,-u)
\right)
}{
\log\sqrt{1-C_{0,\bm{\theta}}(u,-u)}
}
=
\frac{
\log S_{1,\bm{\theta}}^{(2)}(e^{iu},e^{-iu})
}{
\log \left|S_{\bm{\theta}}^{(1)}(e^{iu})\right|
}
=
1-p^\alpha+(1-p)^\alpha.
\]
Since $\alpha=\alpha'$ has already been identified, the equality
$g_{2,\bm{\theta}}=g_{2,\bm{\theta}'}$
implies
$p=p'$.

\textbf{Step 3: Identification of $\delta$.}
Finally, define
\[
g_{3,\bm{\theta}}
\coloneqq 
\log\sqrt{1-C_{0,\bm{\theta}}(\pi,-\pi)}
=
\log\left|S_{\bm{\theta}}^{(1)}(e^{i\pi})\right|
=
-\frac{\delta\,2^\alpha}{1-p^\alpha}.
\]
Since $p=p'$ and $\alpha=\alpha'$, the equality
$g_{3,\bm{\theta}}=g_{3,\bm{\theta}'}$
implies
$\delta=\delta'$.

Therefore, the parameter vector $\bm{\theta}$ is identifiable.

\subsection{Preliminary lemma} 
We provide a useful lemma, providing a closed-form expression for the joint cumulants of $d_n(\lambda;u_1)$.

\begin{lem}\label{lem_cum} Let $\{Z_t\}$ be a strictly stationary process satisfying the condition \eqref{as_sum} with $q=1$. For any $k\in \mathbb N\backslash \{1\}$, it holds that 
\begin{align*}
&{\rm Cum}\skakko{d_n(\lambda_1;u_1),\ldots,d_n(\lambda_k;u_k)}\\
=&
(2\pi)^{k-1}\Delta_n\skakko{\sum_{j=1}^k\lambda_j}
f(\lambda_1,\ldots,\lambda_{k-1};u_1,\dots,u_k)
+O(1),
\end{align*}
where $\Delta_n(\lambda)\coloneqq \sum_{t=1}^n\exp\{-i\lambda t\}$ and the error term is uniform on $(\lambda_1,\ldots,\lambda_{k-1})\in[0,2\pi]^{k-1}$ and $(u_1,\ldots,u_k)\in\mathbb R^k$. 
For $k=1$, it holds that
\begin{align*}
{\rm E}d_n(\lambda;u)=
\begin{cases}
0&\lambda=\frac{2\pi j}{n} \text{ for $j=1,\ldots,n-1$},\\
n{\rm E}e^{iuZ_t}&\lambda=2\pi k\text{ for $k\in\mathbb Z$},\\
\frac{1-e^{-in\lambda}}{1-e^{-i\lambda}}{\rm E}e^{iuZ_t}&\text{otherwise}.\\
\end{cases}
\end{align*}
\end{lem}
{ Note that $f(\lambda_1,\ldots,\lambda_{k-1};u_1,\dots,u_k)$, the generalized spectrum of order $k$ is defined in Section $2$.} This lemma can be shown along the same lines of Theorem 4.3.2 in \cite{b81}.

\subsection{Proof of Lemma \ref{lem1}}
It suffices to show that
\begin{align}\label{ex_D}
{\rm E}D_n(I_n,f_{\bm{\theta}})-\tilde D(f,f_{\bm{\theta}})
=O\skakko{
\skakko{n^{-\min\{\alpha, \tilde{\alpha}\}}+M_n^{-\min\{\beta, \tilde{\beta}\}}}^{{1}/{(\tau+1)}}
+M_n^{-1}+n^{-1}
}
\end{align}
and
\begin{align}\label{var_D}
{\rm Var}\skakko{D_n(I_n,f_{\bm{\theta}})}=O(1/n)
\end{align}
since, for any $\epsilon>0$,
\begin{align*}
&{\rm P}\skakko{\left|D_n(I_n,f_{\bm{\theta}})-\tilde D(f,f_{\bm{\theta}})\right|>\epsilon}\\
\leq&
\frac{1}{\epsilon}
{\rm E}\skakko{\left|D_n(I_n,f_{\bm{\theta}})-{\rm E}D_n(I_n,f_{\bm{\theta}})\right|}
+\frac{1}{\epsilon}\left|{\rm E}D_n(I_n,f_{\bm{\theta}})-\tilde D(f,f_{\bm{\theta}})\right|\\
\leq&
\frac{1}{\epsilon}\sqrt{{\rm Var}\skakko{D_n(I_n,f_{\bm{\theta}})}}
+\frac{1}{\epsilon}\left|{\rm E}D_n(I_n,f_{\bm{\theta}})-\tilde D(f,f_{\bm{\theta}})\right|,
\end{align*}
which tends to zero as $n\to\infty$.

First, we prove \eqref{ex_D}.
{ Recalling that $\lambda_j={2\pi\,j}/{n}$ and that $f(\lambda, u,v)$ is bounded, Lemma \ref{lem_cum} {alongside Theorem 2.3.2 in \cite{b81}} immediately gives that
\begin{align*}
{\rm E} I_n(\lambda_j;u_{i_1},v_{i_2})=f(\lambda_j;u_{i_1},v_{i_2})+O\skakko{{1}/{n}}
\end{align*}
and
\begin{align*}
&{\rm E}\skakko{I_n(\lambda_j;u_{i_1},v_{i_2}) \overline{I_n(\lambda_j;u_{i_1},v_{i_2})}}\\
=&f(\lambda_j;u_{i_1},v_{i_2})f(-\lambda_j;-u_{i_1},-v_{i_2})
+f(\lambda_j;u_{i_1},-u_{i_1})f(-\lambda_j;v_{i_2},-v_{i_2})\\
&+f(\lambda_j;u_{i_1},-v_{i_2})f(-\lambda_j;v_{i_2},-u_{i_1})\Delta_n\skakko{2\lambda_j}\Delta_n\skakko{-2\lambda_j}/n^2
+O\skakko{{1}/{n}},
\end{align*}
where the error term is uniform on $(\lambda_1,\ldots,\lambda_{k-1})\in[0,2\pi]^{k-1}$ and $(u_1,\ldots,u_k)\in\mathbb R^k$. 

Recall that the difference between an integral and the corresponding Riemann sum
for an $\alpha$-H\"older continuous function is $O(L_{\rm H}K^{-\alpha})$, where $K$ is the number of subintervals
and $L_{\rm H}$ denotes the H\"older constant associated with the function (see, e.g., \citealt[Lemma P5.1, p.~415]{b81}).

Let
$u_i\coloneqq -L+\frac{2Li}{M_n}$ and
$v_i\coloneqq -L+\frac{2Li}{M_n}$ for  $i=1,\ldots,M_n$.
Since we suppose ${\rm Leb}(\mathcal S_\varepsilon)\coloneqq O(\varepsilon)$ as $\varepsilon\downarrow 0$, 
$(u_{i_1},v_{i_2})\in \mathcal S_\varepsilon$ implies  
$$\lkakko{u_{i_1}-\frac{L}{M_n},u_{i_1}+\frac{L}{M_n}}\times\lkakko{v_{i_2}-\frac{L}{M_n},v_{i_2}+\frac{L}{M_n}}\subset \mathcal S_{\varepsilon+\frac{\sqrt 2 L}{M_n}}.$$
So, $\#\{(i_1,i_2);(u_{i_1},v_{i_2})\in \mathcal S_\varepsilon\}\frac{4L^2}{M_n^2}\leq {\rm Leb}\skakko{\mathcal S_{\varepsilon+\frac{\sqrt 2 L}{M_n}}}$ and hence 
$\#\{(i_1,i_2);(u_{i_1},v_{i_2})\in \mathcal S_\varepsilon\}
=O(M_n^2\varepsilon +M_n)$.

Since, by the above argument, 
\begin{align*}
&\frac{8\pi L^2}{nM_n^2}\sum_{j=1}^{n-1}\sum_{i_1,i_2=1}^{M_n}
I_{\{(i_1,i_2);(u_{i_1},v_{i_2})\notin \mathcal S_\varepsilon\}}
f(\lambda_j;u_{i_1},v_{i_2})f(-\lambda_j;-u_{i_1},-v_{i_2})\\
=&
\int_{[0,2\pi]}
\int_{[-L,L]^2\backslash \mathcal S_\varepsilon}
f(\lambda;u,v)
f(-\lambda;-u,-v)
{\rm d}u{\rm d}v{\rm d}{\lambda}
+O\skakko{\frac{\varepsilon^{-\tilde{\tau}}}{n^{\tilde\alpha}}+\frac{\varepsilon^{-\tilde{\tau}}}{M_n^{\tilde\beta}}}\\
=&
\int_{[0,2\pi]}
\int_{[-L,L]^2}
f(\lambda;u,v)
f(-\lambda;-u,-v)
{\rm d}u{\rm d}v{\rm d}{\lambda}
+O\skakko{\frac{\varepsilon^{-\tilde{\tau}}}{n^{\tilde\alpha}}+\frac{\varepsilon^{-\tilde{\tau}}}{M_n^{\tilde\beta}}+\varepsilon}
\end{align*}
and 
\begin{align*}
&\frac{8\pi L^2}{nM_n^2}\sum_{j=1}^{n-1}\sum_{i_1,i_2=1}^{M_n}
I_{\{(i_1,i_2);(u_{i_1},v_{i_2})\in S_\varepsilon\}}
f(\lambda_j;u_{i_1},v_{i_2})f(-\lambda_j;-u_{i_1},-v_{i_2})
=
O\skakko{\varepsilon +\frac{1}{M_n}},
\end{align*}
we find that
\begin{align*}
&\frac{8\pi L^2}{nM_n^2}\sum_{j=1}^{n-1}\sum_{i_1,i_2=1}^{M_n}
f(\lambda_j;u_{i_1},v_{i_2})f(-\lambda_j;-u_{i_1},-v_{i_2})\\
=&
\int_{[0,2\pi]}
\int_{[-L,L]^2}
f(\lambda;u,v)
f(-\lambda;-u,-v)
{\rm d}u{\rm d}v{\rm d}{\lambda}
+O\skakko{
\frac{\varepsilon^{-\tilde{\tau}}}{n^{\tilde\alpha}}
+
\frac{\varepsilon^{-\tilde{\tau}}}{M_n^{\tilde\beta}}
+
\varepsilon +\frac{1}{M_n}}.
\end{align*}

Hence, as $n$ and $M_n \rightarrow \infty$ and $\epsilon \rightarrow 0$,
\begin{align*}
&{\rm E}\skakko{D_n(I_n,f_{\bm{\theta}})}\\
=&
\frac{8\pi L^2}{nM_n^2}\sum_{j=1}^{n-1}\sum_{i_1,i_2}^{M_n}{\rm E}| I_n(\lambda_j;u_{i_1},v_{i_2})-f_{\bm{\theta}}(\lambda_j;u_{i_1},v_{i_2}))|^2\\
=&\frac{8\pi L^2}{nM_n^2}\sum_{j=1}^{n-1}\sum_{i_1,i_2=1}^{M_n}
 {\rm E}\skakko{I_n(\lambda_j;u_{i_1},v_{i_2})\overline{I_n(\lambda_j;u_{i_1},v_{i_2})}}\\
&-\frac{8\pi L^2}{nM_n^2}\sum_{j=1}^{n-1}\sum_{i_1,i_2=1}^{M_n}
 {\rm E}\skakko{I_n(\lambda_j;u_{i_1},v_{i_2})}\overline{f_{\bm{\theta}}(\lambda_j;u_{i_1},v_{i_2})}\\
&- \frac{8\pi L^2}{nM_n^2}\sum_{j=1}^{n-1}\sum_{i_1,i_2=1}^{M_n}
f_{\bm{\theta}}(\lambda_j;u_{i_1},v_{i_2}){\rm E}\skakko{\overline{I_n(\lambda_j;u_{i_1},v_{i_2})}}\\
&+\frac{8\pi L^2}{nM_n^2}\sum_{j=1}^{n-1}\sum_{i_1,i_2=1}^{M_n}
f_{\bm{\theta}}(\lambda_j;u_{i_1},v_{i_2})\overline{f_{\bm{\theta}}(\lambda_j;u_{i_1},v_{i_2})}\\
=&\int_{0}^{2\pi}\int_{-L}^L\int_{-L}^L
f(\lambda;u,v)f(-\lambda;-u,-v)
+f(\lambda;u,-u)f(-\lambda;v,-v){\rm d}u{\rm d}v{\rm d}{\lambda}\\
&-\int_{0}^{2\pi}\int_{-L}^L\int_{-L}^L
f(\lambda;u,v)\overline{f_{\bm{\theta}}(\lambda;u,v)}{\rm d}u{\rm d}v{\rm d}{\lambda}\\
&-\int_{0}^{2\pi}\int_{-L}^L\int_{-L}^L
f_{\bm{\theta}}(\lambda;u,v)\overline{f(\lambda;u,v)}{\rm d}u{\rm d}v{\rm d}{\lambda}\\
&+\int_{0}^{2\pi}\int_{-L}^L\int_{-L}^L
f_{\bm{\theta}}(\lambda;u,v)\overline{f_{\bm{\theta}}(\lambda;u,v)}{\rm d}u{\rm d}v{\rm d}{\lambda}\\
&+O(\varepsilon^{-\max\{\tau, \tilde{\tau}\}} n^{-\min\{\alpha, \tilde{\alpha}\}}+\varepsilon^{-\max\{\tau, \tilde{\tau}\}} M_n^{-\min\{\beta, \tilde{\beta}\}} + \varepsilon +M_n^{-1}+n^{-1}
)\\
=&\tilde D(f,f_{\bm{\theta}})
+O(\varepsilon^{-\max\{\tau, \tilde{\tau}\}} n^{-\min\{\alpha, \tilde{\alpha}\}}+\varepsilon^{-\max\{\tau, \tilde{\tau}\}} M_n^{-\min\{\beta, \tilde{\beta}\}} + \varepsilon +M_n^{-1}+n^{-1}
).
\end{align*}
By choosing $\varepsilon=
\skakko{n^{-\min\{\alpha, \tilde{\alpha}\}}+M_n^{-\min\{\beta, \tilde{\beta}\}}}^{{1}/{(\max\{\tau, \tilde{\tau}\}+1)}}
$,
we obtain \eqref{ex_D}.

Next, we show \eqref{var_D}.
From the definition of $D_n(I_n,f_{\bm{\theta}})$, we know that
\begin{align*}
&{\rm Var}\skakko{D_n(I_n,f_{\bm{\theta}})}\\
=&{\rm E}\lkakko{\skakko{D_n(I_n,f_{\bm{\theta}})-{\rm E}D_n(I_n,f_{\bm{\theta}})}}^2\\
=&
{\rm E}\bigg[\bigg\{\frac{8\pi L^2}{nM_n^2}\sum_{j=1}^{n-1}\sum_{i_1,i_2=1}^{M_n}
\big(| I_n(\lambda_j;u_{i_1},v_{i_2})-f_{\bm{\theta}}(\lambda_j;u_{i_1},v_{i_2})|^2\\
&\qquad\qquad\qquad\qquad\qquad
-{\rm E} | I_n(\lambda_j;u_{i_1},v_{i_2})-f_{\bm{\theta}}(\lambda_j;u_{i_1},v_{i_2})|^2
\big)\bigg\}^2\bigg]\\
=&
\frac{8^2\pi^2 L^4}{n^2M_n^4}\sum_{j,j^\prime=1}^{n-1}\sum_{i_1,i_1^\prime,i_2,i_2^\prime=1}^{M_n}
{\rm E}\big\{
\big(| I_n(\lambda_j;u_{i_1},v_{i_2})-f_{\bm{\theta}}(\lambda_j;u_{i_1},v_{i_2})|^2\\
&\qquad\qquad\qquad\qquad\qquad\qquad\qquad
-{\rm E} | I_n(\lambda_j;u_{i_1},v_{i_2})-f_{\bm{\theta}}(\lambda_j;u_{i_1},v_{i_2})|^2\big))\\
&
\times\skakko{| I_n(\lambda_{j^\prime};u_{i_1^\prime},v_{i_2^\prime})-f_{\bm{\theta}}(\lambda_{j^\prime};u_{i_1^\prime},v_{i_2^\prime})|^2
-{\rm E} | I_n(\lambda_{j^\prime};u_{i_1^\prime},v_{i_2^\prime})-f_{\bm{\theta}}(\lambda_{j^\prime};u_{i_1^\prime},v_{i_2^\prime})|^2}\big\}\\
=&
\frac{8^2\pi^2 L^4}{n^2M_n^4}
\sum_{j,j^\prime=1}^{n-1}
\sum_{i_1,i_1^\prime,i_2,i_2^\prime=1}^{M_n}
{\rm Cum}\bigg(
| I_n(\lambda_j;u_{i_1},v_{i_2})-f_{\bm{\theta}}(\lambda_j;u_{i_1},v_{i_2})|^2,\\
&\qquad\qquad\qquad\qquad\qquad\quad\qquad
| I_n(\lambda_{j^\prime};u_{i_1^\prime},v_{i_2^\prime})-f_{\bm{\theta}}(\lambda_{j^\prime};u_{i_1^\prime},v_{i_2^\prime})|^2
\bigg).
\end{align*}
Simple algebra gives
\begin{align*}
&{\rm Cum}\left(
| I_n(\lambda_j;u_{i_1},v_{i_2})-f_{\bm{\theta}}(\lambda_j;u_{i_1},v_{i_2})|^2,
| I_n(\lambda_{j^\prime};u_{i_1^\prime},v_{i_2^\prime})-f_{\bm{\theta}}(\lambda_{j^\prime};u_{i_1^\prime},v_{i_2^\prime})|^2
\right)\\
=&{\rm Cum}\bigg(
I_n(\lambda_j;u_{i_1},v_{i_2})\overline{I_n(\lambda_j;u_{i_1},v_{i_2})}
-I_n(\lambda_j;u_{i_1},v_{i_2})\overline{f_{\bm{\theta}}(\lambda_j;u_{i_1},v_{i_2})}\\
&\qquad\qquad
-f_{\bm{\theta}}(\lambda_j;u_{i_1},v_{i_2})\overline{I_n(\lambda_j;u_{i_1},v_{i_2})},
I_n(\lambda_{j^\prime};u_{i_1^\prime},v_{i_2^\prime})\overline{I_n(\lambda_{j^\prime};u_{i_1^\prime},v_{i_2^\prime})}\\
&\qquad\qquad
-I_n(\lambda_{j^\prime};u_{i_1^\prime},v_{i_2^\prime})\overline{f_{\bm{\theta}}(\lambda_{j^\prime};u_{i_1^\prime},v_{i_2^\prime})}
-f_{\bm{\theta}}(\lambda_{j^\prime};u_{i_1^\prime},v_{i_2^\prime})\overline{I_n(\lambda_{j^\prime};u_{i_1^\prime},v_{i_2^\prime})}\bigg)\\
=&{\rm Cum}\left(I_n(\lambda_j;u_{i_1},v_{i_2})\overline{I_n(\lambda_j;u_{i_1},v_{i_2})},I_n(\lambda_{j^\prime};u_{i_1^\prime},v_{i_2^\prime})\overline{I_n(\lambda_{j^\prime};u_{i_1^\prime},v_{i_2^\prime})}\right)\\
&-{\rm Cum}\left(I_n(\lambda_j;u_{i_1},v_{i_2})\overline{I_n(\lambda_j;u_{i_1},v_{i_2})},I_n(\lambda_{j^\prime};u_{i_1^\prime},v_{i_2^\prime})\right)\overline{f_{\bm{\theta}}(\lambda_{j^\prime};u_{i_1^\prime},v_{i_2^\prime})}\\
&-{\rm Cum}\left(I_n(\lambda_j;u_{i_1},v_{i_2})\overline{I_n(\lambda_j;u_{i_1},v_{i_2})},\overline{I_n(\lambda_{j^\prime};u_{i_1^\prime},v_{i_2^\prime})}\right)f_{\bm{\theta}}(\lambda_{j^\prime};u_{i_1^\prime},v_{i_2^\prime})\\
&-{\rm Cum}\left(I_n(\lambda_j;u_{i_1},v_{i_2}),I_n(\lambda_{j^\prime};u_{i_1^\prime},v_{i_2^\prime})\overline{I_n(\lambda_{j^\prime};u_{i_1^\prime},v_{i_2^\prime})}\right)\overline{f_{\bm{\theta}}(\lambda_j;u_{i_1},v_{i_2})}\\
&+{\rm Cum}\left(I_n(\lambda_j;u_{i_1},v_{i_2}),I_n(\lambda_{j^\prime};u_{i_1^\prime},v_{i_2^\prime})\right)\overline{f_{\bm{\theta}}(\lambda_j;u_{i_1},v_{i_2})}\overline{f_{\bm{\theta}}(\lambda_{j^\prime};u_{i_1^\prime},v_{i_2^\prime})}\\
&+{\rm Cum}\left(I_n(\lambda_j;u_{i_1},v_{i_2}),\overline{I_n(\lambda_{j^\prime};u_{i_1^\prime},v_{i_2^\prime})}\right)\overline{f_{\bm{\theta}}(\lambda_j;u_{i_1},v_{i_2})}f_{\bm{\theta}}(\lambda_{j^\prime};u_{i_1^\prime},v_{i_2^\prime})\\
&-{\rm Cum}\left(\overline{I_n(\lambda_j;u_{i_1},v_{i_2})}
,I_n(\lambda_{j^\prime};u_{i_1^\prime},v_{i_2^\prime})\overline{I_n(\lambda_{j^\prime};u_{i_1^\prime},v_{i_2^\prime})}\right)f_{\bm{\theta}}(\lambda_j;u_{i_1},v_{i_2})\\
&+{\rm Cum}\left(\overline{I_n(\lambda_j;u_{i_1},v_{i_2})}
,I_n(\lambda_{j^\prime};u_{i_1^\prime},v_{i_2^\prime})\right)f_{\bm{\theta}}(\lambda_j;u_{i_1},v_{i_2})\overline{f_{\bm{\theta}}(\lambda_{j^\prime};u_{i_1^\prime},v_{i_2^\prime})}
\\
&+{\rm Cum}\left(\overline{I_n(\lambda_j;u_{i_1},v_{i_2})}
,\overline{I_n(\lambda_{j^\prime};u_{i_1^\prime},v_{i_2^\prime})}
\right)f_{\bm{\theta}}(\lambda_j;u_{i_1},v_{i_2})f_{\bm{\theta}}(\lambda_{j^\prime};u_{i_1^\prime},v_{i_2^\prime}).
\end{align*}
We shall only prove that
\begin{align}\nonumber
&\quad
\frac{8^2\pi^2 L^4}{n^2M_n^4}\sum_{j,j^\prime=1}^{n-1}\sum_{i_1,i_1^\prime,i_2,i_2^\prime=1}^{M_n}{\rm Cum}\bigg(I_n(\lambda_j;u_{i_1},v_{i_2})\overline{I_n(\lambda_j;u_{i_1},v_{i_2})},\\\label{var_case1}
&\qquad\qquad\qquad\qquad\qquad\qquad\qquad\qquad
I_n(\lambda_{j^\prime};u_{i_1^\prime},v_{i_2^\prime})\overline{I_n(\lambda_{j^\prime};u_{i_1^\prime},v_{i_2^\prime})}\bigg)
=O(1/n)
\end{align}
since the other terms can be evaluated in the same manner. Suppose that $(\nu_1,\ldots,\nu_p)$ { represents an} indecomposable partition of the table
\begin{align*}
\begin{pmatrix}
d(\lambda_j;u_{i_1}) &d(-\lambda_j;v_{i_2}) &d(-\lambda_j;-u_{i_1}) &d(\lambda_j;-v_{i_2}) \\
d(\lambda_{j^\prime};u_{i_1^\prime}) &d(-\lambda_{j^\prime};v_{i_2^\prime}) &d(-\lambda_{j^\prime};-u_{i_1^\prime}) &d(\lambda_{j^\prime};-v_{i_2^\prime}) 
\end{pmatrix}
.
\end{align*}
and denote $\nu_j$ by $\{d(\omega_{1\nu_j};q_{1\nu_j}), d(\omega_{2\nu_j};q_{2\nu_j}),\ldots,d(\omega_{\#\nu_j\nu_j};q_{\#\nu_j\nu_j})\}$. Then, Theorem 2.3.2 of \citet[p.21]{b81} and Lemma \ref{lem_cum} yield that 
\begin{align*}
&\frac{8^2\pi^2 L^4}{n^2M_n^4}\sum_{j,j^\prime=1}^{n-1}\sum_{i_1,i_1^\prime,i_2,i_2^\prime=1}^{M_n}{\rm Cum}\bigg(I_n(\lambda_j;u_{i_1},v_{i_2})\overline{I_n(\lambda_j;u_{i_1},v_{i_2})},\\
&\qquad\qquad\qquad\qquad\qquad\qquad\qquad
I_n(\lambda_{j^\prime};u_{i_1^\prime},v_{i_2^\prime})\overline{I_n(\lambda_{j^\prime};u_{i_1^\prime},v_{i_2^\prime})}\bigg)\\
=&
\frac{8^2\pi^2 L^4}{n^2M_n^4}\sum_{j,j^\prime=1}^{n-1}\sum_{i_1,i_1^\prime,i_2,i_2^\prime=1}^{M_n}
\frac{1}{(2\pi n)^4}\times\\
&
\sum_{(\nu_1,\ldots,\nu_p)}
\skakko{(2\pi)^{\# \nu_1-1}\Delta_n\skakko{\sum_{k=1}^{\# \nu_1}\omega_{k\nu_1}}
f(\omega_{1\nu_1},\ldots,\omega_{(\# \nu_1-1)\nu_1};
q_{1\nu_1},\ldots,q_{\#\nu_1\nu_1})+O(1)}
\\
&
\cdots
\skakko{(2\pi)^{\# \nu_p-1}\Delta_n\skakko{\sum_{k=1}^{\# \nu_p}\omega_{k\nu_p}}
f(\omega_{1\nu_p},\ldots,\omega_{(\# \nu_p-1)\nu_p};
q_{1\nu_p},\ldots,q_{\#\nu_p\nu_p})+O(1)},
\end{align*}
where the summation $\sum_{(\nu_1,\ldots,\nu_p)}$ extends over all indecomposable partitions $(\nu_1,\ldots,\nu_p)$  of the { aforementioned} table.
There are two types of partitions which provide the leading order of \eqref{var_case1}.  The first type are indecomposable partitions where the number of set in the partition is four. { Note that then, each set has cardinality two and the frequencies are $\lambda_{j^\prime}=2\pi-\lambda_j$ or $\lambda_j$} so that the sum of the frequencies for each set is zero (mod $2\pi$). The indecomposable partition of the type is, for example, 
\begin{align*}
\{&\{d(\lambda_j;u_{i_1}),d(\lambda_{j^\prime};u_{i_1^\prime})\},\{d(-\lambda_j;v_{i_2}),d(-\lambda_{j^\prime};v_{i_2^\prime})\},
\{d(-\lambda_j;-u_{i_1}),d(-\lambda_{j^\prime};-u_{i_1^\prime}) \},\\
&\{d(\lambda_j;-v_{i_2}),d(\lambda_{j^\prime};-v_{i_2^\prime}) \}\}
\end{align*}
with $\lambda_{j^\prime}=2\pi-\lambda_j$.
For that indecomposable partition, we have
\begin{align*}
&\frac{8^2\pi^2 L^4}{n^2M_n^4}\sum_{j,j^\prime=1}^{n-1}\sum_{i_1,i_1^\prime,i_2,i_2^\prime=1}^{M_n}
\frac{1}{(2\pi n)^4}\\
&\times
\skakko{(2\pi)^{\# \nu_1-1}\Delta_n\skakko{\sum_{k=1}^{\# \nu_1}\omega_{k\nu_1}}
f(\omega_{1\nu_1},\ldots,\omega_{(\# \nu_1-1)\nu_1};
q_{1\nu_1},\ldots,q_{\#\nu_1\nu_1})+O(1)}
\cdots\\
&
\skakko{(2\pi)^{\# \nu_p-1}\Delta_n\skakko{\sum_{k=1}^{\# \nu_p}\omega_{k\nu_p}}
f(\omega_{1\nu_p},\ldots,\omega_{(\# \nu_p-1)\nu_p};
q_{1\nu_p},\ldots,q_{\#\nu_p\nu_p})+O(1)}\\
=&\frac{8^2\pi^2 L^4}{n^2M_n^4}\sum_{j=1}^{n-1}\sum_{i_1,i_1^\prime,i_2,i_2^\prime=1}^{M_n}
\big(
f(\lambda_j;u_{i_1},u_{i_1^\prime})f(-\lambda_j;v_{i_2},v_{i_2^\prime})
f(\lambda_j;-u_{i_1},-u_{i_1^\prime})f(\lambda_j;-v_{i_2},-v_{i_2^\prime})\\
&\qquad\qquad\qquad\qquad\qquad+O(1/n)\big)\\
\leq&
\frac{8^2\pi^2 L^4}{n}
\skakko{\sup_{\lambda,u,v}
\left|f(\lambda;u,v)\right|}^4
+O(1/n^2).
\end{align*}
The { second type of leading-order partition consists in} indecomposable partitions with three sets where two sets have two elements and one set has four elements and the sum of the frequencies for each set is zero (mod $2\pi$) without any restriction of frequencies.
The indecomposable partition satisfying the above is, for example,
\begin{align*}
\{&\{d(\lambda_j;u_{i_1}),d(-\lambda_j;v_{i_2})\},
\{d(\lambda_{j^\prime};u_{i_1^\prime}),d(-\lambda_{j^\prime};v_{i_2^\prime})\},\{d(-\lambda_j;-u_{i_1}),d(-\lambda_{j^\prime};-u_{i_1^\prime})\},\\
& \{d(\lambda_j;-v_{i_2}),d(\lambda_{j^\prime};-v_{i_2^\prime}) \}\}
\end{align*}
Then, for that indecomposable partition, we have
\begin{align*}
&\frac{8^2\pi^2 L^4}{n^2M_n^4}\sum_{j,j^\prime=1}^{n-1}\sum_{i_1,i_1^\prime,i_2,i_2^\prime=1}^{M_n}
\frac{1}{(2\pi n)^4}\\
&\times
\skakko{(2\pi)^{\# \nu_1-1}\Delta_n\skakko{\sum_{k=1}^{\# \nu_1}\omega_{k\nu_1}}
f(\omega_{1\nu_1},\ldots,\omega_{(\# \nu_1-1)\nu_1};
q_{1\nu_1},\ldots,q_{\#\nu_1\nu_1})+O(1)}
\cdots\\
&
\skakko{(2\pi)^{\# \nu_p-1}\Delta_n\skakko{\sum_{k=1}^{\# \nu_p}\omega_{k\nu_p}}
f(\omega_{1\nu_p},\ldots,\omega_{(\# \nu_p-1)\nu_p};
q_{1\nu_p},\ldots,q_{\#\nu_p\nu_p})+O(1)}\\
=&\frac{8^2\pi^2 L^4}{n^2M_n^4}\sum_{j,j^\prime=1}^{n-1}\sum_{i_1,i_1^\prime,i_2,i_2^\prime=1}^{M_n}\\
&\skakko{
\frac{1}{2\pi n}f(\lambda_j;u_{i_1},v_{i_2})f(\lambda_{j^\prime};u_{i_1^\prime},v_{i_2^\prime})
f(-\lambda_j,\lambda_j,-\lambda_j^\prime;-u_{i_1},-v_{i_2},-u_{i_1^\prime},-v_{i_2^\prime})+O(1/n^2)}\\
\leq&
\frac{8^2\pi^2 L^4}{2\pi n}
\skakko{\sup_{\lambda,u,v}\left|f(\lambda;u,v)\right|}^2
\skakko{\sup_{\lambda,\lambda^\prime,u,u^\prime,v,v^\prime}\left|f(-\lambda,\lambda,-\lambda^\prime;-u,-v,-u^\prime,-v^\prime)\right|}\\
&+O(1/n^2).
\end{align*}
Since we impose the condition \eqref{as_sum} and other type indecomposable partitions provide faster convergence rate than $1/n$ , we obtain the desired result. \qed

\subsection{Proof of Theorem \ref{thm1}}
By Theorem 5.7 of \cite{v00}, it suffices to show the uniform laws of large numbers, that is, the fact that 
$\sup_{\bm{\theta}\in\bm{\Theta}}
\left|D_n( I_n,f_{\bm{\theta}})-\tilde D(f,f_{\bm{\theta}})\right|$ converges in probability to zero as $n\to\infty$.
{ Our proof consist in first proving} the stochastic equicontinuity of $\{D_n( I_n,f_{\bm{\theta}})-\tilde D(f,f_{\bm{\theta}}); \bm{\theta}\in\bm{\Theta}\}$,
that is, the fact that for any $\eta>0$ and any $\epsilon>0$, there exists a $\delta>0$ such that
\begin{align*}
\limsup_{n\to\infty}{\rm P}\skakko{\sup_{\bm{\theta}\in\bm{\Theta}}\sup_{\|\bm{\theta}-\bm{\theta}^\prime\|\leq \delta}
\left|D_n( I_n,f_{\bm{\theta}})- \tilde D(f,f_{\bm{\theta}})
-\skakko{D_n( I_n,f_{\bm{\theta^\prime}}) - \tilde D(f,f_{\bm{\theta^\prime}})
}\right|>\eta}<\epsilon.
\end{align*}
{Once we have stochastic equicontinuity, the result holds by the following argument:} by compactness of $\bm{\Theta}$, there exists the finite covering set $\{B({\bm\theta_k}, \delta);k=1,\ldots,K_\delta\}$, where
$B({\bm\theta_k};\delta)\coloneqq \{{\bm\theta}\in{\bm\Theta};\|{\bm\theta_k}-{\bm\theta}\|<\delta\}$. Then, we have
\begin{align*}
&{\rm P}\skakko{\sup_{\bm{\theta}\in\bm{\Theta}}
\left|
D_n( I_n,f_{\bm{\theta}})-\tilde D(f,f_{\bm{\theta}})
\right|>\eta}\\
=&
{\rm P}\skakko{\max_{k=1,\ldots,K_\delta}\sup_{\bm{\theta}\in B({\bm\theta_k};\delta)}
\left|
D_n( I_n,f_{\bm{\theta}})-\tilde D(f,f_{\bm{\theta}})
\right|>\eta}\\
\leq&
{\rm P}\skakko{\max_{k=1,\ldots,K_\delta}\sup_{\bm{\theta}\in B({\bm\theta_k};\delta)}
\left|
D_n( I_n,f_{\bm{\theta}})-\tilde D(f,f_{\bm{\theta}})
-\skakko{D_n( I_n,f_{\bm{\theta_k}})-\tilde D(f,f_{\bm{\theta_k}})}
\right|>\eta/2}\\
&+
{\rm P}\skakko{\max_{k=1,\ldots,K_\delta}
\left|
D_n( I_n,f_{\bm{\theta_k}})-\tilde D(f,f_{\bm{\theta_k}})
\right|>\eta/2},
\end{align*}
which, by the stochastic equicontinuity and pointwise consistency (Lemma \ref{lem1}), tends to zero.

First we show the stochastic equicontinuity of $\{\tilde D(f,f_{\bm{\theta}}); \bm{\theta}\in\bm{\Theta}\}$. Markov's inequality and {the UBPH assumption} of  $f_{\bm{\theta}}(\lambda;u,v)$ yield
{
\begin{align*}
&{\rm P}\skakko{\sup_{\bm{\theta}\in\bm{\Theta}}\sup_{\|\bm{\theta}-\bm{\theta}^\prime\|\leq \delta}
\left|\tilde D(f,f_{\bm{\theta}})
- \tilde D(f,f_{\bm{\theta^\prime}})
\right|>\eta}\\
\leq&
\frac{1}{\eta}\sup_{\bm{\theta}\in\bm{\Theta}}\sup_{\|\bm{\theta}-\bm{\theta}^\prime\|\leq \delta}
\left|\tilde D(f,f_{\bm{\theta}})
- \tilde D(f,f_{\bm{\theta^\prime}})
\right|\\
=&
\frac{1}{\eta}\sup_{\bm{\theta}\in\bm{\Theta}}\sup_{\|\bm{\theta}-\bm{\theta}^\prime\|\leq \delta}
\big|
\int_{0}^{2\pi}\int_{-L}^L\int_{-L}^L
|f(\lambda;u,v)-f_{\bm{\theta}}(\lambda;u,v))|^2\\
&\qquad\qquad\qquad\qquad\qquad\qquad\qquad
-
|f(\lambda;u,v)-f_{\bm{\theta}^\prime}(\lambda;u,v))|^2
{\rm d}u{\rm d}v{\rm d}{\lambda}
\big|\\
\leq&
\frac{1}{\eta}\sup_{\bm{\theta}\in\bm{\Theta}}\sup_{\|\bm{\theta}-\bm{\theta}^\prime\|\leq \delta}
\int_{0}^{2\pi}\int_{-L}^L\int_{-L}^L
|f_{\bm{\theta}}(\lambda;u,v)- f_{\bm{\theta}^\prime}(\lambda;u,v)|\\
&\qquad\qquad\qquad\qquad\qquad\qquad\qquad\times
|2f(\lambda;u,v)-f_{\bm{\theta}}(\lambda;u,v)-f_{\bm{\theta}^\prime}(\lambda;u,v)|
{\rm d}u{\rm d}v{\rm d}{\lambda}
\\
\leq&
\frac{4C^\prime}{\eta}
\skakko{
2\pi (2L)^2 \delta^\gamma
C^\prime
\varepsilon^{-\tau}
+
2C^\prime \varepsilon
},
\end{align*}
which can be arbitrary small, for  $\varepsilon=\delta^{\gamma/(1+\tau)}$, by taking a small $\delta$.}

Next we show the stochastic equicontinuity of 
$\{D_n( I_n,f_{\bm{\theta}}); \bm{\theta}\in\bm{\Theta}\}$.
{Under Assumption (c), $f_{\theta}(\lambda,u,v)$ is uniformly bounded by $C'$ and
\begin{align*}
&{\rm P}\skakko{\sup_{\bm{\theta}\in\bm{\Theta}}\sup_{\|\bm{\theta}-\bm{\theta}^\prime\|\leq \delta}
|D_n( I_n,f_{\bm{\theta}})-
D_n( I_n,f_{\bm{\theta^\prime}})
|>\eta}\\
\leq
&{\rm P}\bigg(
\frac{8\pi L^2}{nM_n^2}\sum_{j=1}^{n-1}\sum_{i_1,i_2=1}^{M_n}
\Big(
|f_{\bm{\theta}}(\lambda_j;u_{i_1},v_{i_2})- f_{\bm{\theta}^\prime}(\lambda_j;u_{i_1},v_{i_2})|\\
&\qquad\qquad\qquad\qquad\qquad
\times
\Big(
2\left|I_n(\lambda_j;u_{i_1},v_{i_2})\right|
+
2\sup_{\bm{\theta},\lambda,u,v}\left|f_{\bm{\theta}}(\lambda;u,v)\right|\Big)\Big)
>\eta\Big)\\
\leq
&
\frac{8\pi L^2}{\eta}
\skakko{
2\sup_{\lambda,u,v}{\rm E}\left|I_n(\lambda;u,v)\right|
+
2C^\prime}
\frac{1}{nM_n^2}
\sum_{j=1}^{n-1}\sum_{i_1,i_2=1}^{M_n}
|f_{\bm{\theta}}(\lambda_j;u_{i_1},v_{i_2})- f_{\bm{\theta}^\prime}(\lambda_j;u_{i_1},v_{i_2})|.
\end{align*}
}
{ Since we showed in the proof of Lemma \ref{lem1} that
\begin{align*}
&{\rm E}\skakko{I_n(\lambda_j;u_{i_1},v_{i_2}) \overline{I_n(\lambda_j;u_{i_1},v_{i_2})}}\\
=&f(\lambda_j;u_{i_1},v_{i_2})f(-\lambda_j;-u_{i_1},-v_{i_2})
+f(\lambda_j;u_{i_1},-u_{i_1})f(-\lambda_j;v_{i_2},-v_{i_2})\\
&+f(\lambda_j;u_{i_1},-v_{i_2})f(-\lambda_j;v_{i_2},-u_{i_1})\Delta_n\skakko{2\lambda_j}\Delta_n\skakko{-2\lambda_j}/n^2
+O\skakko{{1}/{n}},
\end{align*}
where the error term is uniform on $(\lambda_1,\ldots,\lambda_{k-1})\in[0,2\pi]^{k-1}$ and $(u_1,\ldots,u_k)\in\mathbb R^k$,
we have that $\sup_{\lambda,u,v}{\rm E}\left|I_n(\lambda;u,v)\right|$ is bounded. As $f_{{\bm \theta}}$ is ${\rm UBPH}({\mathcal S},\alpha,\beta,\gamma)$, we have that
\begin{align*}
&\frac{1}{nM_n^2}
\sum_{j=1}^{n-1}\sum_{i_1,i_2=1}^{M_n}
I_{(i_1,i_2);(u_{i_1}, v_{i_2}) \in \mathcal{S}_{\epsilon}}|f_{\bm{\theta}}(\lambda_j;u_{i_1},v_{i_2})- f_{\bm{\theta}^\prime}(\lambda_j;u_{i_1},v_{i_2})|=O(\epsilon+\frac{1}{M_n})
\end{align*}
and that
\begin{align*}
&\frac{1}{nM_n^2}
\sum_{j=1}^{n-1}\sum_{i_1,i_2=1}^{M_n}
I_{(i_1,i_2);(u_{i_1}, v_{i_2}) \notin \mathcal{S}_{\epsilon}}|f_{\bm{\theta}}(\lambda_j;u_{i_1},v_{i_2})- f_{\bm{\theta}^\prime}(\lambda_j;u_{i_1},v_{i_2})|\\
\le&\frac{1}{nM_n^2}\sum_{j=1}^{n-1}\sum_{i_1,i_2=1}^{M_n} \delta^{\gamma}\\
=&\delta^{\gamma}
\end{align*}
we now consider $\epsilon=1/M_n$ and, taking $n\rightarrow\infty$ we get stochastic equicontinuity.}
\qed

\subsection{Proof of Theorem \ref{thm2}}
We prove this theorem in a similar manner as Theorem 3.1 of \cite{d00}.
By the mean value theorem 
\begin{align*}
\frac{\partial}{\partial \bm{\theta}}D_n(I_n,f_{{\bm{\theta}}})\Big|_{\bm{\theta}={\hat{\bm{\theta}}_n}}
-\frac{\partial}{\partial \bm{\theta}}D_n(I_n,f_{{\bm{\theta}}})\Big|_{\bm{\theta}={{\bm{\theta}}_0}}
=
\frac{\partial^2}{\partial \bm{\theta}\partial \bm{\theta}^\top}
D_n(I_n,f_{{\bm{\theta}}})\Big|_{\bm{\theta}={{\bm{\theta}}_n^*}}
\skakko{\hat{\bm{\theta}}_n-{\bm{\theta}}_0},
\end{align*}
where $\bm{\theta}_0 \gtreqless\bm{\theta}_n^* \gtreqless {\hat{\bm{\theta}}_n}$. 
Note that $\sqrt n {\partial}/{(\partial \bm{\theta})}D_n(I_n,f_{{\bm{\theta}}})\Big|_{\bm{\theta}={\hat{\bm{\theta}}_n}}=0$ when $\hat{\bm{\theta}}_n$ belongs to the interior of $\bm{\Theta}$.
When $\hat{\bm{\theta}}_n$ belongs to the boundary of $\bm{\Theta}$, $\sqrt n {\partial}/{(\partial \bm{\theta})}D_n(I_n,f_{{\bm{\theta}}})\Big|_{\bm{\theta}={\hat{\bm{\theta}}_n}}$ tends  in probability to zero since, for any $\epsilon>0$, there exists a $\delta>0$ such that $\left|{\hat{\bm{\theta}}}_n-{\bm{\theta}_0}\right|\geq \delta$ and thus
$$
{\rm P}\skakko{\left|\sqrt n \frac{\partial}{\partial \bm{\theta}}D_n(I_n,f_{{\bm{\theta}}})\Big|_{\bm{\theta}={\hat{\bm{\theta}}}_n}\right|>\epsilon}
\leq 
{\rm P}\skakko{\left|{\hat{\bm{\theta}}}_n-{\bm{\theta}_0}\right|\geq \delta}\to 0\quad\text{as $n\to\infty$}.
$$
The proof can now be completed provided we can show that
\begin{enumerate}
\item[(i)]
\begin{align*}
\frac{\partial^2}{\partial \bm{\theta}\partial \bm{\theta}^\top}D_n( I_n,f_{\bm{\theta}})\Big|_{\bm{\theta}={{\bm{\theta}}_n^*}}-
\frac{\partial^2}{\partial \bm{\theta}\partial \bm{\theta}^\top}D_n( I_n,f_{\bm{\theta}})\Big|_{\bm{\theta}={{\bm{\theta}}_0}}=o_p(1)\quad\text{ as $n\to\infty$},
\end{align*}
\item[(ii)]
\begin{align*}
\frac{\partial^2}{\partial \bm{\theta}\partial \bm{\theta}^\top}D_n( I_n,f_{\bm{\theta}})\Big|_{\bm{\theta}={{\bm{\theta}}_0}}
-
{\bm J}=o_p(1)\quad\text{as $n\to\infty$},
\end{align*}
where
\begin{align*}
{\bm J}\coloneqq &
-\int_{0}^{2\pi}\int_{-L}^L\int_{-L}^L
f(\lambda;u,v)\frac{\partial^2}{\partial \bm{\theta}\partial \bm{\theta}^\top}\overline{f_{\bm{\theta}}(\lambda;u,v)}\Big|_{\bm{\theta}={{\bm{\theta}}_0}}
{\rm d}u{\rm d}v{\rm d}{\lambda}\\
&-\int_{0}^{2\pi}\int_{-L}^L\int_{-L}^L
\frac{\partial^2}{\partial \bm{\theta}\partial \bm{\theta}^\top}f_{\bm{\theta}}(\lambda;u,v)\Big|_{\bm{\theta}={{\bm{\theta}}_0}}\overline{f(\lambda;u,v)}
{\rm d}u{\rm d}v{\rm d}{\lambda}\\
&+\int_{0}^{2\pi}\int_{-L}^L\int_{-L}^L
\frac{\partial^2}{\partial \bm{\theta}\partial \bm{\theta}^\top}
\skakko{f_{\bm{\theta}}(\lambda;u,v)\overline{f_{\bm{\theta}}(\lambda;u,v)}}\Big|_{\bm{\theta}={{\bm{\theta}}_0}}
{\rm d}u{\rm d}v{\rm d}{\lambda}.
\end{align*}

\item[(iii)]
$\sqrt n{\partial}/{(\partial \bm{\theta})}D_n(I_n,f_{\bm{\theta}})\big|_{\bm{\theta}={{\bm{\theta}}_0}}$ converges in distribution to a centered normal distribution with variance $\bm I$ as $n\to\infty$,
where 
\begin{align*}
\bm I\coloneqq &(I_{st})_{s,t=1,\ldots,d}, \quad I_{st}\coloneqq L_{1_{st}}+L_{2_{st}}+L_{3_{st}}+L_{4_{st}}
\end{align*}
with, for the delta function $\delta(x)$ and $\textbf{d}{\bm w}={\rm d}u{\rm d}u^\prime{\rm d}v{\rm d}v^\prime{\rm d}\lambda{\rm d}\lambda^\prime$,
\begin{align}\nonumber
L_{1_{st}}\coloneqq 
&{2\pi}
\int_{[0,2\pi]^2}\int_{[-L,L]^4}
\Big[
\delta(\lambda-\lambda^\prime)
f(\lambda;u,v^\prime)
f(-\lambda;v,u^\prime)
\\\nonumber
&+
\delta(2\pi-\lambda-\lambda^\prime)
f(\lambda;u,u^\prime)
f(-\lambda;v,v^\prime)
+
f(\lambda,-\lambda,\lambda^\prime;u,v,u^\prime,v^\prime)\Big]\\\nonumber
&\times
\frac{\partial}{\partial {\theta_s}}\overline{f_{\bm{\theta}}(\lambda;u,v)}\Big|_{\bm{\theta}={{\bm{\theta}}_0}}\frac{\partial}{\partial {\theta_t}}\overline{f_{\bm{\theta}}(\lambda^\prime;u^\prime,v^\prime)}\Big|_{\bm{\theta}={{\bm{\theta}}_0}}
\textbf{d}{\bm w}
\\\nonumber
L_{2_{st}}\coloneqq 
&2\pi
\int_{[0,2\pi]^2}\int_{[-L,L]^4}
\Big[
\delta(\lambda-\lambda^\prime)
f(\lambda;u,-u^\prime)
f(-\lambda;v,-v^\prime)\\\nonumber
&+
\delta(2\pi-\lambda-\lambda^\prime)
f(\lambda;u,-v^\prime)
f(-\lambda;v,-u^\prime)
+
f(\lambda,-\lambda,-\lambda^\prime;u,v,-u^\prime,-v^\prime)\Big]\\\nonumber
&\times
\frac{\partial}{\partial {\theta_s}}\overline{f_{\bm{\theta}}(\lambda;u,v)}\Big|_{\bm{\theta}={{\bm{\theta}}_0}}\frac{\partial}{\partial {\theta_t}}{f_{\bm{\theta}}(\lambda^\prime;u^\prime,v^\prime)}\Big|_{\bm{\theta}={{\bm{\theta}}_0}}
\textbf{d}{\bm w}
\\\nonumber
L_{3_{st}}\coloneqq 
&2\pi
\int_{[0,2\pi]^2}\int_{[-L,L]^4}
\Big[
\delta(\lambda-\lambda^\prime)
f(\lambda;-v,v^\prime)
f(-\lambda;-u,u^\prime)\\\nonumber
&+
\delta(2\pi-\lambda-\lambda^\prime)
f(\lambda;-v,u^\prime)
f(-\lambda;-u,v^\prime)
+
f(-\lambda,\lambda,\lambda^\prime;-u,-v,u^\prime,v^\prime)\Big]\\\nonumber
&\times
\frac{\partial}{\partial {\theta_s}}{f_{\bm{\theta}}(\lambda;u,v)}\Big|_{\bm{\theta}={{\bm{\theta}}_0}}\frac{\partial}{\partial {\theta_t}}\overline{f_{\bm{\theta}}(\lambda^\prime;u^\prime,v^\prime)}\Big|_{\bm{\theta}={{\bm{\theta}}_0}}
\textbf{d}{\bm w}
\\\nonumber
L_{4_{st}}\coloneqq 
&2\pi
\int_{[0,2\pi]^2}\int_{[-L,L]^4}
\Big[
\delta(\lambda-\lambda^\prime)
f(\lambda;-v,-u^\prime)
f(-\lambda;-u,-v^\prime)
\\\nonumber
&
+
\delta(2\pi-\lambda-\lambda^\prime)
f(\lambda;-v,-v^\prime)
f(-\lambda;-u,-u^\prime)
\\\nonumber
&+
f(-\lambda,\lambda,-\lambda^\prime;-u,-v,-u^\prime,-v^\prime)\Big]\\\nonumber
&
\times
\frac{\partial}{\partial {\theta_s}}{f_{\bm{\theta}}(\lambda;u,v)}\Big|_{\bm{\theta}={{\bm{\theta}}_0}}\frac{\partial}{\partial {\theta_t}}{f_{\bm{\theta}}(\lambda^\prime;u^\prime,v^\prime)}\Big|_{\bm{\theta}={{\bm{\theta}}_0}}
\textbf{d}{\bm w}.
\end{align}
Note that the integration over frequencies is defined on the open domain $(0, 2\pi)^2$ and thus the Dirac delta functions $\delta(\cdot)$ do not take values on the boundary of the domain.
\end{enumerate}

\subsubsection{Proof of (i)}

{
Without loss of generality, it can be assumed that the set $\mathcal{S}$, the exponents $\tau,\alpha,\beta,\gamma$ and the uniform constant $C^{\prime}$ are the same for $$
f_{\bm{\theta}}(\lambda;u,v),\quad \frac{\partial}{\partial \bm{\theta}}f_{\bm{\theta}}(\lambda;u,v),\quad \text{and}
\quad
\frac{\partial^2}{\partial \bm{\theta}\partial \bm{\theta}^\top} f_{\bm{\theta}}(\lambda;u,v).$$ Then, straightforward calculations yield
\begin{align*}
&{\rm P}\skakko{\sup_{\bm{\theta}\in\bm{\Theta}}\sup_{\|\bm{\theta}-\bm{\theta}^\prime\|\leq \delta}
\left\|\frac{\partial^2}{\partial \bm{\theta}\partial \bm{\theta}^\top}D_n( I_n,f_{\bm{\theta}})-
\frac{\partial^2}{\partial \bm{\theta^\prime}\partial {\bm{\theta}^\prime}^\top}D_n( I_n,f_{\bm{\theta^\prime}})
\right\|>\eta}\\
\leq
&
\frac{16\pi L^2}{\eta}
\skakko{
2\sup_{\lambda,u,v}{\rm E}\left|I_n(\lambda;u,v)\right|
+
2C^\prime}\\
&\times
\frac{1}{nM_n^2}
\sum_{j=1}^{n-1}\sum_{i_1,i_2=1}^{M_n}
\bigg(\left\|\frac{\partial^2}{\partial \bm{\theta^\prime}\partial {\bm{\theta}^\prime}^\top}f_{\bm{\theta}}(\lambda_j;u_{i_1},v_{i_2})- \frac{\partial^2}{\partial \bm{\theta^\prime}\partial {\bm{\theta}^\prime}^\top}f_{\bm{\theta}^\prime}(\lambda_j;u_{i_1},v_{i_2})\right\|\\
&
+ 2 C^{\prime} |f_{\bm \theta}(\lambda_j,u_{i_{1}}, u_{j_{2}}) - f_{{\bm \theta}^{\prime}}(\lambda_j,u_{i_{1}}, u_{j_{2}})|\\
&
+\Big\| \frac{\partial}{\partial \bm{\theta}}f_{\bm{\theta}}(\lambda_j;u_{i_1},v_{i_2})\frac{\partial}{\partial \bm{\theta}^\top}f_{\bm{\theta}}(\lambda_j;u_{i_1},v_{i_2})
- 
\frac{\partial}{\partial \bm{\theta}^{\prime}}f_{\bm{\theta}^\prime}(\lambda_j;u_{i_1},v_{i_2})\frac{\partial}{\partial{ \bm{\theta}^{\prime}}^\top}f_{\bm{\theta}^\prime}(\lambda_j;u_{i_1},v_{i_2})\Big\|\bigg).
\end{align*}
Now, by the same argument as in the proof of Theorem \ref{thm1} and using the fact that all the functions involved in the bound are ${\rm UBPH}$, the previous expression can be arbitrarily small by choosing small $\delta$.} Thus, the stochastic equicontinuity of $\{{\partial^2}/{(\partial \bm{\theta}\partial \bm{\theta}^\top)} D_n( I_n,f_{\bm{\theta}})\}$ holds. Hence, we obtain (i) since
\begin{align*}
&{\rm P}\skakko{
\left\|\frac{\partial^2}{\partial \bm{\theta}\partial \bm{\theta}^\top}D_n( I_n,f_{\bm{\theta}})\Big|_{\bm{\theta}={{\bm{\theta}}_n^*}}-
\frac{\partial^2}{\partial \bm{\theta^\prime}\partial {\bm{\theta}^\prime}^\top}D_n( I_n,f_{\bm{\theta^\prime}})\Big|_{\bm{\theta}={{\bm{\theta}}_0}}
\right\|>\eta}\\
\leq&
{\rm P}\skakko{\sup_{\bm{\theta}\in\bm{\Theta}}\sup_{\|\bm{\theta}-\bm{\theta}^\prime\|\leq \delta}
\left\|\frac{\partial^2}{\partial \bm{\theta}\partial \bm{\theta}^\top}D_n( I_n,f_{\bm{\theta}})-
\frac{\partial^2}{\partial \bm{\theta^\prime}\partial {\bm{\theta}^\prime}^\top}D_n( I_n,f_{\bm{\theta^\prime}})
\right\|>\eta}\\
&+
{\rm P}\skakko{\|{{\bm{\theta}}_n^*}-{{\bm{\theta}}_0}\|>\delta},
\end{align*}
which tends to zero as $n\to\infty$.

\subsubsection{Proof of (ii)}
{ By using Lemma \ref{lem_cum} and the same arguments as in the proof of Lemma \ref{lem1}, we get that as $n, M_n \rightarrow \infty$ and $\epsilon \rightarrow 0$,
\begin{align*}
&{\rm E}\skakko{\frac{\partial^2}{\partial \bm{\theta}\partial \bm{\theta}^\top}D_n(I_n,f_{\bm{\theta}})\Big|_{\bm{\theta}={{\bm{\theta}}_0}}}\\
=&
\frac{8\pi L^2}{nM_n^2}\sum_{j=1}^{n-1}\sum_{i_1,i_2=1}^{M_n}
{\rm E}\skakko{\frac{\partial^2}{\partial \bm{\theta}\partial \bm{\theta}^\top}
| I_n(\lambda_j;u_{i_1},v_{i_2})-f_{\bm{\theta}}(\lambda_j;u_{i_1},v_{i_2}))|^2\Big|_{\bm{\theta}={{\bm{\theta}}_0}}}\\
=&-\frac{8\pi L^2}{nM_n^2}\sum_{j=1}^{n-1}\sum_{i_1,i_2=1}^{M_n}
{\rm E}\skakko{I_n(\lambda_j;u_{i_1},v_{i_2})}\frac{\partial^2}{\partial \bm{\theta}\partial \bm{\theta}^\top}\overline{f_{\bm{\theta}}(\lambda_j;u_{i_1},v_{i_2})}\Big|_{\bm{\theta}={{\bm{\theta}}_0}}\\
&-\frac{8\pi L^2}{nM_n^2}\sum_{j=1}^{n-1}\sum_{i_1,i_2=1}^{M_n}
\frac{\partial^2}{\partial \bm{\theta}\partial \bm{\theta}^\top}f_{\bm{\theta}}(\lambda_j;u_{i_1},v_{i_2}))\Big|_{\bm{\theta}={{\bm{\theta}}_0}}{\rm E}\skakko{\overline{I_n(\lambda_j;u_{i_1},v_{i_2})}}\\
&+\frac{8\pi L^2}{nM_n^2}\sum_{j=1}^{n-1}\sum_{i_1,i_2=1}^{M_n}
\frac{\partial^2}{\partial \bm{\theta}\partial \bm{\theta}^\top}
\skakko{f_{\bm{\theta}}(\lambda_j;u_{i_1},v_{i_2}))\overline{f_{\bm{\theta}}(\lambda_j;u_{i_1},v_{i_2})}}\Big|_{\bm{\theta}={{\bm{\theta}}_0}}\\
=
&-\int_{0}^{2\pi}\int_{-L}^L\int_{-L}^L
f(\lambda;u,v)\frac{\partial^2}{\partial \bm{\theta}\partial \bm{\theta}^\top}\overline{f_{\bm{\theta}}(\lambda;u,v)}\Big|_{\bm{\theta}={{\bm{\theta}}_0}}{\rm d}u{\rm d}v{\rm d}{\lambda}\\
&
-\int_{0}^{2\pi}\int_{-L}^L\int_{-L}^L
\frac{\partial^2}{\partial \bm{\theta}\partial \bm{\theta}^\top}f_{\bm{\theta}}(\lambda;u,v))\Big|_{\bm{\theta}={{\bm{\theta}}_0}}\overline{f(\lambda;u,v)}{\rm d}u{\rm d}v{\rm d}{\lambda}\\
&+\int_{0}^{2\pi}\int_{-L}^L\int_{-L}^L
\frac{\partial^2}{\partial \bm{\theta}\partial \bm{\theta}^\top}
\skakko{f_{\bm{\theta}}(\lambda;u,v)\overline{f_{\bm{\theta}}(\lambda;u,v)}}\Big|_{\bm{\theta}={{\bm{\theta}}_0}}{\rm d}u{\rm d}v{\rm d}{\lambda}\\
&
+O\skakko{
\skakko{n^{-\min \{\tilde{\alpha}, \alpha\}}+M_n^{-\min\{\tilde{\beta},\beta\}}}^{{1}/{(\max\{\tilde{\tau},\tau\}+1)}}
+M_n^{-1}+n^{-1}
},
\end{align*}
where $\tilde{\tau}$, $\tilde{\alpha}$ and $\tilde{\beta}$ denote respectively the ${\rm PH}$ exponents associated to $\epsilon$, $u$ and $v$ for the generalized spectrum of order 2, $f(\lambda, u, v)$.} We can also show that $${\rm Var}\skakko{{\partial^2}/{(\partial \bm{\theta}\partial \bm{\theta}^\top)}D_n(I_n,f_{\bm{\theta}})}=O(1/n)\quad\text{as $n\to\infty$}$$ along the lines of Lemma \ref{lem1}. 

\subsubsection{Proof of (iii)}
We show the asymptotic normality of $\sqrt n{\partial}/{(\partial \bm{\theta})}D_n(I_n,f_{\bm{\theta}})\big|_{\bm{\theta}={{\bm{\theta}}_0}}$ by the cumulant method. { Recall now that we assumed--without loss of generality--that $\frac{\partial}{\partial \bm{\theta}}f_{\bm{\theta}}(\lambda_j;u_{i_1},v_{i_2})$ is ${\rm UBPH}(\mathcal{S},\tau,\alpha, \beta,\gamma)$.
Then, we observe that,{ by same arguments as in the proof of Lemma \ref{lem1},}
\begin{align*}
&{\rm E}\skakko{\frac{\partial}{\partial \bm{\theta}}D_n(I_n,f_{\bm{\theta}})\Big|_{\bm{\theta}={{\bm{\theta}}_0}}}\\
=&
\frac{8\pi L^2}{nM_n^2}\sum_{j=1}^{n-1}\sum_{i_1,i_2=1}^{M_n}
{\rm E}\skakko{\frac{\partial}{\partial \bm{\theta}}
| I_n(\lambda_j;u_{i_1},v_{i_2})-f_{\bm{\theta}}(\lambda_j;u_{i_1},v_{i_2}))|^2\Big|_{\bm{\theta}={{\bm{\theta}}_0}}}\\
=&-\frac{8\pi L^2}{nM_n^2}\sum_{j=1}^{n-1}\sum_{i_1,i_2=1}^{M_n}
{\rm E}\skakko{I_n(\lambda_j;u_{i_1},v_{i_2})}\frac{\partial}{\partial \bm{\theta}}\overline{f_{\bm{\theta}}(\lambda_j;u_{i_1},v_{i_2})}\Big|_{\bm{\theta}={{\bm{\theta}}_0}}\\
&-\frac{8\pi L^2}{nM_n^2}\sum_{j=1}^{n-1}\sum_{i_1,i_2=1}^{M_n}
\frac{\partial}{\partial \bm{\theta}}f_{\bm{\theta}}(\lambda_j;u_{i_1},v_{i_2}))\Big|_{\bm{\theta}={{\bm{\theta}}_0}}{\rm E}\skakko{\overline{I_n(\lambda_j;u_{i_1},v_{i_2})}}\\
&
+\frac{8\pi L^2}{nM_n^2}\sum_{j=1}^{n-1}\sum_{i_1,i_2=1}^{M_n}
\frac{\partial}{\partial \bm{\theta}}
\skakko{f_{\bm{\theta}}(\lambda_j;u_{i_1},v_{i_2}))\overline{f_{\bm{\theta}}(\lambda_j;u_{i_1},v_{i_2})}}\Big|_{\bm{\theta}={{\bm{\theta}}_0}}\\
=
&-\int_{0}^{2\pi}\int_{-L}^L\int_{-L}^L
f(\lambda;u,v)\frac{\partial}{\partial \bm{\theta}}\overline{f_{\bm{\theta}}(\lambda;u,v)}\Big|_{\bm{\theta}={{\bm{\theta}}_0}}{\rm d}u{\rm d}v{\rm d}{\lambda}\\
&
-\int_{0}^{2\pi}\int_{-L}^L\int_{-L}^L
\frac{\partial}{\partial \bm{\theta}}f_{\bm{\theta}}(\lambda;u,v))\Big|_{\bm{\theta}={{\bm{\theta}}_0}}\overline{f(\lambda;u,v)}{\rm d}u{\rm d}v{\rm d}{\lambda}\\
&+\int_{0}^{2\pi}\int_{-L}^L\int_{-L}^L
\frac{\partial}{\partial \bm{\theta}}
\skakko{f_{\bm{\theta}}(\lambda;u,v)\overline{f_{\bm{\theta}}(\lambda;u,v)}}\Big|_{\bm{\theta}={{\bm{\theta}}_0}}{\rm d}u{\rm d}v{\rm d}{\lambda}\\
&
+O\skakko{
\skakko{n^{-\min \{\tilde{\alpha}, \alpha\}}+M_n^{-\min\{\tilde{\beta},\beta\}}}^{{1}/{(\max\{\tilde{\tau},\tau\}+1)}}
+M_n^{-1}+n^{-1}
}\\
=
&\int_{0}^{2\pi}\int_{-L}^L\int_{-L}^L
\frac{\partial}{\partial \bm{\theta}}| f(\lambda;u,v)-f_{\bm{\theta}}(\lambda;u,v))|^2\Big|_{\bm{\theta}={{\bm{\theta}}_0}}{\rm d}u{\rm d}v{\rm d}{\lambda}\\
&
-\int_{0}^{2\pi}\int_{-L}^L\int_{-L}^L
\frac{\partial}{\partial \bm{\theta}}
\bigg(f(\lambda;u,v)\overline{f(\lambda;u,v)}\bigg)\Big|_{\bm{\theta}={{\bm{\theta}}_0}}
{\rm d}u{\rm d}v{\rm d}{\lambda}\\
&
+O\skakko{
\skakko{n^{-\min \{\tilde{\alpha}, \alpha\}}+M_n^{-\min\{\tilde{\beta},\beta\}}}^{{1}/{(\max\{\tilde{\tau},\tau\}+1)}}
+M_n^{-1}+n^{-1}
}\\
=
&\frac{\partial}{\partial \bm{\theta}} D(f,f_{\bm{\theta}}) \Big|_{\bm{\theta}={{\bm{\theta}}_0}}+O\skakko{
\skakko{n^{-\min \{\tilde{\alpha}, \alpha\}}+M_n^{-\min\{\tilde{\beta},\beta\}}}^{{1}/{(\max\{\tilde{\tau},\tau\}+1)}}
+M_n^{-1}+n^{-1}
}.
\end{align*}

Since $D(f,f_{\bm{\theta}})$ has a unique minimum at ${\bm{\theta}}_0\in \mathring{{\bm{\Theta}}}$, 
\begin{align*}
&\sqrt n{\rm E}\skakko{{\partial}/{(\partial \bm{\theta})}D_n(I_n,f_{\bm{\theta}})\Big|_{\bm{\theta}=\bm{\theta}_0}}\\
&=O\skakko{n^{1/2}\skakko{n^{-\min \{\tilde{\alpha}, \alpha\}}+M_n^{-\min\{\tilde{\beta},\beta\}}}^{{1}/{(\max\{\tilde{\tau},\tau\}+1)}}
+n^{1/2}M_n^{-1}+n^{-1/2}}
.\end{align*}
Now, note that if $$\frac{\min\{\alpha, \tilde{\alpha}, \beta, \tilde{\beta}\}}{\max\{\tau, \tilde{\tau}\}+1}> 1/2$$ and $M^{-1}_n=o(n^{-1/2})$, then $\sqrt n{\rm E}\skakko{{\partial}/{(\partial \bm{\theta})}D_n(I_n,f_{\bm{\theta}})\Big|_{\bm{\theta}=\bm{\theta}_0}}=o(1)$.} 

Next, we consider the covariance of ${\partial}/{(\partial \bm{\theta})}D_n(I_n,f_{\bm{\theta}})\big|_{\bm{\theta}={{\bm{\theta}}_0}}$. Elementary calculation yield
\begin{align*}
&{\rm Cov}\skakko{
\frac{\partial}{\partial{\theta_s}}D_n(I_n,f_{\bm{\theta}})\Big|_{\bm{\theta}={{\bm{\theta}}_0}},
\frac{\partial}{\partial{\theta_t}}D_n(I_n,f_{\bm{\theta}})\Big|_{\bm{\theta}={{\bm{\theta}}_0}}}\\
=&
\frac{8^2\pi^2 L^4}{n^2M_n^4}\sum_{j,j^\prime=1}^{n-1}\sum_{i_1,i_1^\prime,i_2,i_2^\prime=1}^{M_n}
{\rm Cum}\Bigg(
\frac{\partial}{\partial {\theta_s}}| I_n(\lambda_j;u_{i_1},v_{i_2})-f_{\bm{\theta}}(\lambda_j;u_{i_1},v_{i_2})|^2\Big|_{\bm{\theta}={{\bm{\theta}}_0}},\\
&\qquad\qquad\qquad\qquad\qquad\qquad\qquad
\frac{\partial}{\partial {\theta_t}}| I_n(\lambda_{j^\prime};u_{i_1^\prime},v_{i_2^\prime})-f_{\bm{\theta}}(\lambda_{j^\prime};u_{i_1^\prime},v_{i_2^\prime})|^2\Big|_{\bm{\theta}={{\bm{\theta}}_0}}
\Bigg)\\
=&L_{1_{st}n}+L_{2_{st}n}+L_{3_{st}n}+L_{4_{st}n},
\end{align*}
where
\begin{align*}
L_{{1_{st}}n}\coloneqq &\frac{8^2\pi^2 L^4}{n^2M_n^4}\sum_{j,j^\prime=1}^{n-1}\sum_{i_1,i_1^\prime,i_2,i_2^\prime=1}^{M_n}{\rm Cum}\left(I_n(\lambda_j;u_{i_1},v_{i_2}),I_n(\lambda_{j^\prime};u_{i_1^\prime},v_{i_2^\prime})\right)\\
&\times\frac{\partial}{\partial {\theta_s}}\overline{f_{\bm{\theta}}(\lambda_j;u_{i_1},v_{i_2})}
\frac{\partial}{\partial {\theta_t}}\overline{f_{\bm{\theta}}(\lambda_{j^\prime};u_{i_1^\prime},v_{i_2^\prime})},\\
L_{{2_{st}}n}\coloneqq &\frac{8^2\pi^2 L^4}{n^2M_n^4}\sum_{j,j^\prime=1}^{n-1}\sum_{i_1,i_1^\prime,i_2,i_2^\prime=1}^{M_n}{\rm Cum}\left(I_n(\lambda_j;u_{i_1},v_{i_2}),\overline{I_n(\lambda_{j^\prime};u_{i_1^\prime},v_{i_2^\prime})}\right)\\
&\times
\frac{\partial}{\partial {\theta_s}}\overline{f_{\bm{\theta}}(\lambda_j;u_{i_1},v_{i_2})}
\frac{\partial}{\partial {\theta_t}}f_{\bm{\theta}}(\lambda_{j^\prime};u_{i_1^\prime},v_{i_2^\prime}),\\
L_{{3_{st}}n}\coloneqq &\frac{8^2\pi^2 L^4}{n^2M_n^4}\sum_{j,j^\prime=1}^{n-1}\sum_{i_1,i_1^\prime,i_2,i_2^\prime=1}^{M_n}{\rm Cum}\left(\overline{I_n(\lambda_j;u_{i_1},v_{i_2})}
,I_n(\lambda_{j^\prime};u_{i_1^\prime},v_{i_2^\prime})\right)\\
&\times
\frac{\partial}{\partial {\theta_s}}f_{\bm{\theta}}(\lambda_j;u_{i_1},v_{i_2})
\frac{\partial}{\partial {\theta_t}}\overline{f_{\bm{\theta}}(\lambda_{j^\prime};u_{i_1^\prime},v_{i_2^\prime})},\\
L_{{4_{st}}n}\coloneqq &\frac{8^2\pi^2 L^4}{n^2M_n^4}\sum_{j,j^\prime=1}^{n-1}\sum_{i_1,i_1^\prime,i_2,i_2^\prime=1}^{M_n}{\rm Cum}\left(\overline{I_n(\lambda_j;u_{i_1},v_{i_2})}
,\overline{I_n(\lambda_{j^\prime};u_{i_1^\prime},v_{i_2^\prime})}
\right)\\
&\times
\frac{\partial}{\partial {\theta_s}}f_{\bm{\theta}}(\lambda_j;u_{i_1},v_{i_2})
\frac{\partial}{\partial {\theta_t}}f_{\bm{\theta}}(\lambda_{j^\prime};u_{i_1^\prime},v_{i_2^\prime}).
\end{align*}
Here, we assume without loss of generality that in the order four spectrum $$f(\lambda,-\lambda,\lambda^\prime;u,v,u^\prime,v^\prime)$$ all the H\"{o}lder exponents associated to variables $(u,v,u^\prime,v^\prime)$ are $\tilde{\alpha}$ and that the H\"{o}lder exponent associated to $\epsilon$ is $\tilde{\tau}$. This fact has no importance as the exponents do not play any role in the statement of the final result. Lemma \ref{lem_cum} and Lemma P5.1 of \cite[p.415]{b81} { alongside the assumption that the fourth-order generalized spectrum is Partially H\"{o}lder  yield that }
\begin{align*}
L_{1n}=&
\frac{8^2\pi^2 L^4}{n^2M_n^4}\sum_{j=1}^{n-1}\sum_{i_1,i_1^\prime,i_2,i_2^\prime=1}^{M_n}
f(\lambda_j;u_{i_1},u_{i_1^\prime})
f(-\lambda_j;v_{i_2},v_{i_2^\prime})\\
&\qquad\qquad\qquad\qquad\qquad\times
\frac{\partial}{\partial {\theta_s}}
\overline{f_{\bm{\theta}}(\lambda_j;u_{i_1},v_{i_2})}\Big|_{\bm{\theta}={{\bm{\theta}}_0}}
\frac{\partial}{\partial {\theta_t}}\overline{f_{\bm{\theta}}(2\pi-\lambda_{j};u_{i_1^\prime},v_{i_2^\prime})}\Big|_{\bm{\theta}={{\bm{\theta}}_0}}\\
&+\frac{8^2\pi^2 L^4}{n^2M_n^4}\sum_{j=1}^{n-1}\sum_{i_1,i_1^\prime,i_2,i_2^\prime=1}^{M_n}
f(\lambda_j;u_{i_1},v_{i_2^\prime})
f(-\lambda_j;v_{i_2},u_{i_1^\prime})\\
&\qquad\qquad\qquad\qquad\qquad\times
\frac{\partial}{\partial {\theta_s}}
\overline{f_{\bm{\theta}}(\lambda_j;u_{i_1},v_{i_2})}\Big|_{\bm{\theta}={{\bm{\theta}}_0}}
\frac{\partial}{\partial {\theta_t}}\overline{f_{\bm{\theta}}(\lambda_{j};u_{i_1^\prime},v_{i_2^\prime})}\Big|_{\bm{\theta}={{\bm{\theta}}_0}}\\
&+\frac{8^2\pi^2 L^4}{n^2M_n^4}\sum_{j,j^\prime=1}^{n-1}\sum_{i_1,i_1^\prime,i_2,i_2^\prime=1}^{M_n}
\frac{2\pi}{n}
f(\lambda_j,-\lambda_j,\lambda_{j^\prime};u_{i_1},v_{j_2},u_{i_1^\prime},v_{j_2^\prime})\\
&\qquad\qquad\qquad\qquad\qquad\times
\frac{\partial}{\partial {\theta_s}}\overline{f_{\bm{\theta}}(\lambda_j;u_{i_1},v_{i_2})}\Big|_{\bm{\theta}={{\bm{\theta}}_0}}
\frac{\partial}{\partial {\theta_t}}\overline{f_{\bm{\theta}}(\lambda_{j^\prime};u_{i_1^\prime},v_{i_2^\prime})}\Big|_{\bm{\theta}={{\bm{\theta}}_0}}\\
&+O(1/n^2)\\
=&
\frac{2\pi}{n}
\int_{0}^{2\pi}\int_{[-L,L]^4}
f(\lambda;u,u^\prime)
f(-\lambda;v,v^\prime)
\frac{\partial}{\partial {\theta_s}}\overline{f_{\bm{\theta}}(\lambda;u,v)}
\Big|_{\bm{\theta}={{\bm{\theta}}_0}}\\
&\qquad\qquad\qquad\qquad\qquad\times
\frac{\partial}{\partial {\theta_t}}\overline{f_{\bm{\theta}}(2\pi-\lambda;u^\prime,v^\prime)}\Big|_{\bm{\theta}={{\bm{\theta}}_0}}
{\rm d}u{\rm d}u^\prime{\rm d}v{\rm d}v^\prime{\rm d}\lambda\\
&+\frac{2\pi}{n}
\int_{0}^{2\pi}
\int_{[-L,L]^4}
f(\lambda;u,v^\prime)
f(-\lambda;v,u^\prime)
\frac{\partial}{\partial {\theta_s}}\overline{f_{\bm{\theta}}(\lambda;u,v)}\Big|_{\bm{\theta}={{\bm{\theta}}_0}}\\
&\qquad\qquad\qquad\qquad\qquad\times
\frac{\partial}{\partial {\theta_t}}\overline{f_{\bm{\theta}}(\lambda;u^\prime,v^\prime)}\Big|_{\bm{\theta}={{\bm{\theta}}_0}}
{\rm d}u{\rm d}u^\prime{\rm d}v{\rm d}v^\prime{\rm d}\lambda\\
&+\frac{2\pi}{n}
\int_{[0,2\pi]^2}
\int_{[-L,L]^4}
f(\lambda,-\lambda,\lambda^\prime;u,v,u^\prime,v^\prime)\\
&\qquad\qquad\qquad\qquad\times
\frac{\partial}{\partial {\theta_s}}\overline{f_{\bm{\theta}}(\lambda;u,v)}\Big|_{\bm{\theta}={{\bm{\theta}}_0}}
\frac{\partial}{\partial {\theta_t}}\overline{f_{\bm{\theta}}(\lambda^\prime;u^\prime,v^\prime)}\Big|_{\bm{\theta}={{\bm{\theta}}_0}}
{\rm d}u{\rm d}u^\prime{\rm d}v{\rm d}v^\prime{\rm d}\lambda{\rm d}\lambda^\prime\\
&+{ O\skakko{n^{-1}\skakko{n^{-\min \{\tilde{\alpha}, \alpha\}}+M_n^{-\min\{\tilde{\beta},\beta\}}}^{{1}/{(\max\{\tilde{\tau},\tau\}+1)}}
+(nM_n)^{-1}+n^{-2}}}
\end{align*}
and thus $n L_{1n}$ converges to $L_{1}$ as $n\to\infty$.
The convergence of $n L_{in}$ to $L_{i}$ for $i=2,3,4$ can be shown in the same way above.

For any integer $\ell\geq3$, the $\ell$-th order cumulant of ${\partial}/{\partial \bm{\theta}}D_n(I_n,f_{\bm{\theta}})\big|_{\bm{\theta}={{\bm{\theta}}_0}}$ can be evaluated as follows: For notational simplicity, define $W_{k1}\coloneqq I_n(\lambda_{j_k};u_{i_k},v_{i_k^\prime})$ and 
$W_{k2}\coloneqq \overline{I_n(\lambda_{j_k};u_{i_k},v_{i_k^\prime})}$. Suppose that $(\nu_1,\ldots,\nu_p)$ is an indecomposable partition of the table
\begin{align*}
\begin{pmatrix}
d((-1)^{t_1-1}\lambda_{j_1};(-1)^{t_1-1}u_{i_1}) &d((-1)^{t_1}\lambda_{j_1};(-1)^{t_1-1}v_{i_1^\prime})\\
\vdots&\vdots\\
d((-1)^{t_\ell-1}\lambda_{j_\ell};(-1)^{t_\ell-1}u_{i_\ell}) &d((-1)^{t_\ell}\lambda_{j_\ell};(-1)^{t_\ell-1}v_{i_\ell^\prime})\\
\end{pmatrix}
\end{align*}
and denote $\nu_j$ by $\{d(\omega_{1\nu_j};q_{1\nu_j}), d(\omega_{2\nu_j};q_{2\nu_j}),\ldots,d(\omega_{\#\nu_j\nu_j};q_{\#\nu_j\nu_j})\}$. For any $s_1,\ldots,s_\ell\in\{1,\ldots,d\}$, it can be seen that
\begin{align*}
&{\rm Cum}\skakko{
\frac{\partial}{\partial {\theta_{s_1}}}D_n(I_n,f_{\bm{\theta}})\Big|_{\bm{\theta}={{\bm{\theta}}_0}},
\frac{\partial}{\partial {\theta_{s_2}}}D_n(I_n,f_{\bm{\theta}})\Big|_{\bm{\theta}={{\bm{\theta}}_0}},
\cdots,
\frac{\partial}{\partial {\theta_{s_\ell}}}D_n(I_n,f_{\bm{\theta}})\Big|_{\bm{\theta}={{\bm{\theta}}_0}}}\\
=&\skakko{\frac{8\pi L^2}{nM_n^2}}^\ell\sum_{j_1,\ldots,j_\ell=1}^{n-1}\sum_{i_1\ldots,i_\ell,i_1^\prime,\ldots,i_\ell^\prime=1}^{M_n}\\
&\quad{\rm Cum}\bigg(
\frac{\partial}{\partial {\theta_{s_1}}}| I_n(\lambda_{j_1};u_{i_1},v_{i_1^\prime})-
f_{\bm{\theta}}(\lambda_{j_1};u_{i_1},v_{i_1^\prime})|^2\Big|_{\bm{\theta}={{\bm{\theta}}_0}}\\
&\qquad\qquad\qquad\qquad\qquad\qquad\qquad
,\cdots,
\frac{\partial}{\partial {\theta_{s_\ell}}}| I_n(\lambda_{j_\ell};u_{i_\ell},v_{i_\ell^\prime})-
f_{\bm{\theta}}(\lambda_{j_\ell};u_{i_\ell},v_{i_\ell^\prime})|^2\Big|_{\bm{\theta}={{\bm{\theta}}_0}}
\bigg)\\
=&
\skakko{\frac{8\pi L^2}{nM_n^2}}^\ell\sum_{j_1,\ldots,j_\ell=1}^{n-1}\sum_{i_1\ldots,i_\ell,i_1^\prime,\ldots,i_\ell^\prime=1}^{M_n}
\sum_{t_1,\ldots,t_\ell=1,2}
\left|{\rm Cum}\left(W_{1t_1},\cdots,W_{\ell t_\ell}\right)\right|\\
&\times
\prod_{k=1}^\ell\left|\frac{\partial}{\partial {\theta_{s_k}}}f_{\bm{\theta}}(\lambda_{j_k};u_{i_k},v_{i_k^\prime})\Big|_{\bm{\theta}={{\bm{\theta}}_0}}\right|
\\
=&
\skakko{\frac{8\pi L^2}{nM_n^2}}^\ell\sum_{j_1,\ldots,j_\ell=1}^{n-1}\sum_{i_1\ldots,i_\ell,i_1^\prime,\ldots,i_\ell^\prime=1}^{M_n}
\sum_{t_1,\ldots,t_\ell=1,2}\frac{1}{(2\pi n)^\ell}\\
&\times
\Bigg|\sum_{(\nu_1,\ldots,\nu_p)}
\Bigg((2\pi)^{\# \nu_1-1}\Delta_n\skakko{\sum_{k=1}^{\# \nu_1}\omega_{k\nu_1}}
f(\omega_{1\nu_1},\ldots,\omega_{(\# \nu_1-1)\nu_1};
q_{1\nu_1},\ldots,q_{\#\nu_1\nu_1})\\
&\qquad\qquad\qquad\qquad\qquad\qquad\qquad\qquad\qquad\qquad\qquad\qquad\qquad\qquad\qquad
+O(1)\Bigg)
\cdots\\
&\times
\skakko{(2\pi)^{\# \nu_p-1}\Delta_n\skakko{\sum_{k=1}^{\# \nu_p}\omega_{k\nu_p}}
f(\omega_{1\nu_p},\ldots,\omega_{(\# \nu_p-1)\nu_p};
q_{1\nu_p},\ldots,q_{\#\nu_p\nu_p})+O(1)}\Bigg|\\
&\times
\prod_{k=1}^\ell\left|\frac{\partial}{\partial {\theta_{s_k}}}f_{\bm{\theta}}(\lambda_{j_k};u_{i_k},v_{i_k^\prime})\Big|_{\bm{\theta}={{\bm{\theta}}_0}}\right|
,
\end{align*}
where the summation $\sum_{(\nu_1,\ldots,\nu_p)}$ extends over all indecomposable partitions $(\nu_1,\ldots,\nu_p)$  of the { aforementioned} table.
There exists, for any $k\in{0,\ldots,\ell-1}$, the indecomposable partition with the number of sets is $\ell-k$ but $\ell-1-k$ constraints for frequencies are required so that all $\Delta_n\skakko{\sum_{k=1}^{\# \nu_1}\omega_{k\nu_1}},\ldots,\Delta_n\skakko{\sum_{k=1}^{\# \nu_{\ell-k}}\omega_{k\nu_{\ell-k}}}$ equal to $n$.
However, the number of restrictions cannot be improved by construction. Therefore, the $\ell$-th order cumulant of ${\partial}/{\partial \bm{\theta}}D_n(I_n,f_{\bm{\theta}})\big|_{\bm{\theta}={{\bm{\theta}}_0}}$ is of order $n^{-\ell+1}$ and we conclude that 
\begin{align*}
&{\rm Cum}\skakko{
\sqrt n\frac{\partial}{\partial {\theta_{s_1}}}D_n(I_n,f_{\bm{\theta}})\Big|_{\bm{\theta}={{\bm{\theta}}_0}},
\sqrt n\frac{\partial}{\partial {\theta_{s_2}}}D_n(I_n,f_{\bm{\theta}})\Big|_{\bm{\theta}={{\bm{\theta}}_0}},
\cdots,
\sqrt n\frac{\partial}{\partial {\theta_{s_\ell}}}D_n(I_n,f_{\bm{\theta}})\Big|_{\bm{\theta}={{\bm{\theta}}_0}}}\\
&=O(n^{-\ell/2+1}).
\end{align*}
\qed

\subsection{Proof of Theorem \ref{thm3}}
The proof is omitted since it is {(a simpler version of) the proof of Theorem \ref{thm4}}.\qed

\subsection{Proof of Theorem \ref{thm4}}
We follow the idea of the proofs of Theorem 4.1 of  \cite{gkvvdh22} and Theorem 3 of \cite{dw16}.

First, we show our test has level ${\varphi}$.
Define $\delta_i\coloneqq {\theta}_{0i}$. 
Under the null hypothesis, we have $\skakko{\sqrt n-\sqrt b}(\delta_i-\kappa_i)\leq0$ and thus
\begin{align*}
p_{n,3}=&\frac{1}{n-b+1}\sum_{t=1}^{n-b+1}\mathbb I
\left\{
\sqrt b \skakko{\hat{\theta}_{b,t}^{(i)}-\delta_i}
>
\sqrt n \skakko{\hat{\theta}_{n}^{(i)}-\delta_i}+\skakko{\sqrt n-\sqrt b}(\delta_i-\kappa_i)
\right\}\\
\geq&
\frac{1}{n-b+1}\sum_{t=1}^{n-b+1}\mathbb I
\left\{
\sqrt b \skakko{\hat{\theta}_{b,t}^{(i)}-\delta_i}
>
\sqrt n \skakko{\hat{\theta}_{n}^{(i)}-\delta_i}
\right\}\\
=&1-H_{n,b}\skakko{\sqrt n \skakko{\hat{\theta}_{n}^{(i)}-\delta_i}},
\end{align*}
where
\begin{align*}
H_{n,b}\skakko{x}\coloneqq 
\frac{1}{n-b+1}\sum_{t=1}^{n-b+1}\mathbb I
\left\{
\sqrt b \skakko{\hat{\theta}_{b,t}^{(i)}-\delta_i}
\leq
x
\right\}.
\end{align*}
Define $Y_n\coloneqq \sqrt n \skakko{\hat{\theta}_{n}^{(i)}-\delta_i}$ and $Y\coloneqq \lim_{n\to\infty}Y_n$ and
Let $R_n$ and $R$ denote the distribution function of $Y_n$ and $Y$, respectively.
From Proposition 7.3.1 of \cite{prw99}, $\rho_L(H_{n,b},R_n)$ converges in probability to zero as $n\to\infty$,
where $\rho_L$ is the bounded Lipschitz metric on the space of distribution function on $\mathbb R$. 
Theorem \ref{thm2} gives that $\rho_L(R_n,R)$  converges in probability to zero as $n\to\infty$.
Hence, $\rho_L(H_{n,b},R)$ converges in probability to zero as $n\to\infty$.
In the same manner as the proof of Theorem 4.1 of  \cite{gkvvdh22}, we have
$\sup_{x\in\mathbb R}\left|H_{n,b}(x)-R(x)\right|=o_p(1)$. 
Then, it holds that 
\begin{align*}
{\rm P}\skakko{p_{n,3}\leq {\varphi}}
\leq& {\rm P}\skakko{1-H_{n,b}\skakko{\sqrt n \skakko{\hat{\theta}_{n}^{(i)}-\delta_i}}\leq {\varphi}}\\
\leq& {\rm P}\skakko{1-{\varphi} \leq H_{n,b}\skakko{\sqrt n \skakko{\hat{\theta}_{n}^{(i)}-\delta_i}}}\\
=& {\rm P}\skakko{1-{\varphi} \leq R\skakko{Y_n}+o_p(1)},
\end{align*}
which, by the continuous mapping theorem, converges to ${\rm P}\skakko{1-{\varphi} \leq R\skakko{Y}}={\varphi}$.
Thus, our test has { asymptotic} level ${\varphi}$, that is, $\lim_{n \rightarrow \infty}{\rm P}\skakko{p_{n,3}\leq {\varphi}}\leq{\varphi}$.

Next, we show the consistency, that is, under the alternative hypothesis
${\rm P}\skakko{p_{n,3}\leq {\varphi}}\to1$ as $n\to\infty$, or equivalently, ${\rm P}\skakko{p_{n,3}>{\varphi}}\to0$  as $n\to\infty$. Since $Y_n=O_p(1)$, for any $\epsilon>0$ there exists $M_\epsilon$ such that ${\rm P}\skakko{|Y_n|\geq M_\epsilon}<\epsilon$ for all $n$. Then, it holds
\begin{align*}
{\rm P}\skakko{p_{n,3}> {\varphi}}
=&
{\rm P}\skakko{1-H_{n,b}\skakko{\sqrt n \skakko{\hat{\theta}_{n}^{(i)}-\delta_i}+\skakko{\sqrt n-\sqrt b}(\delta_i-\kappa_i)}> {\varphi}}\\
\leq&
{\rm P}\skakko{1-H_{n,b}\skakko{-M_\epsilon+\skakko{\sqrt n-\sqrt b}(\delta_i-\kappa_i)}> {\varphi}}
+
{\rm P}\skakko{|Y_n|\geq M_\epsilon}\\
\leq&
{\rm P}\skakko{1-{ R }\skakko{-M_\epsilon+\skakko{\sqrt n-\sqrt b}(\delta_i-\kappa_i)}+o_p(1)> {\varphi}}
+
\epsilon,
\end{align*}
{where we used $\sup_{x\in\mathbb R}\left|H_{n,b}(x)-R(x)\right|=o_p(1)$ again on the last line.} Now, since $\delta_i-\kappa_i>0$ under the alternative and $b/n\to0$ as $n\to\infty$, { the previous display is arbitrarily close to} $\epsilon$ for large $n$. \qed

\section{Additional plots of generalized spectra}\label{sec11}

This section presents heatmaps of the generalized spectra for the integer-valued models discussed in Section 3. Figures \ref{Fig_spec_all_real_MA} and \ref{Fig_spec_all_imaginary_MA} illustrate the real and imaginary parts of the spectrum for the INMA(1) model (Example 5). Similarly, Figures \ref{Fig_spec_all_real_AR} and \ref{Fig_spec_all_imaginary_AR} depict the spectral components for the INAR(1) model (Example 6).

\begin{figure}[htbp]
\centering
\includegraphics[width=0.8\textwidth,page=1]{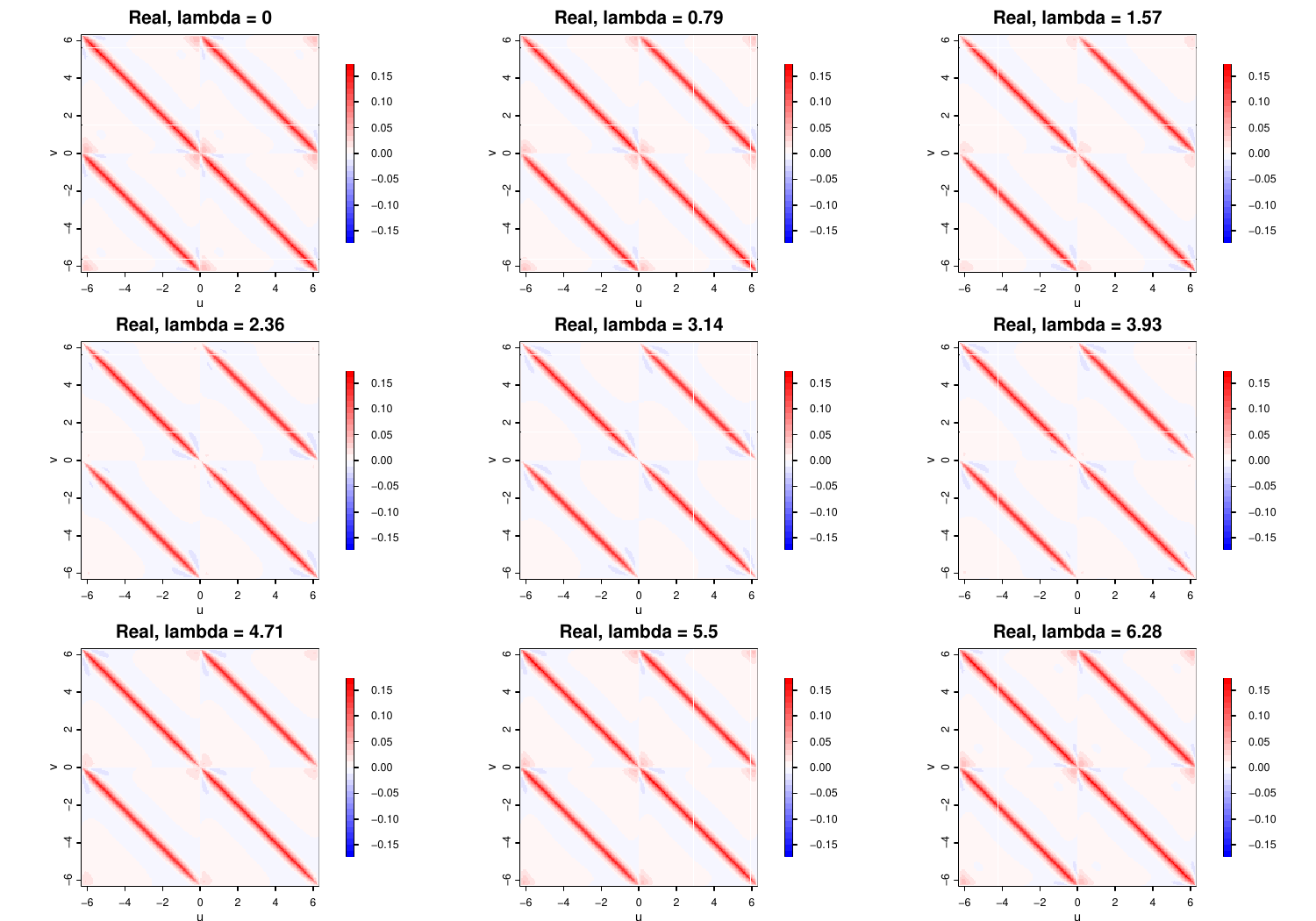}
  \caption{
  Heatmaps of the real part of the spectrum for the causal INMA(1) model in Example~5, with parameters $\delta = 2$, $p = 0.3$, $\alpha = 0.7$, and $\lambda$ values equally spaced over $[0, 2\pi]$ with nine points.}
  \label{Fig_spec_all_real_MA}
\end{figure}

\begin{figure}[htbp]
\centering
\includegraphics[width=0.8\textwidth,page=2]{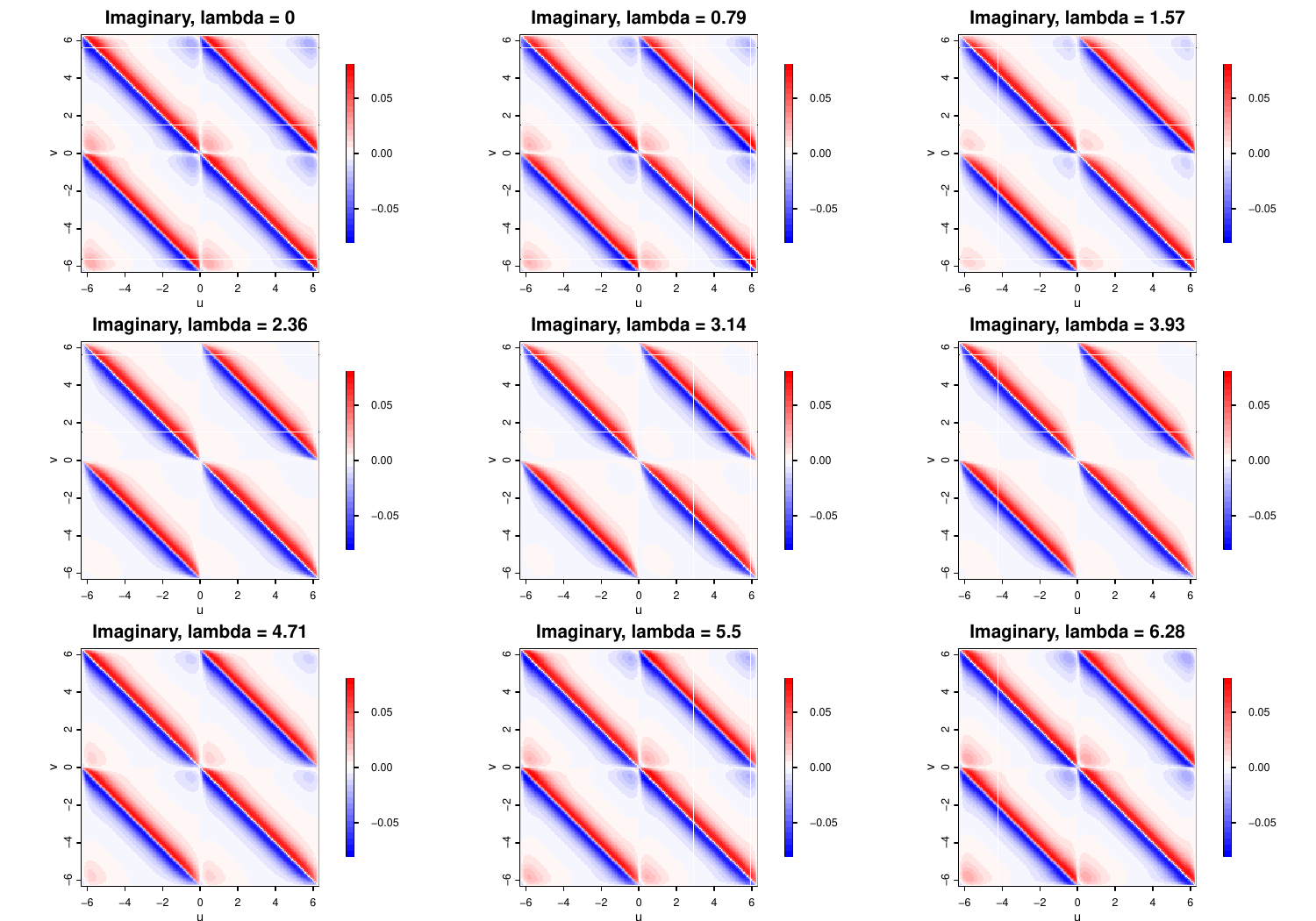}
  \caption{
  Heatmaps of the imaginary part of the spectrum for the causal INMA(1) model in Example~5, with parameters $\delta = 2$, $p = 0.3$, $\alpha = 0.7$, and $\lambda$ values equally spaced over $[0, 2\pi]$ with nine points.}
  \label{Fig_spec_all_imaginary_MA}
\end{figure}

\begin{figure}[htbp]
\centering
\includegraphics[width=0.8\textwidth,page=1]{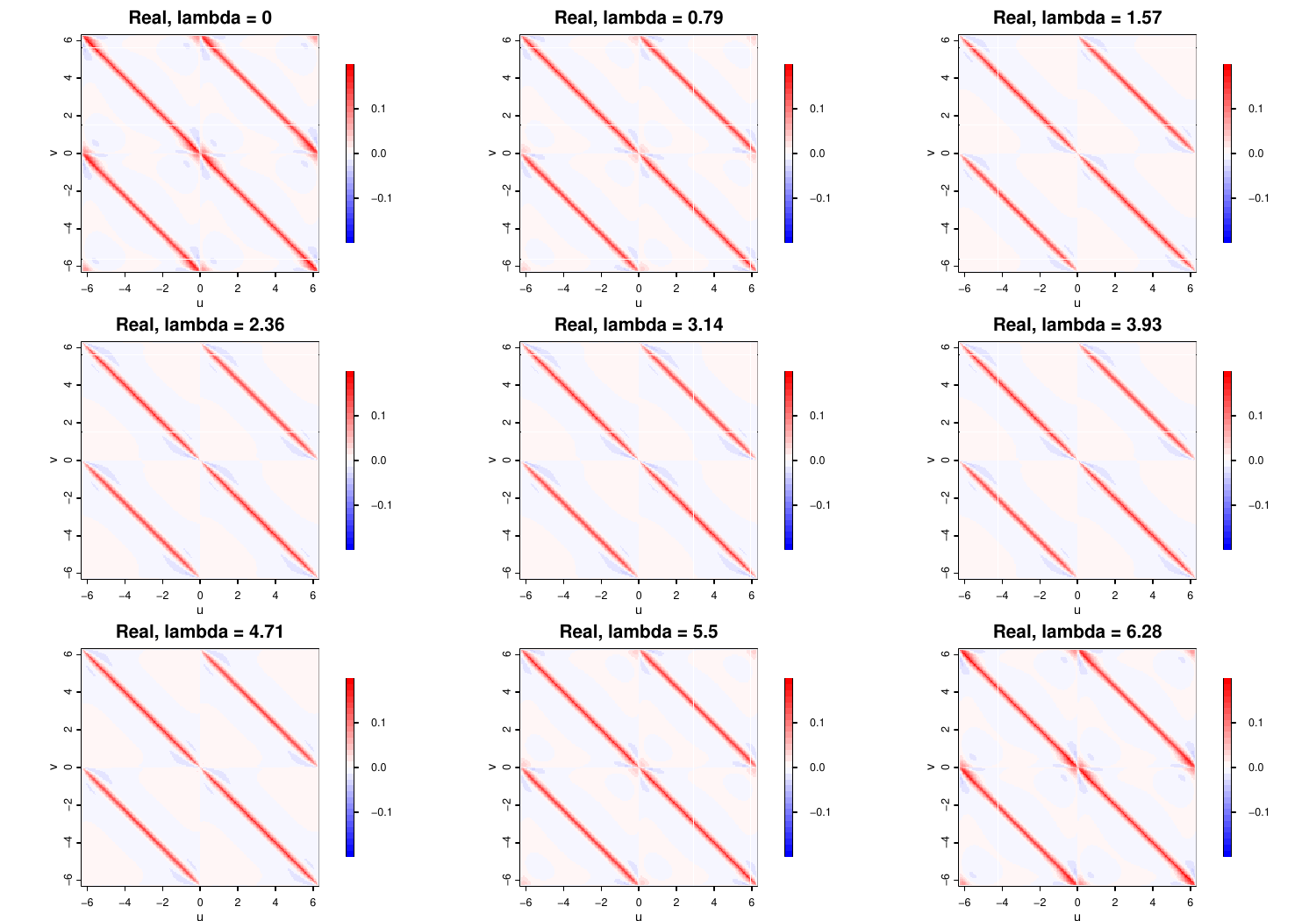}
  \caption{
  Heatmaps of the real part of the spectrum for the causal integer-valued AR(1) model in Example~6, with parameters $\delta = 2$, $p = 0.3$, $\alpha = 0.7$, and $\lambda$ values equally spaced over $[0, 2\pi]$ with nine points.}
  \label{Fig_spec_all_real_AR}
\end{figure}

\begin{figure}[htbp]
\centering
\includegraphics[width=0.8\textwidth,page=2]{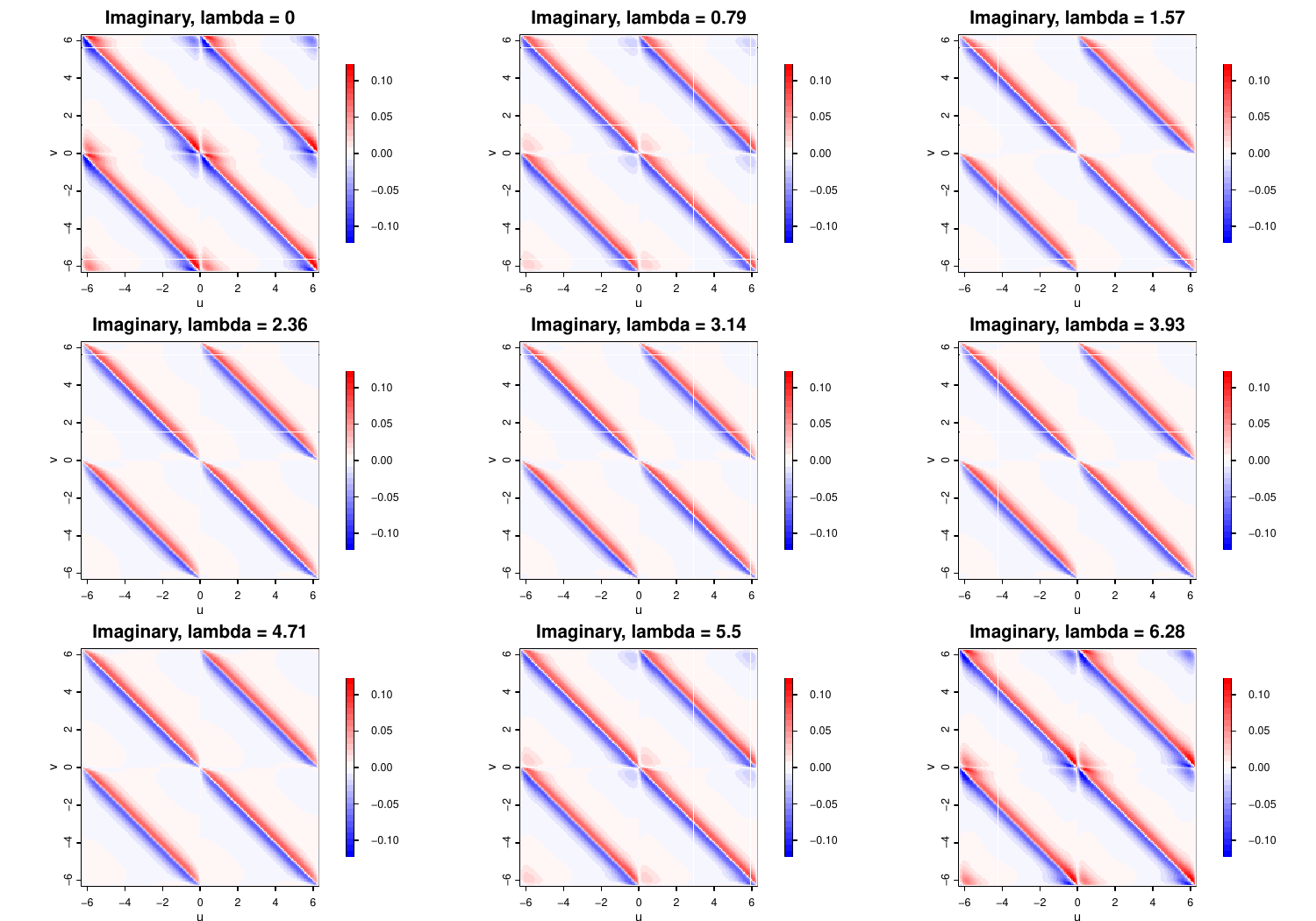}
  \caption{
  Heatmaps of the imaginary part of the spectrum for the causal integer-valued AR(1) model in Example~6, with parameters $\delta = 2$, $p = 0.3$, $\alpha = 0.7$, and $\lambda$ values equally spaced over $[0, 2\pi]$ with nine points.}
  \label{Fig_spec_all_imaginary_AR}
\end{figure}

\clearpage

\section{Additional simulation results}\label{sec12}

Figures \ref{Fig1} and \ref{Fig2bis} present the simulation results for the MA(1) and causal AR(1) models discussed in Section 7.1.1 of the main text. Figure \ref{Fig4} provides the simulation results for the INMA(1) model described in Section 7.1.2.

\begin{figure}[!htbp]
\vspace{-0mm}
\centering
\includegraphics[width=0.8\textwidth]{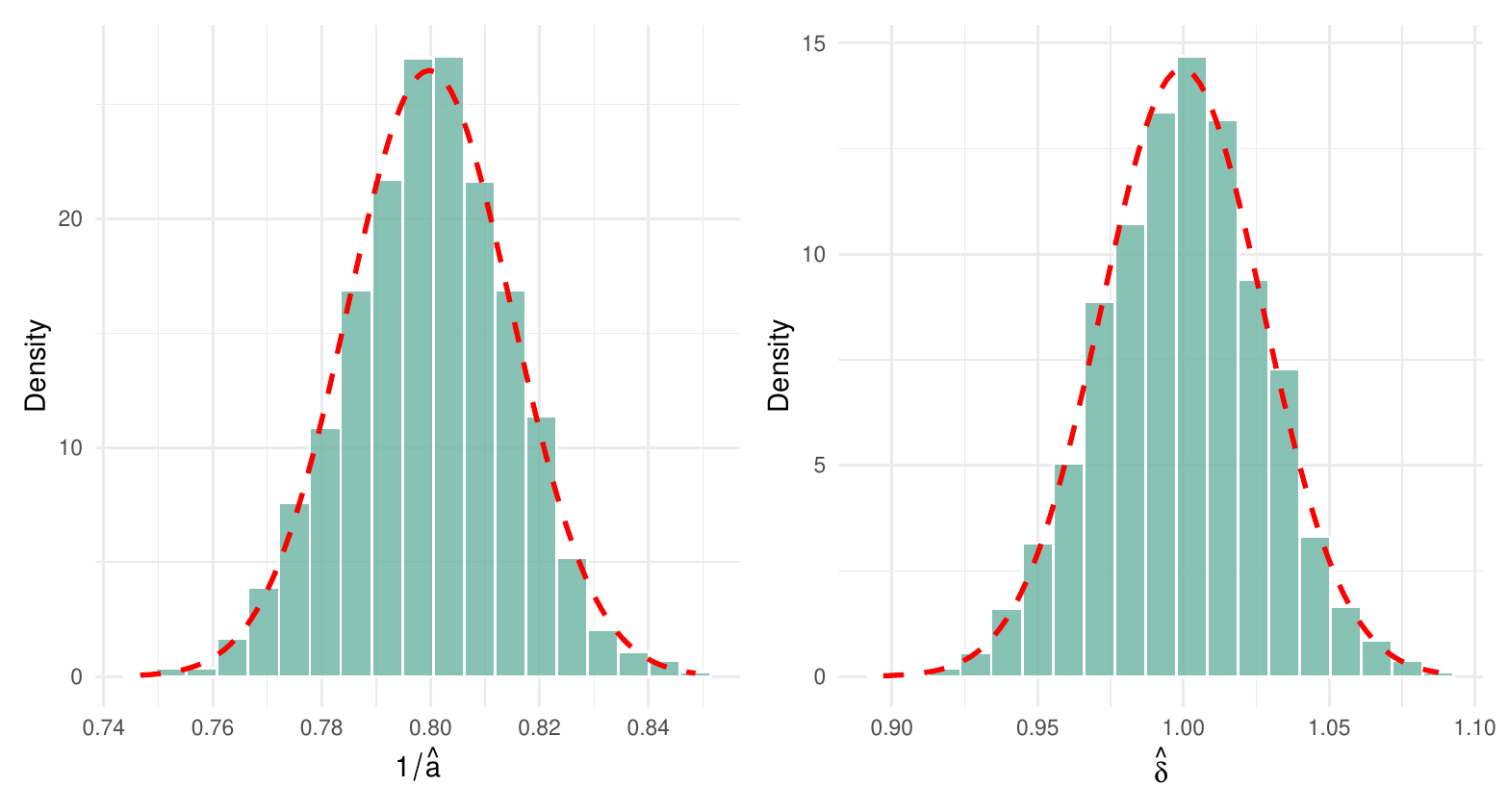}
\vspace*{-0mm}
  \caption{Histograms of ${\hat{a}_n^{-1}}$ and $\hat{\delta}_n$ (in green) with respect to limiting Gaussian distribution (in red) in the MA$(1)$ model with Cauchy innovations. $2000$ replications, $n=500$ and $(a, \delta)=({0.8^{-1}}, 1)$.}
  \label{Fig1}
\end{figure}

\begin{figure}[!htbp]
\vspace{-0mm}
\centering
\includegraphics[width=0.8\textwidth]{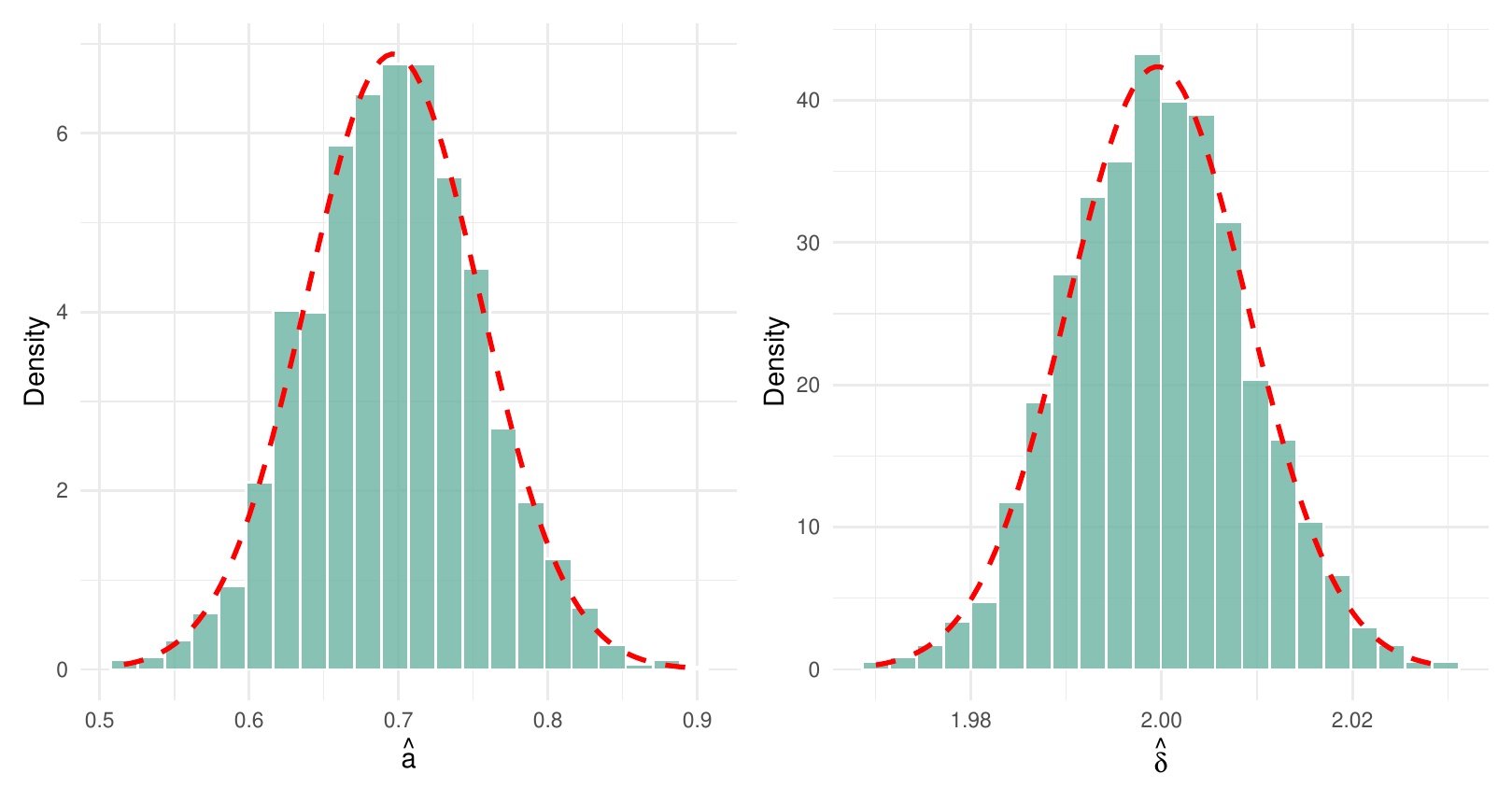}
\vspace*{-0mm}
  \caption{Histograms of ${\hat{a}_n}$ and $\hat{\delta}_n$ (in green) with respect to limiting Gaussian distribution (in red) in the causal AR$(1)$ model with Cauchy innovations. $2000$ replications, $n=500$ and $(a, \delta)=(0.7, 2)$.}
  \label{Fig2bis}
\end{figure}

\begin{figure}[!htbp]
\vspace{-0mm}
\centering
\includegraphics[width=0.8\textwidth]{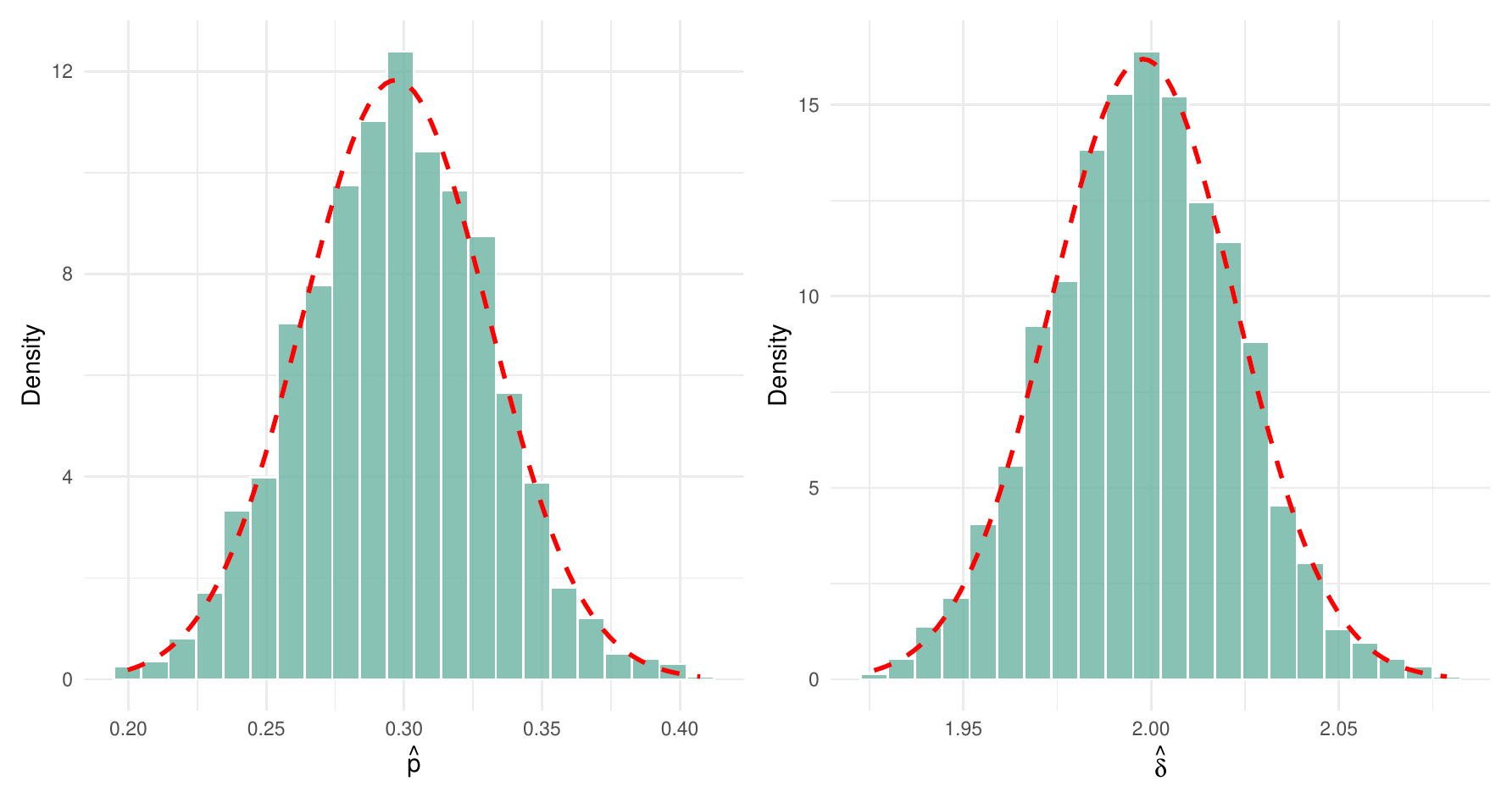}
\vspace*{-0mm}
  \caption{Histograms of ${\hat{p}_n}$ and $\hat{\delta}_n$ (in green) with respect to limiting Gaussian distribution (in red) in the INMA$(1)$ model with discrete stable innovations. $2000$ replication, $n=500$ and $(p, \delta)=(0.3, 2)$. The rate of selection of the proper $\alpha=0.7$ is a hundred percent.}
  \label{Fig4}
\end{figure}


\end{document}